\documentclass[11pt]{amsart}
\usepackage{amsfonts}
\usepackage{amsmath}
\usepackage{amssymb, esint}
\usepackage{amsthm,color}
\usepackage{enumerate}
\usepackage{mathtools}
\usepackage[hidelinks, bookmarksdepth=3]{hyperref}

\usepackage{cancel}
\usepackage{longtable}

\usepackage{geometry}
\newcommand{\R}{{\mathbb R}}

\newcommand{\N}{{\mathbb N}}

\newcommand{\cM}{{\mathcal M}}
\newcommand{\cN}{{\mathcal N}}

\newcommand{\cE}{{\mathcal E}}

\newcommand{\cH}{{\mathcal H}}

\newcommand{\cR}{{\mathcal R}}

\def\0{{\mathbf 0}}

\def\loc{{\textup{loc}}}

\newcommand{\id}{\mathrm{id}}
\newcommand{\e}{\varepsilon}
\newcommand{\vp}{\varphi}
\newcommand{\osc}{\operatornamewithlimits{osc}}

\newcommand{\supp}{\operatorname{supp}}
\newcommand{\ddiv}{\operatorname{div}}
\newcommand{\diam}{\operatorname{diam}}
\newcommand{\dist}{\operatorname{dist}}

\newcommand{\Hess}{\operatorname{Hess}}

\newcommand{\means}[1]{-\hskip-1.00em\int_{#1}}

\theoremstyle{plain}
\newtheorem{thm}{Theorem}[section]

\newtheorem{cor}[thm]{Corollary}
\newtheorem{lem}[thm]{Lemma}
\newtheorem{prop}[thm]{Proposition}

\theoremstyle{definition}
\newtheorem{defn}[thm]{Definition}

\theoremstyle{remark}

\newtheorem{case}{Case}

\newtheorem*{claim*}{Claim}

\newtheorem{rem}[thm]{Remark}

\numberwithin{equation}{section}

\title[Constraint maps, Bernoulli type free boundary]{Constraint maps with free boundaries:\\ the Bernoulli case}

\author[A.\ Figalli]{Alessio Figalli}
\email{alessio.figalli@math.ethz.ch}
\address{Department of Mathematics, ETH Z\"urich,  Raemistrasse 101, 8092 Z\"urich, Switzerland }

\author[A.\ Guerra]{Andr\'e Guerra}
\email{andre.guerra@eth-its.ethz.ch}
\address{Institute for Theoretical Studies, ETH Z\"urich,  CLV, Clausiusstrasse 47, 8006 Z\"urich, Switzerland}

\author[S.\ Kim]{Sunghan Kim}
\email{sunghan.kim@math.uu.se}
\address{Department of Mathematics, Uppsala University, S-751 06 Uppsala, Sweden}

\author[H.\ Shahgholian]{Henrik Shahgholian}
\email{henriksh@kth.se}
\address{Department of Mathematics, KTH Royal Institute of Technology, 100 44 Stockholm, Sweden}

\begin{document}

\begin{abstract}
In this manuscript, we delve into the study of maps $u\in W^{1,2}(\Omega;\overline M)$
that minimize the Alt--Caffarelli energy functional
$$
\int_\Omega (|Du|^2 + q^2 \chi_{u^{-1}(M)})\,dx, 
$$
under the condition that the image $u(\Omega)$ is confined within $\overline M$. Here, $\Omega$ denotes a bounded domain in the ambient space $\R^n$ (with $n\geq 1$), and $M$ represents a smooth domain in the target space $\R^m$ (where $m\geq 2$).

Since our minimizing constraint maps coincide with harmonic maps in the interior of the coincidence set, ${\rm int}(u^{-1}(\partial M))$, such maps are prone to developing discontinuities due to their inherent nature. This research marks the commencement of an in-depth analysis of potential singularities that might arise within and around the free boundary.

Our first significant contribution is the validity of a $\varepsilon$-regularity theorem. This theorem is founded on a novel method of Lipschitz approximation near points exhibiting low energy. Utilizing this approximation and extending the analysis through a bootstrapping approach, we show Lipschitz continuity of our maps whenever the energy is small.

Our subsequent key finding reveals that,
whenever the complement of $M$ is uniformly convex and of class $C^3$, the maps minimizing the Alt--Caffarelli energy with a positive parameter $q$ exhibit Lipschitz continuity within a universally defined neighborhood of the non-coincidence set $u^{-1}(M)$. In particular, this Lipschitz continuity extends to the free boundary.

A noteworthy consequence of our findings is the smoothness of flat free boundaries and of the resulting image maps.

\end{abstract}

\maketitle
\setcounter{tocdepth}{1}
\tableofcontents


\section{Introduction}\label{sec:intro}

Let $\Omega$ be a bounded domain in $\R^n$, $n\geq 1$, and consider a smooth domain $M$ in $\mathbb{R}^m$, $m\geq 2$, whose boundary is of class $C^3$. 
We consider maps $u\in W^{1,2}(\Omega;\overline{M})$ which are minimizers of the energy functional
\begin{equation}\label{eq:main}
\cE_q[v] := \int_\Omega (|Dv|^2 + q^2 \chi_{v^{-1}(M)}) \, dx,
\end{equation} 
where $q\geq 0$ is a fixed constant. We refer to such maps as \textit{minimizing constraint maps}: precisely, this means that $u\in W^{1,2}(\Omega;\overline M)$ is such that
$$\cE_q[u]\leq \cE_q[v] \quad \text{for all }  v\in W^{1,2}(\Omega;\overline M)\cap (u+W^{1,2}_0(\Omega;\R^m)).$$

It should be noted that constraint maps, when restricted to the interior of the coincidence set ${\rm int}(u^{-1}(\partial M))$, are harmonic maps into $\partial M$. Thus, similarly to how harmonic maps develop singularities, so do constraint maps, and their singularities are denoted   throughout the paper by 
\begin{equation}\label{eq:Sing}
\Sigma(u) := \{ x\in\Omega: u\text{ is not continuous at }x\}.
\end{equation}
At the same time, each component of a constraint map is harmonic in the interior of the non-coincidence set ${\rm int}(u^{-1}(M))$, and so constraint maps give rise to a free boundary $\Omega\cap \partial u^{-1}(M)$. For $q>0$, this free boundary problem can be seen as a natural extension of the scalar Bernoulli problem into the vectorial setting. 

This paper is motivated by the following question: {\it Do singularities appear on or around the free boundary?} Note that free boundaries provide further information, as they impose additional conditions on the solutions. Even though it is a highly delicate issue to classify or locate the singularities of a harmonic map,  in this paper we show that when the functional $\cE_q$ has positive weight $q$ and the target constraint $M$ have uniformly convex complement (e.g., the exterior of a ball), the free boundary pushes away the singularities by a uniform distance to the interior of the coincidence set $u^{-1}(\partial M)$. 
Moreover, there is an {\it a priori} interior estimate of the Lipschitz norm of the mapping in this neighborhood. 

Along the way, we realized that the partial regularity theory itself (without any geometric condition on the target constraints)  deserves attention on its own, as it turns out to be highly nontrivial. Here, we establish a sharp partial regularity theorem, which is optimal in both dimensions of the singular set and the degree of regularity of the mapping in the complement of the singular set. 

Our results can be lifted to the manifold setting, i.e., to the case where $\overline M$ is a Riemannian manifold with boundary and of nonpositive sectional curvature, and $\Omega$ is equipped with a Riemannian metric. Moreover, our analysis allows for sufficiently regular (such as Lipschitz) but non-constant weights $q$. However, for the sake of clarity, we chose not to include such extensions here. 

In the forthcoming paper \cite{FGKS}, we will study the same question for minimizing constraint maps of the Dirichlet energy, $\cE_0$, where the maps can behave more wildly near their free boundaries. 

\subsection{A short overview of the literature}

The problem studied in this paper can be investigated from two different perspectives: the one coming from the theory of harmonic maps, and the one coming from the theory of free boundary problems of Bernoulli-type. 

Let us begin by briefly overviewing the theory of minimizing constraint maps for $\cE_0$, since they also play an important role in the study of constraint maps for $\cE_q$ when $q>0$. Minimizing constraint maps for $\cE_0$ were studied many decades ago, by F. Duzaar, M. Fuchs and many others. To mention a few results among a long list of literature regarding such variational problems, Duzaar showed in \cite{D} 
(in the manifold setting) 
that $\cH^{n-2}(\Sigma(u)) = 0$, where $\cH^k$ is the $k$-dimensional Hausdorff measure. Later in \cite{DF}, parallel to the development of the theory for minimizing harmonic maps, Duzaar and Fuchs established the optimal bound $\dim_\cH(\Sigma(u))\leq n-3$. Several interesting results were established, including (but by no means limited to) the free boundary regularity \cite{F} and the absence of singularities under suitable geometric or topological assumptions on the target \cite{F2, FW}. The theory was also  generalized to energy functionals with $p$-growth by Fuchs \cite{F3, F4}, and then to almost minimizing maps by Luckhaus \cite{L}. The constraint maps can also be understood as harmonic maps {\it into} manifolds-with-boundary; see \cite{CM} for further development beyond the energy-minimizing maps. 

Another fruitful perspective on our problem stems from the theory of free boundary problems. When $q>0$,  the energy functional in \eqref{eq:main} can be seen as a vectorial version of the Alt--Caffarelli functional \cite{AC}, and leads to a Bernoulli-type free boundary problem. Recently, the vectorial extension of the Bernoulli problem has received a lot of interest, see e.g.\ \cite{CSY,DPESV,KL1,KL2,MTV1,MTV2}. In our terminology, these papers are concerned with  minimizing constraint maps for $\cE_q$, $q> 0$, when $M = \R^m\setminus \{0\}$, i.e.,  when $\partial M = \{0\}$. These references  also contain the case $M = \{ x \in \R^m   : x_i > 0\}$.
In these cases, the norm of the map, $|u|$, plays an essential role, and shares important characteristics with the minimizers for the scalar problem. 

Nevertheless, this is the first paper, to the best of our knowledge, that treats the case where $\partial M$ is non-flat and $q > 0$.  In fact, the case $q=0$ was only recently revisited by three of the current authors in \cite{FKS} in the context of vectorial obstacle problems, with joint perspectives from both harmonic maps and free boundary problems. The latter work focuses, however, on analysis of maps away from the singularities.


\subsection{Main results}\label{sec:main}

We start with an $\e$-regularity theorem, which states that minimizing constraint maps are Lipschitz in the vicinity of any point with small energy.

\begin{thm}[$\e$-regularity]\label{thm:e-reg}
Let $M$ be a $C^3$-domain in $\R^m$ with bounded principal curvatures  and let $u\in W^{1,2}(\Omega;\overline M)$ be a minimizing constraint map for the functional $\cE_q$, with $q\geq 0$. There are constants $\bar \e$, $\bar r,$ and $c$, which depend at most on $n$, $m$, $q$ and the $C^3$-character of $\partial M$, such that for every ball $B_{4r}(x_0)\Subset\Omega$ with $r \in(0,\bar r)$, the following holds: if
$$
r^{2-n} \int_{B_{4r}(x_0)} |Du|^2\,dx \leq \bar\e^2
$$ then $u \in C^{0,1}(B_r(x_0))$ with
$$
[u]_{C^{0,1}(B_r(x_0))} \leq \frac{c}{r}.
$$
\end{thm}

Theorem \ref{thm:e-reg} is optimal\footnote{Our methods rely on the smoothness of $\partial M$ up to $C^3$ regularity, but it is an interesting question whether one can still work with less regular manifolds. In a nutshell, the $C^3$-assumption is needed to obtain $C^2$-regularity of the nearest projection map $\Pi:\overline M\to \partial M$ which plays an important role in our analysis (see Corollary \ref{cor:proj-pde2}). 
} in the sense that $u$ may be Lipschitz but not $C^1$, see  Remark \ref{rem:e-reg}. To the best of our knowledge, no H\"older estimates, or indeed any type of continuity estimate, were previously established in the context of Theorem \ref{thm:e-reg}. Notably, a significant disparity exists in the regularity theory between minimizers of $\cE_0$ and $\cE_q$ for $q>0$. While the Dirichlet energy $\cE_0$ has been thoroughly examined several decades ago in \cite{D}, the Alt--Caffarelli energy $\cE_q$ for $q > 0$ entails the introduction of a characteristic function that subsequently gives rise to a  measure within the Euler--Lagrange system, as exemplified by \eqref{eq:main-sys} below. 
 Consequently, we are compelled to rely directly on the minimality assumption and cannot exclusively depend on PDE methods.

Through a suitable almost monotonicity formula and the usual dimension-reduction argument, Theorem \ref{thm:e-reg} yields the following optimal partial regularity result for the vectorial Alt--Caffarelli energy:

\begin{cor}\label{cor:e-reg3}
Let $M$ be a $C^3$-domain in $\R^m$ with bounded principal curvatures, and let $u\in W^{1,2}(\Omega;\overline M)$ be a minimizing constraint map for $\cE_q$, with $q > 0$. Then $\dim_\cH(\Sigma(u)) \leq n-3$ and $u\in C_\textup{loc}^{0,1}(\Omega\setminus\Sigma(u))$. 
\end{cor}

\begin{rem}\label{rem:e-reg}
Let us discuss the sharpness of Corollary \ref{cor:e-reg3}. The statement is optimal in the dimension of the singular set, see e.g. \cite{LW, Simon}, as minimizing harmonic maps are also minimizing constraint maps when the complement of $M$ is convex (since the projection map onto any closed convex set decreases the energy).

The sharpness on the regularity front follows the optimality of Lipschitz regularity for minimizers of the scalar Alt--Caffarelli energy \cite{AC}, as one can always lift the examples of the scalar minimizers to the vectorial setting by taking $M$ to be a half-space.
\end{rem}

\begin{rem}
As a consequence of Theorem \ref{thm:e-reg}, we shall see in Corollary \ref{cor:main-sys} that minimizing constraint maps solve the system
\begin{equation}\label{eq:main-sys}
\Delta u = A_u(Du,Du)\chi_{u^{-1}(\partial M)} + q \nu_u \cH^{n-1}|_{\partial_{\rm red} u^{-1}(M)}\quad\text{in }\Omega\setminus\Sigma(u), 
\end{equation}
(in a suitable sense which we shall clarify in Section \ref{sec:e-reg}), where $A_p$ is the second fundamental form of $\partial M$ at $p$, see \eqref{eq:identityhessrho}. In this sense, Theorem \ref{thm:e-reg} provides a good picture of the general behaviour of minimizing constraint maps, without any further assumptions on the target $M$. 
\end{rem}

Theorem \ref{thm:e-reg} raises the question of whether there are singularities present apart from those resulting from the part of the domain, ${\rm int}(u^{-1}(\partial M))$, where the constraint map is a harmonic map. To answer this question, we assume that the complement of $M$ is {\it uniformly convex}.\footnote{Uniform convexity of $M^c$ means that there exists $R>0$ such that the following holds: For any point $x \in \partial M$ there exists a ball $B_R(y)\supset M^c$ such that $x \in \partial B_R(y)$.}
As we show in Section \ref{sec:contrho}, mere convexity guarantees that the non-coincidence set $u^{-1}(M)$ is open and, therefore, there is a well-defined free-boundary $\partial u^{-1}(M)\cap \Omega$.

\begin{thm}\label{thm:main}
Let $M$ be a $C^3$-domain with uniformly convex complement. Let $u \in W^{1,2}(\Omega;\overline M)$ be a minimizing constraint map for $\cE_q$ with $q>0$, such that $\int_\Omega |Du|^2\,dx\leq \Lambda$. For every $\eta > 0$, there exist constants $\delta_\eta \in (0,\eta)$ and $c_\eta > 1$, both depending only on $n$, $m$, $q$, $\eta^{-n}\Lambda$, and $\partial M$, such that 
$$
\|Du\|_{L^\infty(B_{\delta_\eta}(u^{-1}(M))\cap\Omega_\eta)} \leq c_\eta,
$$
where $\Omega_\eta := \{ x \in\Omega:\dist(x,\partial \Omega)>\eta\}$.
\end{thm}

\begin{rem}
As our map is a minimizing harmonic map into $\partial M$ in the interior of the coincidence set, ${\rm int}(u^{-1}(\partial M))$, the above theorem indicates that the only singularities of the minimizing constraint maps are those arising from harmonic maps. 
\end{rem}

\begin{rem}
We do not know the minimal geometric conditions on the target constraint $\overline M$ under which the above theorem holds. Leaving  aside constraints with non-convex complement, even the case where  the complement is merely convex already raises some subtle issues.

Nonetheless, our proof of Theorem \ref{thm:main} extends to the case where the complement of $M$ is locally uniformly convex and has bounded principal curvatures, which may then include unbounded sets, e.g.\ a paraboloid. We note that, in the special case when $M$ is graphical and concave, it follows from the results of \cite{F2} that singularities do {\it not} appear in such scenarios, see Appendix \ref{sec:remove} for further details.    
\end{rem}

\begin{rem}
An analogue of the above theorem for the
case $q = 0$ poses additional significant challenges and is the subject of our forthcoming research in \cite{FGKS}.
\end{rem}

Since, by the result above, minimizing constraint maps are regular in a neighborhood of the free boundary, we are led to consider the regularity of the free boundary itself. An important ingredient in the proof of Theorem \ref{thm:main}, which we will discuss in more detail below, is that the distance of $u$ to $\partial M$ is a minimizer for a scalar Alt--Caffarelli-type functional, so we deduce from \cite{GS} (see also \cite{DeSilva2011} for a viscosity approach)
that flat free boundaries are locally $C^{1,\alpha}$-graphs. 
Then by a bootstrapping argument analogous to the one in \cite{FKS}, we obtain that the graphs are indeed of class $C^\infty$. More precisely, we have the following theorem.

\begin{thm}\label{thm:fb-reg}
Assume, in addition to the assumptions in Theorem \ref{thm:main}, that $\partial M$ is of class $C^{k+1}$, $k\geq 2$. Then there exists a constant $\e > 0$, depending only on $n$, $m$, $q$, $\Lambda,$ and the $C^3$-character of $\partial M$, such that if for some $e \in \partial B_1$, $x_0\in \Omega_\eta\cap \partial u^{-1}(M)$, and $r\in(0,\frac{1}{4}\delta_\eta)$, 
\begin{equation}\label{eq:flat}
\begin{aligned}
& B_{4r}(x_0)\cap \{ x: (x-x_0)\cdot e > \e r\} \subset u^{-1}(M), \\
& B_{4r}(x_0)\cap \{ x: (x-x_0)\cdot e < -\e r\} \subset u^{-1}(\partial M), 
\end{aligned}
\end{equation}
then $B_r(x_0)\cap \partial u^{-1}(M)$ is a $C^{k,\gamma}$-hypersurface, for every $\gamma\in(0,1)$, whose $C^{k,\gamma}$-character depends, in addition to the aforementioned quantities, on $k$, $\gamma$, and $r$. 
\end{thm}

By  \cite{CJK,DSJ,JS}, the only scalar $1$-homogeneous global minimizers are half-space solutions for low dimensions. If we denote by $\rho$ the signed distance function to $\partial M$ then $\rho\circ u$ is a minimizer of a scalar Alt-Caffarelli energy with lower-order perturbation, which vanishes along $1$-homogeneous rescalings, and so we can verify the flatness condition \eqref{eq:flat} for every free boundary point, in such low dimensions; the argument is the same as in \cite[Theorem 2.17 (a)]{GS} (see also \cite[Theorem 8.3]{AC}). Consequently, we obtain full regularity of the free boundaries in low dimensions.

\begin{cor}\label{cor:fb-reg}
There is a number $n_* \in \{5,6,7\}$ such that if $\Omega$ lies in $\R^n$ with $n \leq n_*$, and the assumptions of Theorem~\ref{thm:fb-reg} are satisfied, then $\Omega\cap \partial u^{-1}(M)$ is a $C^{k,\gamma}$-hypersurface.
\end{cor}

In our work, an important role is played by the decomposition
$$u=\Pi \circ u + (\rho\circ u)\nu\circ u,$$
cf.\ \eqref{eq:decompu}, where as above $\rho$ is the signed distance to $\partial M$, $\Pi$ is the nearest-point projection onto $\partial M$ and $\nu$ is the inwards pointing unit normal to $\partial M$. In order to prove the theorems stated above, we will need to study in detail the regularity properties of $\rho\circ u$ (see e.g.\ Lemma \ref{lem:min} and Theorem \ref{thm:dist-reg}). It is then natural to also investigate the regularity properties of $\Pi \circ u$. The regularity of $\Pi \circ u$ was studied in-depth in \cite{FKS} for minimizing constraint maps for the Dirichlet energy $\cE_0$. Notably, the analysis of the projection map $\Pi\circ u$ exhibits vectorial characters, and opens up new, interesting problems to be further explored. Here, we shall study these issues in the Bernoulli setting, obtaining the following result:

\begin{thm}\label{thm:proj-reg}
In the setting of Theorem \ref{thm:main}, we have
$$
\| D^2 (\Pi \circ u) \|_{\textup{BMO}(B_{\delta_\eta}(u^{-1}(M))\cap\Omega_\eta)} \leq c_\eta.
$$
Moreover, if \eqref{eq:flat} holds for some $x_0\in\Omega_\eta\cap \partial u^{-1}(M)$ and $ r\in (0,\frac{1}{4}\delta_\eta)$ and $\partial M$ is of class $C^{3,\gamma}$ for some $\gamma\in(0,1)$, then 
$$
\| D^2 (\Pi \circ u) \|_{L^\infty(B_r(x_0))} \leq c_{\eta,\gamma},
$$
where $c_{\eta,\gamma}$ depends on $\gamma$, in addition to the quantities that determine $c_\eta$.
\end{thm}

In light of Corollary \ref{cor:fb-reg}, it follows that $|D^2(\Pi\circ u)| \in L^\infty(B_r(x_0))$ at every free boundary point $x_0\in \partial u^{-1}(M)\cap \Omega$ in low dimensions, that is, whenever $n\leq n_*$. This leaves the study of $C^{1,1}$ regularity close to singular points as a tantalizing open problem. 

We also remark that when dealing with the obstacle problem (specifically when $q = 0$), a significant advancement was made in two dimensions in the work of \cite{FKS} regarding the sharp result for $\Pi\circ u$, as long as the distance component $\rho\circ u$ allows for a certain geometric approximation (see \cite{FKS} Definition 2.8).
The question of achieving optimal regularity in higher dimensions for the case $q=0$ remains unresolved.


\subsection{Idea of the proofs}

Here, we delve into the key concepts in our methodology intended for those well-versed in the subject matter. We shall also illuminate aspects of the theory and approach that are entirely new, both at the technical and conceptual levels.
Readers less familiar with the topic may skip this subsection on a first reading.

\subsubsection{Theorem \ref{thm:e-reg}}

Our first main result is the $\e$-regularity theorem, which asserts that the map is regular near points of small normalized energy. In the classical setting of harmonic maps, smallness of the normalized energy entails that the map is approximately an unconstrained minimizer of the Dirichlet energy. Such minimizers are of course regular, and hence the original map inherits (almost) the same regularity. 

Our problem is substantially different, as the characteristic function $\chi_{u^{-1}(M)}$ involved in the Alt--Caffarelli energy $\cE_q$ is discontinuous in $u$.
Our key observation is that, at points of small energy, the map can be approximated by maps that satisfy {\it a priori} Lipschitz estimates (Proposition \ref{prop:apprx}). The approximating maps are not solutions to a single minimization problem, but instead, there are two types of approximating maps:  either (i) each of their components is harmonic,  or (ii) they are minimizing maps of a new Alt--Caffarelli energy subject to a flat half-space constraint. We now detail each of these cases:

\begin{enumerate}[(i)]
    \item in this case, the target constraint disappears and  we are essentially in the classical setting of harmonic maps, hence the Lipschitz regularity is a consequence of Luckhaus' extension lemma \cite[Lemma 1]{L};
    \item in this case, the flat constraint induces a decoupling in the problem: the components of the map parallel to the half-space are harmonic, while the normal component is a  minimizer for a \textit{scalar} Alt--Caffarelli energy, for which the optimal regularity is known from \cite{AC}. To construct the approximating sequence, we again rely on Luckhaus' lemma, but we also have to use boundary flattening maps in a very precise way to ensure that the limit constraint is flat.
\end{enumerate}

The above Lipschitz approximation result implies $\alpha$-H\"older regularity of minimizing constraint maps near points of small normalized energy, for every $\alpha \in (0,1)$. To reach the endpoint exponent $\alpha=1$, we have to further exploit the precise structure of $\cE_q$. 

More precisely, we first work in Sobolev spaces, and we prove $W^{1,p}$-estimates for any $p<\infty$ (Lemma \ref{lem:e-reg}). These estimates seem to be new even in the classical case $q=0$, and the main challenge here is to combine the almost-monotonicity formula with maximal function estimates. To conclude the proof and get a $W^{1,\infty}$-estimate, we observe that the tangential part $\Pi\circ u$ has one more derivative than $u$, so it belongs to $W^{2,p}$ (Lemma \ref{lem:proj-W2p}). In addition, the normal part $\rho\circ u$ turns out to be a minimizer of a scalar Alt--Caffarelli functional, modulo some lower order perturbation (Lemma \ref{lem:min}), so that it is of class $W^{1,\infty}$. Combining these two ingredients, we establish the optimal $\e$-regularity theorem.

\subsubsection{Theorem \ref{thm:main}}

Our second result, which is in fact what motivated this paper, concerns the improved regularity of the map near the free boundary. This is a topic that has not been explored in the literature on harmonic maps.\footnote{We should remark here that our mappings in the manifold level are minimizing harmonic maps {\it into} manifolds-with-boundary, rather than those {\it from } manifolds-with-boundary. There is a large amount of literature on the latter subject, in which case the {\it free boundary} is the {\it image} of the prescribed boundary in the source. In this case, the singularities appear in general on the Dirichlet boundary. In our case, the {\it free boundary} is the topological boundary of the {\it preimage} of the interior of the target constraint, whence the problem is very different.} 
Our argument is based on compactness techniques, and thus can be applied to a wide variety of energy functionals involving other types of nonlinearities\footnote{E.g., the Alt--Phillips functionals, where the characteristic function $\chi_{u^{-1}(M)}$ will then be replaced with some power of the distance map, $(\rho\circ u)^\gamma$.}. However, as we shall specify below, there is one component of our argument that cannot be carried over to the case of the Dirichlet energy $\cE_0$, as it relies crucially on the strict positivity of the weight $q > 0$ for $\cE_q$. The case $q = 0$, which is much more delicate and involved, will be studied in our forthcoming paper \cite{FGKS}. 

As in the statement of Theorem \ref{thm:main}, we now assume that the complement $M^c$ of our target constraint, which we will call the obstacle, is uniformly convex.
It is unclear how much this condition can be relaxed, although we believe that the conclusion of Theorem \ref{thm:main} is generally false if $M^c$ is not convex. 

A first important observation is that, even though $u$ is discontinuous in general, the distance map $\rho\circ u$ is \textit{always continuous} whenever $q\geq 0$  (Theorem \ref{thm:dist-reg}). This follows from the fact that $\rho\circ u$ is subharmonic (Corollary \ref{cor:subharm}) together with an analysis of the tangent maps (Definition \ref{defn:tan}). The continuity of $\rho \circ u$ has several immediate but important consequences:
\begin{enumerate}[(i)]
    \item there is a well-defined free boundary $\Omega\cap\partial u^{-1}(M) = \Omega\cap \partial\{ \rho\circ u > 0\}$;
    \item each component of $u$ is a scalar harmonic function in the open set $u^{-1}(M)$;
    \item the singularities of $u$ are contained in the coincidence set $u^{-1}(\partial M)$ and, at singular points, tangent maps are minimizing harmonic maps into $\partial M$. 
\end{enumerate}

In order to carry our analysis further, we upgrade (iii) to a much stronger result about the original map: we show that the density of the non-coincidence set $u^{-1}(M)$ has to vanish at every singular point (Proposition \ref{prop:dense}). 
This indicates that the free boundary (if any) around a singular point can be, at most, cusp-like. When combined with the subharmonicity of the distance map $\rho\circ u$, the vanishing density of the non-coincidence set implies that $\rho\circ u$ decays to infinite order at singular points (Corollary \ref{cor:dense}).  
Now, to complete the proof, we rule out the possibility of having singularities at the free boundary by using the assumption $q>0$ to establish the linear nondegeneracy of $\rho \circ u$ at any free boundary point (Proposition \ref{prop:ndeg}).

Let us make a final comment on Theorem \ref{thm:main}, which asserts that singularities are pushed away from the free boundary towards the interior of the coincidence set. In fact, by the compactness method, we can establish universal {\it a priori} estimates for both the minimal distance of the singularities to the free boundary and for the Lipschitz norm of our map near the free boundary. It is worth comparing our result with the boundary regularity theory for minimizing harmonic maps \cite{SU2}, since minimizing constraint maps $u$ are minimizing harmonic maps when restricted to the interior of the coincidence set, ${\rm int}(u^{-1}(\partial M))$. However, we do not have any {\it a priori} control over the regularity of the free boundary $\Omega \cap \partial u^{-1}(M)= \Omega \cap \partial u^{-1}(\partial M)$, and in general the free boundary has singularities of its own, even for the scalar problem. Therefore, our paper yields a striking result in the context of the regularity theory of harmonic maps, and at the same time it opens up uncharted areas from the free boundary perspective.

\subsubsection{Theorems \ref{thm:fb-reg} and \ref{thm:proj-reg}}
With Theorem \ref{thm:main} at our disposal, we can study further properties of the free boundaries by combining the partial hodograph-Legendre transformation with the regularity estimates for elliptic transmission problems. 

We remark that once we are at a point of continuity of the map, the minimality of the map does not play a central role anymore. For this reason, most of the arguments for Theorems \ref{thm:fb-reg} and \ref{thm:proj-reg} can be generalized to constraint maps that are not necessarily minimizers.


\subsection{Organization of the paper}

The main body of our paper can be divided into three parts: 
\begin{enumerate}
    \item[(i)] $\e$-regularity theory (Sections \ref{sec:basicanalysis}--\ref{sec:e-reg});
    \item[(ii)] regularity improvement near free boundaries (Sections \ref{sec:contrho}--\ref{sec:reg}); 
    \item[(iii)] regularity of  free  boundaries, and image map  (Sections \ref{sec:fb-reg}--\ref{sec:proj}). 
\end{enumerate}

More specifically, our paper is organized as follows: In Section \ref{sec:prelim}, we provide the precise setting of our problem, introduce the notation and terminology, and present some preliminary lemmas. In Section \ref{sec:basicanalysis}, we derive basic equations and inequalities satisfied by constraint maps with a sufficient amount of regularity. In Section \ref{sec:apprx}, we prove the Lipschitz approximation, Proposition \ref{prop:apprx}. In Section \ref{sec:e-reg}, we study the sharp $\e$-regularity theorem, and prove Theorem \ref{thm:e-reg}. Up to this point, we do not assume any geometric condition on $\partial M$ other than regularity. 

We impose the convexity assumption starting from Section \ref{sec:contrho} onwards and build up our argument for Theorem \ref{thm:main}. More specifically, in Section \ref{sec:contrho}, we establish the continuity of $\rho\circ u$. Section \ref{sec:dense} studies the (vanishing) density of the non-coincidence set at singular points. Then, in Section \ref{sec:ndeg}, we establish the non-degeneracy of $\rho\circ u$ at free boundary points. In Section \ref{sec:reg}, we combine the previous results to prove Theorem \ref{thm:main}. 

Section \ref{sec:fb-reg} concerns the regularity of the flat free boundaries and contains the proof of Theorem \ref{thm:fb-reg}. Finally, we close this paper with Section \ref{sec:proj}, where we take a closer look at the projected image of the map and prove Theorem \ref{thm:proj-reg}.



\subsection{Notation}

Here is a list of notation to be used throughout this paper.

\begin{longtable}[l]{l l }
\,\,\,\,$\Omega$ &\qquad  a bounded domain in $\R^n$ ($n\geq 1$)  \medskip \\ 
\,\,\,\,$M$  &\qquad smooth domain in $\R^m$ ($m\geq 2$)  \medskip \\  
\,\,\,\,$\cE_q[v]$  &\qquad $\int_\Omega (|Dv|^2 + q^2 \chi_{v^{-1}(M)}) \, dx$,  \quad $q\geq 0$    \medskip \\   
\,\,\,\,$\Sigma(u)$   &\qquad  $ \{ x\in\Omega: u\text{ is not continuous at }x\} $    \medskip \\   
\,\,\,\,$\Omega_\eta $ &\qquad  $\{ x \in\Omega:\dist(x,\partial \Omega)>\eta\}$      \medskip \\   
\,\,\,\,$A_p (\cdot , \cdot)$  &\qquad     the second fundamental form of $\partial M$ at $p$, see \eqref{eq:identityhessrho}   \medskip \\   
\,\,\,\,$\cH^k$  &\qquad  $k$-dimensional Hausdorff measure      \medskip\\   
\,\,\,\,$\hbox{dim}_\cH$  &\qquad  Hausdorff dimension     \medskip\\   
\,\,\,\,$\rho$  &\qquad    signed distance function to $\partial M$    \medskip \\   
\,\,\,\,$\nu=\nabla \rho$  &\qquad inward unit normal to $M$ \medskip \\   
\,\,\,\,$\cN(\partial M)$  &\qquad    tubular neighborhood  of $\partial M$ 
    \medskip\\   
\,\,\,\,$\Pi$   &\qquad  nearest point projection: $ \overline{M}\to \partial M$
       \medskip\\ 
\,\,\,\,$\nu\otimes \nu$   &\qquad        tensor product $(\nu^i\nu^j)_{ij}$ \medskip\\   
\,\,\,\,$\xi^\top$  &\qquad   $
  (I - \nu\otimes \nu)\xi
$,   \ orthogonal component of $\xi$  \medskip\\   
\,\, Tangent maps &\qquad  see Definition \ref{defn:tan}      \medskip\\   
\,\,\,\,$E(u,x_0,r)$  &\qquad   $ r^{2-n} \int_{B_r(x_0)} |Du|^2\,dx$, \ normalized energy     \medskip\\   
\,\,\,\,$F:G$  &\qquad  $ f_\alpha^i g_\alpha^i$, \ where  $F = (f_\alpha^i)$ and  $G= (g_\alpha^i)$     \medskip \\   
\,\, $\cM f(x_0) $ &\qquad $  \sup_{r> 0} \left[ \frac{1}{r^n}\int_{B_r(x_0)} f\,dx \right] $ \quad (maximal function)    \medskip \\ 
\,\,\,\,$W^{1,p}(\Omega;\overline M)$&\qquad $\{\phi\in W^{1,p}(\Omega;\R^m): \phi\in \overline M\text{ a.e.\ in }\Omega\}$      \medskip \\ 
\end{longtable}


\section{Preliminaries}\label{sec:prelim}

\subsection{Problem setting}

In this paper $M\subset\R^m$ denotes a domain whose complement is a closed set of class $C^3$. We denote by $\rho$ the  signed distance function to $\partial M$, so that $\rho$ is positive in $M$ and is $C^3$ in a tubular neighborhood $\cN(\partial M)$ of $\partial M$.  One can thus extend the inward unit normal $\nu$ (i.e., the normal which points towards $M$) to $\cN(\partial M)$, for instance through $\nu:= \nabla \rho$.  Let $\Pi\colon \cN(M)\to \partial M$ denote the nearest-point projection to $\partial M$. 
Then $\nu = \nu\circ\Pi$ in $\cN(\partial M)$ and thus the decomposition 
\begin{equation}\label{eq:w-V}
\id := \Pi + \rho\nu = \Pi + \rho (\nu\circ \Pi)
\end{equation}
holds everywhere in $\cN(\partial M)$. 

The linear map $I - \nu\otimes \nu$, when evaluated at a point $x^0 \in \partial M$, corresponds to the orthogonal projection to the tangent space of $\partial M$ at $x^0$.
For a vector $\xi\in\R^m$, we denote by
$$
\xi^\top  := (I - \nu\otimes \nu)\xi
$$
its orthogonal component, as in \cite{Simon}. We define, through the same formula, $X^\top$ for a matrix $X\in \R^{m\times n}$: this is the matrix that results from $X$ by applying $\cdot^\top$ to each of its rows.  

The second fundamental form $A$ of $\partial M$ then can be expressed in terms of the Hessian of $\rho$ via
\begin{equation}
    \label{eq:identityhessrho}
    \Hess\rho(\xi,\xi) = \Hess\rho(\xi^\top,\xi^\top) = -\nu\cdot A(\xi^\top,\xi^\top).
\end{equation}
We remark that these definitions are canonically extended to the tubular neighborhood $\cN(\partial M)$ through the identity $\nu = \nu \circ \Pi$. 

Let us also record a useful identity: 
\begin{equation}\label{eq:DPi-Dnu}
(\nabla \Pi)( \nabla \nu) = (I - \rho \nabla\nu)\nabla\nu,
\end{equation}
in $\cN(\partial M)$; in coordinates, this reads as $\partial_i\Pi^k \partial_j \nu^k = (\delta_{ik}- \rho\partial_i\nu^k)\partial_j\nu^k$, where $\delta_{ik}$ is Kronecker's delta.

Let $M$ be a $C^2$-domain in $\R^m$ with uniformly bounded principal curvatures, and set $\rho$, $\nu$, $\Pi$ and $\cN(\partial M)$ as above. Then there exists some $\kappa > 1$, such that 
\begin{equation}\label{eq:curv}
|\Hess\rho| \leq \kappa\quad\text{in }\cN(\partial M),
\end{equation}
where the tubular neighborhood $\cN(\partial M)$ satisfies
\begin{equation}\label{eq:width}
|\rho| \geq \frac{1}{c\kappa}\quad\text{on }\partial \cN(\partial M);
\end{equation}
see \cite[Lemmas 14.16--17]{GT} for a proof. 
Let $\cN(M) := \cN(\partial M)\cup M$. Then we can define the nearest point projection $\pi \colon \cN(M) \to \overline M$ by 
\begin{equation}\label{eq:pi}
\pi (a) := \begin{cases}
a, &\text{if }a \in \overline M,\\
\Pi(a), &\text{if }a \in \cN(\partial M)\setminus\overline M. 
\end{cases}
\end{equation} 
The tangent cone  $T_a\overline M$ is also defined by 
\begin{equation}\label{eq:tan}
T_a\overline M := \begin{cases}
    \R^m, &\text{if }a\in M,\\
    \{\xi \in \R^m: \nu_a \cdot \xi \geq 0\}, &\text{if }a\in \partial M;
\end{cases}
\end{equation}
where $\nu_a$ is the inward unit normal to $\partial M$ at $a$. Here, and throughout the paper, we adopt the convention of using a subscript to denote the point in $\R^m$ where a function is evaluated. Now by \eqref{eq:curv} (along with \eqref{eq:pi} and \eqref{eq:w-V}), we have, for each $i\in \{0,1\}$,
\begin{equation}\label{eq:id-pi}
|\nabla^i(\id - \pi)| \leq c\kappa \e^{2-i}\quad \text{in }B_\e(a)\cap (a + T_a\overline M),
\end{equation}
for all $a\in \overline M$ and every $\e>0$ such that $B_\e(a)\subset \cN(M)$; here $c\equiv c(m)$. 

For definiteness, let $\eta$ be the width of $\cN(\partial M)$; by \eqref{eq:width}, $c\kappa \eta > 1$. Given a point $a \in \partial M$, we can also consider the boundary flattening map, i.e., a $C^1$-diffeomorphism $\Phi\colon B_\eta(a)\to U$, for some open neighborhood $U\subset \R^m$ of the origin, such that  
\begin{equation}\label{eq:Phi}
\Phi (M \cap B_\eta(a)) = U^+ := \{ y \in U: y^m > 0\},\quad \Phi(a) = 0. 
\end{equation}
We may also choose $\Phi$ such that  
\begin{equation}\label{eq:id-Phi}
|\nabla^i(\cR -\cR a - \Phi)| \leq c\kappa \e^{2-i}\quad\text{in }B_\e(a),
\end{equation}
for every $\e \in (0,\eta)$ and $i\in \{0,1\}$, for a rotation $\cR$ which maps $\nu_a$ to $e_m$.

As is customary in the literature on harmonic maps \cite{Simon}, by $A(Du,Du)$ we denote the summation of $A(D_\alpha u, D_\alpha u)$ over $\alpha \in \{1,2,\ldots,n\}$; here and thereafter we follow the convention of summation with repeated indices and $\alpha,\beta$ (and $i,j$) will run through $\{1,2,\ldots,n\}$ (resp. $\{1,2,\ldots,m\}$).

The decomposition \eqref{eq:w-V} induces a similar decomposition for $u$ in the pre-image $u^{-1}(\cN(\partial M))$ of the tubular neighborhood:
\begin{equation}
    \label{eq:decompu}
    u = \Pi\circ u + (\rho\circ u) \nu\circ u
    = \Pi\circ u + (\rho\circ u) \nu\circ{\Pi\circ u }.
\end{equation}
A straightforward yet valuable insight, as presented by \cite{FKS}, is that $\Pi\circ u$ and $\rho\circ u$ exhibit different levels of regularity. Specifically, the former possesses an additional derivative compared to the latter, as indicated in Corollary \ref{cor:proj-pde} below. This implies that the tangential component of $Du$ exhibits greater regularity than the normal component, as evidenced by the identity
\begin{equation}
    \label{eq:tangpart}
    (Du)^\top=D(\Pi\circ u) + (\rho\circ u)D(\nu\circ u),
\end{equation}
where derivatives of $\rho\circ u$ are absent.

It is important to emphasize that the first part  of this paper, encompassing Sections \ref{sec:basicanalysis} through \ref{sec:e-reg} and leading to the $\e$-regularity theorem (Theorem \ref{thm:e-reg}), applies to any domain with a $C^3$ smoothness condition and uniformly bounded principal curvatures. Notably, this includes cases where the complement of $M$ can be unbounded. To be more specific, our analysis will primarily be conducted under the assumption specified in equation \eqref{eq:curv}.

In the rest of the paper (Sections \ref{sec:contrho}--\ref{sec:proj}) on the regularity improvement near free boundaries, we shall assume that the complement of $M$ is uniformly convex (in addition to the $C^3$-regularity). More precisely, we will typically work under the hypothesis that (assuming $\kappa > 1$ in \eqref{eq:curv})
\begin{equation}
    \label{eq:boundshessrho}
    \frac{1}{\kappa}|\xi^\top|^2 \leq \Hess \rho_a(\xi,\xi) \leq \kappa |\xi^\top|^2
\end{equation}
for all $\xi\in \R^m$ and $a\in \cN(\partial M)$.

\subsection{Compactness results for the constraint maps}\label{sec:cpt}
In our analysis we will often study sequences of rescalings of a given minimizing constraint map. The rescalings compatible with our problem are 0-homogeneous and, as in the literature of harmonic maps \cite{GM,LW,Simon}, we shall call the limit maps obtained by this blowup procedure {\it tangent maps}: 

\begin{defn}[Tangent maps]\label{defn:tan}
Let $u\in W^{1,2}(\Omega;\overline M)$ be a minimizing constraint map for the functional $\cE_q$, where $q\geq 0$. For $x_0\in\Omega$ we call $\phi \in W_\loc^{1,2}(\R^n;\overline M)$ a {\it tangent map of $u$ at $x_0$}, and write $\phi \in T_{x_0}u$, if there is a sequence $r_k\to 0$ of radii and a rescaled sequence 
$$u_{x_0,r_k}(y) := u(x_0+ r_k y)$$ such that $u_{x_0,r_k} \rightharpoonup \phi$ weakly in $W_{\loc}^{1,2}(\R^n;\R^m)$. 

\end{defn}

Associated to the above rescalings we have a scale-invariant quantity that we refer to as the \textit{normalized energy} and that we denote by
$$
E(u,x_0,r) := r^{2-n} \int_{B_r(x_0)} |Du|^2\,dx.
$$
Indeed, it is easy to see that $E(u,x_0,rs)=E(u_{x_0,r},0,s).$

Any minimizing constraint map is necessarily stationary upon domain variations.  Thus, by a standard argument which for completeness we present in Appendix \ref{sec:tech}, we obtain an almost monotonicity of the normalized energy; when $q=0$, as is well-known, this energy becomes fully monotone.

\begin{lem}[Almost monotonicity formula]\label{lem:monot}
Let $u\in W^{1,2}(\Omega;\overline M)$ be a minimizing constraint map for $\cE_q$, $q\geq 0$, and let $x_0\in\Omega$. Then 
$$
\begin{aligned}
E(u,x_0,s) - E(u,x_0,r) \geq 
2\int_{B_s(x_0)\setminus B_r(x_0)} R^{2-n} \left|\frac{\partial u}{\partial N}\right|^2 \,dx 
 - \omega_n q^2 s^2, 
\end{aligned}
$$ 
whenever  $0<r<s<\dist(x_0,\partial\Omega)$, where $R = |x-x_0|$ and $\partial/\partial N$ denotes the directional derivative in the radial direction $R^{-1}(x-x_0)$.  Moreover, if $E(u,x_0,r) = E(u,x_0,s)$ and $q=0$, then $u$ is $0$-homogeneous about $x_0$ in the annulus $B_s(x_0)\setminus B_r(x_0)$, i.e.,  $$|(x-x_0)\cdot Du(x)| = 0 \text{ for a.e.\ } x\in B_s(x_0)\setminus B_r(x_0).$$ 

\end{lem}

Minimizers (or almost minimizers) of $\cE_0$ 
exhibit crucial compactness characteristics, as confirmed in \cite{L} and \cite{SU}, and these traits are fundamental in the examination of tangent maps. We will now present several lemmas that mirror these compactness outcomes for the Alt-Caffarelli energy.

We will commence with the compactness result. Although Luckhaus's findings in \cite{L} may not be directly applicable, it is apparent that the methodologies employed to establish these results remain applicable in our context. For the readers' convenience, we have included a proof in  Appendix \ref{app:compact}.

\begin{lem}[Compactness {\cite{L}}]\label{lem:strong}
Consider a bounded sequence $\{u_k\}_{k=1}^\infty$ of minimizing constraint maps in the space $W^{1,2}(\Omega;\overline M)$ with respect to the functionals $\cE_{q_k}$, where $q_k\geq 0$. If $q_k \to q$ for some $q\geq 0$, then there exists a minimizing constraint map $u \in W^{1,2}(\Omega;\overline M)$ for the functional $\cE_q$ such that $u_k\to u$ strongly in $W^{1,2}(\Omega;\R^m)$ along a subsequence.
\end{lem}

 A typical application of Lemma \ref{lem:strong} concerns sequences $\{u_{x_0,r_k}\}_{k=1}^\infty$ with $r_k\to 0$. By Lemma \ref{lem:monot}, the normalized energy density for  minimizing constraint maps for the functional $\cE_q$ is well-defined:
$$
\Theta_u(x_0) := \lim_{r\to 0}\left[r^{2-n}\int_{B_r(x_0)} |Du|^2\,dx \right].
$$
 As a useful byproduct of Lemmas \ref{lem:monot} and \ref{lem:strong} we obtain the upper semicontinuity of the normalized energy density:

 \begin{cor}\label{cor:semi}
 Let $\{u_k\}_{k=1}^\infty$ be a sequence of minimizing constraint maps in $W^{1,2}(\Omega;\overline M)$ for the functionals $\cE_{q_k}$, where $q_k\geq 0$, such that $u_k\to u$ strongly in $W^{1,2}(\Omega;\R^m)$ and $q_k\to 0$.  For any sequence $\{y_k\}_{k=1}^\infty$ in $\Omega$, if $y_k \to y_0\in\Omega$ as $k\to\infty$, then 
 $$
 \Theta_u(y_0) \geq \limsup_{k\to\infty} \Theta_{u_k}(y_k). 
 $$
 \end{cor}

The following is another useful corollary of Lemmas \ref{lem:monot} and \ref{lem:strong}:

\begin{cor}\label{cor:blowup}
Let $u\in W^{1,2}\cap L_\loc^\infty(\Omega;\overline M)$ be a minimizing constraint map for the functional $\cE_q$, where $q\geq 0$. Then $T_{x_0}u \neq \emptyset$, whenever $x_0\in\Omega$. In fact, every $\phi\in T_{x_0}u$ is a minimizing constraint map in $\R^n$ for the functional $\cE_0$. Moreover,
$$
r^{2-n}\int_{B_r(0)} |D\phi|^2\,dy = \Theta_u(x_0)
$$
for all $r>0$. In particular, any $\phi \in T_{x_0} u$ is homogeneous of degree $0$, i.e., $|(y\cdot D)\phi| = 0$ a.e. in $\R^n$.

\end{cor}

We remark that in case $u$ is continuous the tangent cone at $x_0$ contains only one element, namely $u(x_0)$. Otherwise elements like $x/|x|$ are inside the tangent cone.

\begin{proof}
Let $r_k\to 0$ be an arbitrary sequence of radii and write $\phi_k (y) := u_{x_0,r_k}$ for short. Clearly, for every ball $B_r(0)\subset\R^n$, $\phi_k \in W^{1,2}\cap L^\infty(B_r(0);\overline M)$ for all $k$ sufficiently large. One can also directly verify that $\phi_k$ is a minimizing constraint map associated with $q_k:= r_k q \to 0$. Moreover, by Lemma \ref{lem:monot}, 
$$
\begin{aligned}
r^{2-n} \int_{B_r(0)} |D\phi_k|^2\,dy &= (r_kr)^{2-n} \int_{B_{r_k r}(x_0)} |Du|^2\,dx \leq R_0^{2-n}\int_{B_{R_0}(x_0)} |Du|^2\,dx + \omega_n q^2 r_0^2, 
\end{aligned}
$$
for all large $k$, where $R_0 := \dist(x_0,\partial\Omega)$. Thus, Lemma \ref{lem:strong} (along with the Rellich compactness theorem) yields a minimizing constraint map $\phi\in W_{\loc}^{1,2}(\R^n;\overline M)$ for  $\cE_0$ such that $\phi_k \to \phi$ strongly in $W^{1,2}(B_r(0);\R^m)$ up to a subsequence, which can be chosen independent of $r$. This proves $T_{x_0}u \neq \emptyset$. The remaining assertions follow easily from Lemmas \ref{lem:strong} and  \ref{lem:monot}.
\end{proof}

\subsection{Some basic lemmas}

In this short subsection we collect, for later use, a few elementary results of general character.

The first such result states that $0$-homogeneous weakly subharmonic and bounded functions are necessarily constants.

\begin{lem}\label{lem:0-hom}
Let $B$ be a ball centered at the origin, and $w \in W^{1,2}(B)$ a nonnegative, weakly subharmonic function. If $w$ is $0$-homogeneous in $B$, then $w$ is constant a.e.\ in $B$.
\end{lem}

\begin{proof}
Since $w$ is $0$-homogeneous in $B$, we have   \begin{equation}\label{eq:0-hom}
x\cdot Dw = 0\quad\text{for a.e.\ } x \in  B.
\end{equation}
Now let $r \in (0,1)$ be given, and let us consider a radially-symmetric and smooth cutoff function $\eta \in C_0^\infty (B)$ such that $\eta = 1$ in $r B$ and $0\leq \eta\leq 1$ in $B$. 
By the radial symmetry of $\eta$, there exists some $\tilde \eta\in C_c^\infty([0,\infty))$ such that $D\eta(x) = \tilde\eta(|x|)x$. Testing the weak 
subharmonicity of $w$ against $\eta w \in W_0^{1,2}(B)$, which is nonnegative, we obtain 
$$
\int_B \eta |Dw|^2 \,dx \leq -\int_B w Dw \cdot D\eta \,dx.
$$ 
However, by \eqref{eq:0-hom}, $Dw\cdot D\eta = \tilde\eta(|x|) x\cdot Dw = 0$ a.e.\ in $B$. Thus, since $0\leq\eta\leq 1$ with $\eta=1$ in $rB$, we have
$$
|Dw|^2 = 0\quad\text{a.e.\ in } r B.
$$
As $r\in(0,1)$ was arbitrary, we deduce that $|Dw| = 0$ a.e.\ in $B$, proving the assertion. 
\end{proof}

The next lemma is an elementary inequality for Sobolev functions. 

\begin{lem}\label{lem:W11}
Let $w \in W^{1,1}(B_r(x_0))$ be a nonnegative function. Then 
\begin{equation}\label{eq:W11}
\int_{\partial B_r(x_0)} w \,d\sigma \leq \frac{n}{r} \int_{B_r(x_0)} w\,dx + \int_{B_r(x_0)} |Dw|\,dx.
\end{equation}
\end{lem}

\begin{proof}
Without loss of generality we take $r=1$ and $x_0=0$. After integrating by parts in $r$, we have

\begin{align*}
   n \int_{B_1} w \, dx & = 
   n \int_0^1 r^{n-1}\int_{\partial B_1}
   w(r \theta) \, d\sigma_\theta \,dr\\
    & = 
    \int_{\partial B_1} w \, d \sigma_\theta -  \int_0^1 r^{n-1} \int_{\partial B_1}  \theta\cdot D w(r\theta)\, d\sigma_\theta \, dr\\
    &\geq \int_{\partial B_1} w \, d \sigma_\theta 
     - \int_{B_1} |Dw| \,dx,
\end{align*}
and the conclusion follows by rearranging.
\end{proof}

Let us introduce the double-dot-product
$$F : G = f_\alpha^i g_\alpha^i$$ whenever $F = (f_\alpha^i)$ and  $G= (g_\alpha^i)$

Let us record a simple inequality.

\begin{lem}\label{lem:con}
Let $\partial M$ be a convex $C^2$-graph in $B_\eta(a) \subset \cN(\partial M)$ for some $a\in\partial M$, and suppose that all principal curvatures are bounded by 
$\kappa$. If $\eta\kappa\leq\frac{1}{2}$, then for any $v\in W^{1,1}(B;B_\eta(a)\cap\overline M)$ for a ball $B$, we have 
\begin{equation}\label{eq:DPi-Dnu-re}
D(\Pi\circ v): D(\nu\circ v) \geq 0\quad \text{a.e.\ in }B.
\end{equation} 

As a consequence, we also have 
\begin{equation}\label{eq:con}
|D(\Pi\circ v)|^2 + |D(\rho\circ v)|^2 \leq |Dv|^2\quad \text{a.e.\ in }B.
\end{equation}
\end{lem}

\begin{proof}
By \cite[eq. 14.99]{GT}, $\nabla \nu \leq 2\kappa I$ in $B_\eta(a)$ if $\eta\kappa \leq \frac{1}{2}$. By the assumption on the convexity of $\partial M$, $\nabla\nu = \Hess\rho \geq 0$ in $B_\eta(a)$. Thus, by \eqref{eq:DPi-Dnu}, $(\nabla\Pi)(\nabla\nu) \geq 0$ in $B_\eta(a)$, if and only if $I - \rho \nabla\nu \geq 0$ in $B_\eta(a)$. The latter inequality holds as $\rho < \eta$ in $B_\eta(a)$ and $2\eta\kappa \leq 1$. Then \eqref{eq:DPi-Dnu-re} follows from the chain rule. 

Now \eqref{eq:con} is a direct consequence of \eqref{eq:DPi-Dnu-re}.
In fact, by \eqref{eq:decompu}, a direct computation yields
\begin{equation}\label{Dv-expansion}
\begin{aligned}
 &   |Dv|^2 - |D(\Pi\circ v)|^2 - |D(\rho\circ v)|^2\\
& =  (\rho\circ v)^2 |D(\nu\circ v)|^2 + 2(\rho\circ v)D(\Pi\circ v): D(\nu\circ v) \geq 0,
\end{aligned}
\end{equation}

and the last inequality follows from $\rho\circ v \geq 0$ a.e.\ (as $|\Omega\setminus v^{-1}(\overline M)| = 0$) and \eqref{eq:DPi-Dnu-re}.
\end{proof}

Given a  measurable function $f\colon \R^n\to[0,\infty)$, let $\cM f$ denote its \textit{maximal function}:
$$
\cM f(x_0) := \sup_{r> 0} \left[ \frac{1}{r^n}\int_{B_r(x_0)} f\,dx \right]. 
$$
We refer the reader to \cite{Stein} for a rather comprehensive account of the properties of the maximal function. Here we will require the following result:

\begin{thm}[Maximal theorem]\label{thm:max} Let $f \in L^p(\R^n)$ be a nonnegative function, for some $p\in[1,\infty]$.
\begin{enumerate}[(i)]
    \item\label{it:stronlp} (Strong $(p,p)$-inequality) If $p\in (1,\infty],$ then $$
    \|\cM f\|_{L^\infty(\R^n)} \leq c\| f \|_{L^p(\R^n)},
    $$
    for some $c\equiv c(n,p)$.
    \item\label{it:weak1} (Weak $(1,1)$-inequality) If $p = 1$, then 
    $$
    | \{ \cM f > t \}| \leq \frac{c}{t} \int_{\R^n} f\,dx,
    $$
    for some $c\equiv c(n)$. 
\end{enumerate}
\end{thm}

 Finally, we have the following Vitali covering lemma, see \cite{W}. 

 \begin{lem}\label{lem:CZ-pre}
 Let $F$ and $G$ be measurable sets in the unit ball $B_1$, and let $\delta > 0$ be given. Assume the following:  
 \begin{enumerate}[(i)]
 \item $|F| \leq \delta|B_1|$.
 \item For every $x_0\in B_1$ and $r\in(0,1]$, if $|F\cap B_r(x_0)| > \delta |B_r(x_0)|$, then $B_r(x_0)\cap B_1\subset G$. 
 \end{enumerate}
 Then $|F| \leq 10^n\delta |G|$. 
 \end{lem}


\section{Basic properties of  constraint maps} \label{sec:basicanalysis}

Here we collect some basic variational inequalities satisfied by minimizing constraint maps. These (in)equalities were already proved in \cite{D} for the Dirichlet energy $\cE_0$. The proofs for the Alt--Caffarelli energy $\cE_q$ are essentially the same, but we include them for completeness. 

\begin{rem}
Here we shall call $w \in W^{1,2}(B)$ a weak solution to 
$$
\Delta w \geq f,
$$
for some $f \in L^1(B)$, if 
$$
- \int_B Dw \cdot D\vp \,dx \geq \int_B f \vp\,dx,
$$
for all $\vp\in W^{1,2}\cap L^\infty(B)$ with $\vp\geq 0$. 
\end{rem}

We begin with the inequality satisfied by the distance of $u$ to the obstacle. We remark that a similar statement for the case $q = 0$ can be found in \cite[Lemma 2.2]{D}, and our proof is essentially the same. However, the argument in \cite[Lemma 2.2]{D} is more complicated as it covers the manifold setting. This brought us to present a simpler proof in the Euclidean setting for the benefit of readers. 

\begin{lem}[Essentially the same as {\cite[Lemma 2.2]{D}}]\label{lem:subharm}
Let $u\in W^{1,2}(\Omega;\overline M)$ be a minimizing constraint map for $\cE_q$, $q\geq 0$. Then in any ball $B\subset u^{-1}(\cN(\partial M))$, $\rho\circ u \in W^{1,2}(B)$ and 
\begin{equation}\label{eq:w-pde}
\Delta (\rho\circ u) \geq \Hess\rho_u(Du,Du)\chi_{u^{-1}(M)}
\end{equation}
in the weak sense. 
\end{lem}

\begin{proof}
Let $\beta \in C^\infty([0,\infty))$ be such that $\beta' \geq 0$, $\beta = 0$ in $[0,\frac{1}{2})$ and $\beta = 1$ in $[1,\infty)$. Let $\delta > 0$ be given, and define $\beta_\delta(t) := \beta(\delta^{-1}t)$, and let $0\leq \eta\in C_c^\infty(B)$ be arbitrary. Let $\e > 0$ be arbitrarily small constant such that $2\e\|\eta\|_{L^\infty(B)}<\delta$, and set 
$$
u_\e := u - \e \eta \beta_\delta(\rho\circ u)\nu\circ u. 
$$
Then $\supp(u_\e - u)\subset\supp\eta\subset B\subset\Omega$. Moreover, since 
$$
\rho\circ u_\e = \rho\circ u - \e\eta\beta_\delta(\rho\circ u) \geq \frac{1}{2}\rho\circ u\quad\text{in }\Omega,
$$ 
we have $u_\e(\Omega)\subset\overline M$. Thus, $u_\e \in W^{1,2}(\Omega;\overline M)$ and it is an admissible competitor for $u$. We also observe that $u_\e^{-1}(\partial M)\subset u^{-1}(\partial M)$, since $\rho\circ u_\e \leq \rho\circ u$. Therefore, the minimality of $u$ yields 
$$
\begin{aligned}
0 &\leq \frac{1}{\e} \int_\Omega (|Du_\e|^2 - |Du|^2)\,dx \\
&=  -2\int_B D(\eta \beta_\delta(\rho\circ u)\nu\circ u) : Du \,dx + O(\e) \\
&= -2\int_B [\eta \beta_\delta'(\rho\circ u)|D(\rho\circ u)|^2 +\beta_\delta(\rho\circ u) D (\rho\circ u)\cdot  D \eta]\,dx \\
&\quad -2 \int_B \eta \beta_\delta(\rho\circ u)\Hess\rho_u(Du,Du)]\,dx + O(\e) \\
& \leq -2 \int_B \beta_\delta(\rho\circ u) [D(\rho\circ u)\cdot D\eta + \eta \Hess\rho_u(Du,Du)]\,dx + O(\e)
\end{aligned}
$$ 
where in the derivation of the the third line  we used the identities $D_\alpha (\rho\circ u) =  D_\alpha u^i (\nu^i\circ u)$ and $\Hess\rho_u (Du,Du) = D(\nu\circ u) : Du$; the last line follows from $\beta' ,\eta \geq 0$. Thus, taking $\e\to 0$, while keeping $\delta$ fixed, we obtain   
$$
0\leq - \int_B \beta_\delta(\rho\circ u)[D(\rho\circ u)\cdot D\eta + \eta \Hess\rho_u(Du,Du) ]\,dx.
$$
Letting $\delta\to 0$, utilizing $D(\rho\circ u) = D(\rho\circ u)_+$ a.e., we finally arrive at
$$
- \int_B D(\rho\circ u)\cdot D\eta \,dx \geq \int_{B\cap\{\rho\circ u>0\}} \eta \Hess\rho_u(Du,Du) \, dx .
$$
As $\eta$ was an arbitrary nonnegative function in $C_c^\infty(B)$ and $\{\rho\circ u > 0\} = u^{-1}(M)$, our claim is proved. 
\end{proof} 

\begin{rem}
    In Lemma \ref{lem:min} below, we will derive the equation satisfied by $\rho\circ u$ away from the singular set $\Sigma(u)$.
\end{rem}

\begin{cor}\label{cor:subharm}
If the complement of $M$ is convex, then $\rho\circ u \in W^{1,2}(\Omega)$ and it is weakly subharmonic in $\Omega$. 
\end{cor}

\begin{proof}   
    The convexity assumption on 
   $M^c$ guarantees that 
    $\rho$, $\Pi$ and $\nu$ are smooth in $\overline M$ and the decomposition $\id = \Pi + \rho\nu$ holds everywhere in $\overline M$. Therefore, the above lemmas holds for all $u\in W^{1,2}(\Omega;\overline M)$ when $M^c$ is convex,
   and that $\Hess \rho_a(\xi,\xi)\geq 0$ for all $\xi\in \R^m$ and all $a\in\cN(\partial M)$, hence the conclusion follows from the lemma. 
\end{proof}

Next we derive the variational equality for $f\circ u$ whenever $\nabla f$ is a tangential vector-field to $\partial M$. Let us stress that the argument here is exactly the same as that of \cite[Lemma 2.1 (iii)]{D}, since the term $\chi_{u^{-1}(M)}$ involved in $\cE_q[u]$ does not play any role when the variation respects the geometry of $\partial M$. Here we present the proof only to keep our exposition self-contained. 

\begin{lem}[Essentially the same as {\cite[Lemma 2.1 (iii)]{D}}]\label{lem:projectedPDE}
Let $f \in C_c^2(\R^m)$ be such that $\nabla f\cdot \nu = 0$ on $\partial M$, and let $u\in W^{1,2}(\Omega;\overline M)$ be a minimizing constraint map for $\cE_q$, $q\geq 0$. Then $f \circ u$ is a weak solution to 
$$
\Delta (f\circ u) = \Hess f_u(Du,Du)\quad\text{in }\Omega.
$$
\end{lem}

\begin{proof}
As $\nabla f \in C_c^1(\R^m;\R^m)$ is a vector-field, we can generate a global flow $\Phi :(-\infty,\infty)\times \R^m\to\R^m$ corresponding to $\nabla f$. Let $\eta \in C_c^\infty(\Omega)$ be given and consider, for each $\e > 0$,
$$
u_\e := \Phi(\e \eta, u), 
$$
which is clearly an admissible map in $W^{1,2}(\Omega;\overline M)$. Now as $\nabla f\cdot \nu = 0$ on $\partial M$, $u_\e^{-1}(M) = u^{-1}(M)$. Thus, 
$$
\cE_0[u_\e] - \cE_0[u] =\cE_q [u_\e] - \cE_q[u] \geq 0,
$$
for every small $\e > 0$. 
Now by \cite[Lemma 2.1(i)]{D}, 
$$
\left.\frac{d}{d\e}\right|_{\e= 0} \cE_q[u_\e] = \int_\Omega Du : D(\eta \nabla f_u) \,dx.
$$ 
Therefore, combining this with the above inequality, and taking into account the fact that $\eta$ may change sign, we arrive at 
$$
\int_\Omega Du : D(\eta \nabla f_u)\,dx = 0.
$$
The rest of the proof now follows by the chain rule. 
\end{proof}

By taking $f$ to be a component of $\Pi$ we obtain, after some vector calculus (cf.\ \cite[(4.1)]{FKS}), the following:

\begin{cor}\label{cor:proj-pde2}
Let $u\in W^{1,2}(\Omega;\overline M)$ be a minimizing constraint map for $\cE_q$, $q\geq 0$. In any ball $B\subset u^{-1}(\cN(\partial M))$ we have
    \label{cor:proj-pde}
 \begin{equation}\label{eq:proj-pde}
 \Delta (\Pi\circ u) = -2D(\rho\circ u)\cdot D(\nu\circ u) + \Hess\Pi_u((Du)^\top,(Du)^\top).
 \end{equation}
\end{cor}
Recalling \eqref{eq:tangpart}, we see that the projected image $\Pi\circ u$ has, at least formally, one extra derivative compared to $\rho\circ u$.

For the rest of this section, we shall observe some useful regularity improvement under stronger assumptions on $u$ or its projected image $\Pi\circ u$.

\begin{lem}\label{lem:proj-W2p}
Let $p\in (1,\infty)$ and suppose that $u\in W^{1,2p}(B_2;\overline M)$ is a  minimizing constraint map for the functional $\cE_q$. Let $f \in C_c^2(\R^m)$ be a function such that $\nu \cdot \nabla f= 0$ on $\partial M$. Then $|D^2 (f\circ u)| \in L^p(B_1)$ with the estimate 
$$
\| D^2 (f\circ u) \|_{L^p(B_1)} \leq c \left[\osc_{B_2} (f\circ u) + \| Du \|_{L^{2p}(B_2)}^2\right],  
$$
where $c$ depends only on $n$, $m$ and the $C^2$-character of $f$. 
\end{lem}

\begin{rem}\label{rem:proj-W2p}
We emphasize that, in this lemma, the assumption $\nu\cdot\nabla f= 0$ on $\partial M$ is crucial. In the proof of Theorem \ref{thm:e-reg}, we will take $f$ to be a $C^2$-extension of $\Pi$ from $\cN(\partial M)$ to $\overline M$. 
\end{rem}

\begin{proof}[Proof of Lemma \ref{lem:proj-W2p}]
By Lemma \ref{lem:projectedPDE} we have 
$$
\Delta (f\circ u) = \Hess f_u(Du,Du)\quad\text{in }B_2,
$$
in the weak sense. Since $f \in C_c^2(\R^m)$, there is a constant $c_f$ such that $|\Hess f| \leq c_f$ in $\R^m$. Thus, we have $|\Hess f_u(Du,Du)| \leq c_f|Du|^2$ a.e.\ in $B_2$. Since we assume $|Du|^2 \in L^p(B_{2r}(x_0))$, with $p > 1$, the conclusion follows from the classical $L^p$-theory. 
\end{proof}

By making use of the $W^{2,p}$-regularity of $\Pi\circ u$, we are now in a position to interpret the distance component of our mapping as the optimization of a certain Alt-Caffarelli-type functional with an added lower-order variation. In fact, this concept originates from the expansion mentioned in \eqref{Dv-expansion}, suggesting that an appropriate functional to examine is the following
$$
\begin{aligned}
\hat\cE_q[w] &:= \int_B (|Dw|^2 +  |D(\nu\circ u)|^2 w^2 + 2(D(\Pi\circ u): D(\nu\circ u))w_+ + q^2\chi_{\{w >0\}})\,dx,
\end{aligned}
$$
where $w_+ := \max\{w,0\}$. This functional is well-defined for every (scalar) function $w\in W^{1,2}(B)$ whenever $u(B)\subset \cN(\partial M)$ and $\Pi\circ u \in C^{0,1}(B)$. Indeed, since $\nu= \nu\circ \Pi$ in $\cN(\partial M)$, it is enough to assume regularity on the projected image.

\begin{lem}\label{lem:min}
Let $u\in W^{1,2}({B};\overline{M})$ be a minimizing constraint map for the functional $\cE_q$, and assume that $u(B)\subset \cN(\partial M)$ for some ball $B$. Suppose that $\Pi\circ u\in C^{0,1}(B;\partial M)$. Then $\rho\circ u \in W^{1,2}({B})$ is a minimizer of $\hat\cE_q$. In particular, $\rho\circ u\in C_{\loc}^{0,1}({B})$, and 
$$
\Delta (\rho\circ u) - \Hess\rho_u(Du,Du)\chi_{u^{-1}(M)} = q \cH^{n-1}|_{\partial_{{\rm red}} u^{-1}(M)},
$$
in the weak sense\footnote{That is, for every test function $\vp\in C_c^\infty({B})$, $$-\int_{B} D(\rho\circ u)\cdot D\vp\,dx + \int_{{B}\cap u^{-1}(M)} \vp \Hess \rho_u(Du,Du)\,dx = \int_{{B}\cap\partial_{\rm red} u^{-1}(M)} q \vp \,d\cH^{n-1}.$$} in ${B}$.  
Moreover, for every ${B}'\Subset{B}$ with $\dist({B}',\partial{B}) > \eta$, 
$$
\| \rho\circ u \|_{C^{0,1}({B}')} \leq c_\eta,
$$
where $c_\eta$ depends only on $n$, $m$, $\eta$, the $C^2$-character of $\partial M$, and $\| D(\Pi\circ u)\|_{L^\infty({B})}$. 
\end{lem}

\begin{rem}
We shall see in Section \ref{sec:e-reg} that the assumptions $u(B)\subset \cN(\partial M)$ and $\Pi\circ u\in C^{0,1}(B;\partial M)$ will be fulfilled in any neighborhood of small energy. 
\end{rem}

\begin{proof}
Let $0\leq w\in (\rho\circ u) + W^{1,2}_0({B})$, and set $v := \Pi\circ u + w \nu\circ u$ in ${B}$. Since $w\geq 0$, we see that $v \in W^{1,2}({B};\overline M)$, and moreover $v-u\in W^{1,2}_0({B};\overline{M})$.  Therefore, by the minimality of $u$ for the functional $\cE_q$, we have $\cE_q[u]\leq \cE_q[v]$ in ${B}$. Thus, by a direct computation, 
\begin{equation*}\label{eq:scalarbernoulli}
\hat\cE_{q}[\rho\circ u] = \cE_q[u] - \int_{B}  |D(\Pi\circ u) |^2 \,dx 
\leq \cE_q[v] - \int_{B}  |D(\Pi\circ u) |^2\,dx = \hat\cE_{q}[w],
\end{equation*}
which proves the minimality of $\rho\circ u$ for the functional $\hat\cE_{q}$ in ${B}$. Note that in this calculation we used the fact that, as $|\nu\circ V|=1$, we have $|(\nu^i\circ V) D(\nu^i\circ V)|=0$ for any admissible map $V$; in particular, the cross-term
$ \langle w D(\nu\circ V), (\nu\circ V) \otimes Dw\rangle = (w D_\alpha(\nu^i\circ V))( (\nu^i\circ V)D_\alpha w)=0$
vanishes.

The conclusion now follows essentially from \cite{GS}. The only difference here is that the functional $\hat \cE_{q}[w]$ involves the semilinear term $w^2 |D(\nu\circ u)|^2$. However, this is a lower-order perturbation, as $|D(\nu\circ u)|\in L^\infty({B})$ by our assumption that $\Pi\circ u\in C^{0,1}(B;\partial M)$. Thus, we may repeat the proof of \cite[Theorem 2.13]{GS} without any major modification. As a result, we obtain that $\rho\circ u\in C_{\loc}^{0,1}({B})$, as well as the interior estimate shown in the statement, and that
$$
\Delta (\rho\circ u) - F\chi_{\{\rho\circ u > 0\}} = q \cH^{n-1}|_{\partial_{\rm red}\{\rho\circ u > 0\}}, 
$$
in the weak sense in ${B}$, where 
$$
\begin{aligned}
F &:= (|D(\nu\circ u)|^2(\rho\circ u) + D(\Pi\circ u) : D\nu_u)\\
& = D(\nu\circ u) : D((\Pi + \rho\nu)\circ u) = D(\nu\circ u) : Du = \Hess\rho_u(Du,Du). 
\end{aligned}
$$
Note also that $u^{-1}(M) = \{\rho\circ u>0\}$. This completes the proof.
\end{proof}

\begin{cor}\label{rmk:proj-W2BMO}
Under the same assumption in Lemma \ref{lem:min}, $D^2(\Pi\circ u) \in {BMO}_\textup{loc}({B})$ and for every ${B}'\Subset{B}$ with $\dist({B}',\partial{B}) > \eta$, 
$$
[ D^2(\Pi\circ u) ]_{{BMO}({B}')} \leq c_\eta,
$$
where $c_\eta$ depends only on $n$, $m$, $\eta$, the $C^3$-character of $\partial M$, and $\| D(\Pi\circ u)\|_{L^\infty({B})}$. 
\end{cor}

\begin{proof}
This is an immediate consequence of standard elliptic regularity theory, Corollary \ref{cor:proj-pde} (especially, Eq.~\eqref{eq:proj-pde}) and Lemma \ref{lem:min}, which altogether yields $|\Delta(\Pi\circ u)| \in L_\loc^\infty(B)$. The {\it a priori} estimate follows easily, so we omit the details.
\end{proof}


\section{Lipschitz approximation}\label{sec:apprx}

This section is devoted to the study of Lipschitz approximation of minimizing constraint maps of the Alt--Caffarelli functional, when the energy is small. 

\begin{prop}[Lipschitz approximation]
\label{prop:apprx}
Let $M\subset \R^m$ be a $C^2$-domain with uniformly bounded principal curvatures  $\kappa$, and $\Omega\subset \R^n$ be a bounded domain with $B_{4r}(x_0)\Subset \Omega$. For each $\eta \in (0,1)$ and $\Lambda > 0$, there is a positive constant $\e_{\eta,\Lambda}$, depending only on $n$, $m$, $\kappa$, $\eta$ and $\Lambda$, such that the following holds: for every $\e \in (0,\e_{\eta,\Lambda})$ and every $q\in [0,\Lambda\e]$, if $u\in W^{1,2}(\Omega;\overline M)$ is a minimizing constraint map of $\cE_q$ satisfying 
    $$
    (4r)^{2-n}\int_{B_{4r}(x_0)} |Du|^2\,dx \leq \e^2,
    $$
    then there is a map $h_{\e}\in W^{1,2}(B_{2r}(x_0);\R^m)$ with
    $$
    \| Dh_{\e} \|_{L^\infty(B_r(x_0))} \leq \frac{c_\Lambda\e}{r},
    $$
    such that 
$$  
(2r)^{2-n}\int_{B_{2r}(x_0)} | D(u - h_{\e})|^2\, dx \leq \e^2\eta^2,
$$
where $c_\Lambda$ depends only on $n$, $m$ and $\Lambda$.
\end{prop}

The proof of Proposition \ref{prop:apprx} is somewhat long, and therefore we   divide the argument  into several lemmas. Let us remark that by rescaling, it suffices to prove the lemma for $x_0 = 0$ and $r= 1$. Throughout this section, we shall write by $c$ a generic positive constant, which depends at most on $n$, $m$, $\kappa$ and $\Lambda$, and in particular, it will be independent on the running index $k\in\N$. 

Let us begin with the setting of our proof. Suppose towards a contradiction that, for each $k\in\N$, we can find a $C^2$-domain $M_k\subset \R^m$ with uniformly bounded principal curvatures (by $\kappa$),  and a minimizing constraint map $u_k\in W^{1,2}(B_4;\overline M_k)$ of $\cE_{q_k}$, with $q_k\leq \Lambda\e_k$, for some sequence $\e_k\to 0$ such that 
        \begin{equation}\label{eq:Duk-L2}
        \int_{B_4} |Du_k|^2\,dx \leq \e_k^2,
        \end{equation}
        but for any map $h \in W^{1,2}(B_4;\R^m)$, satisfying 
        \begin{equation}\label{eq:Dh-Linf}
        \| Dh \|_{L^\infty(B_1)} \leq c_\Lambda \e_k,
        \end{equation}
        we must have 
        \begin{equation}\label{eq:Duk-Dh-L2}
        \int_{B_2} |D(u_k - h)|^2\,dx \geq \e_k^2\eta_0,
        \end{equation}
        for some $\eta_0 > 0$, independent of $k$. In \eqref{eq:Dh-Linf} the constant  $c_\Lambda $  (depending on $n$, $m$ and $\Lambda$) relates to  the scalar case  of the Alt-Caffarelli energy, see \cite[Corollary 3.3]{AC}, and will be determined later in Lemma \ref{lem:apprx}.

        Since $M_k$ is a $C^2$-domain with bounded principal curvatures (by $\kappa$), there is a tubular neighborhood $\cN(\partial M_k)$, whose width $2\eta$ is at least $1/(c\kappa)$, such that the signed distance function $\rho_k$ to (and the nearest point projection $\Pi_k$ onto) $\partial M_k$ is of class $C^2$ (resp.\ $C^1$); more explicitly, both $\Hess\rho_k$ and $\nabla\Pi_k$ are uniformly bounded in $\cN(\partial M_k)$; see Section \ref{sec:prelim}. We shall write $\nu_k := \nabla \rho_k$ as in Section \ref{sec:prelim}. 

        Due to the assumption $q_k \in [0,\e_k\Lambda]$, we have  
        \begin{equation}\label{eq:Qk-Q}
        Q_k := \frac{q_k}{\e_k} \to Q \in [ 0,\Lambda],
        \end{equation}
        along a subsequence. We shall assume here without loss of generality that the convergence holds along the full sequence.

Now, without loss of generality, we may assume 
\begin{equation}\label{eq:uk-av}  \means{B_4} u_k \,dx = 0, \end{equation}
which follows by a shift $u_k - a_k$, $M_k - a_k$, where $a_k$ is the integral average of $u_k$ over $B_4$.

        From here the Poincar\'e inequality implies

        \begin{equation}\label{eq:uk-L2}
        \int_{B_4} |u_k|^2\,dx 
        \leq c \int_{B_4} |Du_k|^2\,dx  \leq c\e_k^2.
        \end{equation}
        By both \eqref{eq:Duk-L2} and \eqref{eq:uk-L2}, we find some $v\in W^{1,2}(B_4;\R^m)$ such that 
        \begin{equation}\label{eq:vk-v}
        v_k := \frac{1}{\e_k} u_k \rightharpoonup v \quad\text{ weakly in }W^{1,2}(B_4;\R^m),
        \end{equation}
        along a subsequence. By the Rellich compactness theorem, $v_k \to v$ strongly in $L^2(B_4;\R^m)$ and a.e.\ in $B_4$, along a further subsequence. We shall denote this convergent subsequence with the same subscript in order to reduce the notation.

        Although  $0$ may not lie on $\overline{M}_k$ and 
        the signed distance  $\rho_k (0) $ maybe negative, we may still conclude,  from \eqref{eq:uk-L2}, that 
        \begin{equation}\label{eq:rk-0}
        \liminf_{k\to\infty} \rho_k(0) \geq 0.
        \end{equation}

        It does not harm us to assume that $\rho_k(0)$ converges (in the extended real line) along the full sequence. Then we can also find some $\gamma \in [-\infty,\infty]$ such that 
        \begin{equation}\label{eq:rk-ek}
        \limsup_{k\to\infty} \frac{\rho_k(0)}{\e_k} = \gamma.
        \end{equation}
        At this stage, a simple observation is  
        the next  lemma.

        \begin{lem}\label{lem:gammainfty}
            In the setting above, $\gamma>-\infty.$
        \end{lem}

        \begin{proof}
        Let $\tilde \rho_k=\e_k^{-1}\rho_k$ be the signed distance to $N_k:= \e_k^{-1} M_k$. 
        If $\gamma=-\infty$ then we would have $\tilde \rho_{k'}(0)\to -\infty$ along a subsequence $k'\to \infty$. Under this hypothesis, there is $\ell_{k'}\to \infty$ such that $\overline N_{k'}\cap B_{\ell_{k'}}(0)=\emptyset$. 
        However, for a.e.\ $x$ in $B_4$ we have $\rho_{k'}(v_k(x))\to \rho(v(x))$. Thus we can find at least one point $x_0\in B_4$ such that  the sequence $(v_{k'}(x_0))_{k'}$ is bounded. On the other hand  $v_{k'}(x_0)\in \overline N_{k'}$, 
        implies $|v_{k'}(x_0)| \geq \ell_{k'} \to \infty$, and hence 
         a contradiction is reached.
        \end{proof}

        In fact, in the proof of Lemma \ref{lem:AC-min} below (see \eqref{eq:signgamma}) we will reach    the stronger conclusion $\gamma \geq 0$.
        
        We shall divide the rest of the proof into the cases (i) $\gamma \in \R$ and (ii) $\gamma = \infty$. Intuitively, the rescaled domains $\e_k^{-1} M_k$ converge to a half-space in case (i) and to the full space in case (ii). Thus, in the limit, one either sees a flat constraint in case (i) or no constraint whatsoever in case (ii), and accordingly the limit map $v$ is either a minimizer of the Alt--Caffarelli energy $\cE_Q$ with a flat constraint, or simply a harmonic map into $\R^m$. In either case, a contradiction is derived from the Lipschitz regularity of $v$.

        \begin{lem}\label{lem:AC-min}
        Suppose that $\gamma \in \R$ in \eqref{eq:rk-ek}, and let $v_k$, $v$ and $Q$ be as in \eqref{eq:vk-v} and \eqref{eq:Qk-Q}, respectively. Then there is a unit normal $\nu_0$ such that $v\in W^{1,2}(B_4;\overline H)$ with 
        \begin{equation}
            \label{eq:H}
            H := \{ a\in\R^m: a\cdot\nu_0 > -\gamma\},
        \end{equation}
        and $v$ is a minimizing constraint map of $\cE_Q$ in $B_4$.
        Moreover, $v_k \to v$ strongly in $W_\loc^{1,2}(B_4;\R^m)$. 
        \end{lem}

        \begin{proof}        
    
        Since the proof is rather long, we split it into five steps.

        \medskip 
        
        \textbf{Step 1} (geometric setup and estimates)\textbf{.} As $\e_k\to 0$, \eqref{eq:rk-ek} and $\gamma \in \R$ imply that $\rho_k(0) \to 0$. Let $a_k \in\partial M_k$ be a boundary point such that  
        \begin{equation}\label{eq:ak}
        |a_k| = |\rho_k(0)|.
        \end{equation}
        This, combined with \eqref{eq:rk-ek}, clearly yields (for large enough $k$)
        \begin{equation}
            \label{eq:boundak}
            |a_k|\leq 2|\gamma|\e_k.
        \end{equation}
        
        As $M_k$ is a $C^2$-domain with principal curvatures  bounded uniformly and independently of $k$, and $\cN(\partial M_k)$ is the tubular neighborhood of width $2\eta$, there is an open neighborhood $U_k\subset\R^m$ of the origin, a rotation $\cR_k$ that maps $(\nu_k)_{a_k}$ to $e_m$, and a $C^1$-diffeomorphism $\Phi_k \colon B_\eta(a_k) \to U_k$ such that \eqref{eq:Phi} and \eqref{eq:id-Phi} hold. 
        In particular, 
        \begin{equation}\label{eq:Rk-Phik}
        |\nabla^i ((\cR_k-\cR_k a_k) - \Phi_k)| = o(\e_k^{1-i})\quad\text{in }B_{\e_kR_k}(a_k),
        \end{equation}
        for $i\in\{0,1\}$, as $k\to\infty$, where 
        \begin{equation}\label{eq:Rk}
        R_k := (\e_k|\log\e_k|^2)^{-1/2} \to \infty. 
        \end{equation}
        Up to subsequences, $(\nu_k)_{a_k} \to \nu_0$ and $\cR_k\to \cR$, where $\nu_0$ is a unit vector and $\cR$ is a rotation that maps $\nu_0$ to $e_m$. 

        Let us now take
        \begin{equation}\label{eq:Rk-R}
        \tilde \e_k := |\nabla(\cR_k - \cR)| = o(1)
        \end{equation}
        where we see $\nabla(\cR_k-\cR)$ as an $m\times m$ matrix.
        Let us choose a sequence of radii $\tilde R_k \to \infty$, $\tilde R_k \leq R_k$ such that 
        
        \begin{equation}
            \label{eq:tek-tRk}
            B_{\e_k\tilde R_k}(0)\subset \Phi_k(B_{\e_k  R_k}(a_k))
        \end{equation} for all large $k$. We now claim that
        \begin{equation}
            \label{eq:Rk-Phik-inverses}
            |\nabla^i(\id + a_k- \Phi_k^{-1}\circ \cR)| = o(\e_k^{1-i}) \quad \text{in } B_{\e_k \tilde R_k}(0)
        \end{equation}
        for $i\in \{0,1\}.$ Let us just prove the case $i=0$, as the other one is identical. If $y\in B_{\e_k \tilde R_k}(0)$, then by \eqref{eq:tek-tRk} we have $y = \cR^{-1} \Phi_k(a)$ for some $a\in B_{\e_kR_k}(a_k)$. We then estimate, using \eqref{eq:Rk-Phik},
        \begin{align*}
        |(\id + a_k- \Phi_k^{-1}\circ \cR)y|
        & \leq |(\id - \cR_k^{-1}\circ \cR)y| + |(a_k+\cR_k^{-1}\circ \cR- \Phi_k^{-1}\circ \cR)y|\\
        & = |(\cR_k-\cR)y|+|\cR_k a_k + \Phi_k(a)-\cR_k a|\\
        & \leq \e_k \tilde \e_k \tilde R_k + o(\e_k) =  o(\e_k),
        \end{align*}
        as claimed.
        
        Observe that by \eqref{eq:tek-tRk}, \eqref{eq:Rk-Phik-inverses} and the Inverse Function Theorem we have 
        \begin{equation}
            \label{eq:phikinverse}
            |\nabla \Phi_k^{-1}| \leq c\quad \text{ in } B_{\e_k \tilde R_k}(0).
        \end{equation}
        for all $k$ large.

         Let us also note that the nearest point projection $\pi_k$ onto $\overline M_k$ is well-defined via \eqref{eq:pi} in $\cN(M_k):= \cN(\partial M_k)\cup M_k$, and satisfies \eqref{eq:id-pi}; especially, 
        \begin{equation}\label{eq:id-pk}
        |\nabla^i (\id - \pi_k)| = o (\e_k^{1-i})\quad\text{in }B_{\e_k R_k}(a_k)\cap (a_k + T_{a_k}\overline{M_k}),
        \end{equation}
        for $i\in\{0,1\}$, with $R_k$ as in \eqref{eq:Rk}. Note from \eqref{eq:tan} that since $a_k \in\partial M_k$, $T_{a_k}\overline{M_k}$ is the half-space orthogonal to $(\nu_k)_{a_k}$, whose boundary passes through the origin.

        \medskip
        \textbf{Step 2} (image of the limit map)\textbf{.} In this step we show that
        \begin{equation}\label{eq:v-W12}
        v\in W^{1,2}(B_4;\overline H), 
        \end{equation}
        where $v$ is the weak limit given by \eqref{eq:vk-v} and $H$ is as in \eqref{eq:H}.

        To see this, fix any large $R > 2|\gamma|$ and small $\theta,\delta > 0$. By \eqref{eq:rk-ek},  \eqref{eq:boundak} and the convergence $(\nu_k)_{a_k}\to \nu_0$, we have 
        \begin{equation}\label{eq:ekMk-H}
        \e_k^{-1}M_k \cap B_{2R}(0) \subset -\theta \nu_0 + H
        \end{equation}
        for all $k$ large. Since $v_k\to v$ strongly in $L^2(B_4;\R^m)$, Egorov's theorem yields a closed set $F_\delta\subset B_4$ with $|B_4\setminus F_\delta| \leq\delta$, such that $v_k\to v$ uniformly on $F_\delta$. Then $\{|v| \leq R\}\cap F_\delta \subset \{|v_k| \leq 2R\}\cap F_\delta$ for all $k$ sufficiently large. 
        Thus, it follows from $|B_4\setminus u_k^{-1}(\overline{M_k})| = 0$, \eqref{eq:ekMk-H} and $u_k = \e_kv_k$ that   
        \begin{equation*}\label{eq:vknu0}
        v_k \cdot \nu_0 > - \gamma -\theta\quad\text{in }\{|v| \leq R\}\cap F_\delta,
        \end{equation*}
        for all large $k$. Letting $k\to\infty$, we obtain 
        \begin{equation*}
        v \cdot \nu_0 \geq -\gamma - \theta\quad\text{in }\{|v| \leq R\}\cap F_\delta. 
        \end{equation*}
        Since $R > 2|\gamma|$ can be taken arbitrarily large, while $\theta,\delta > 0$ arbitrarily small, we arrive at 
        \begin{equation}\label{eq:vnu0}
        v\cdot \nu_0 \geq - \gamma\quad\text{a.e.\ in }B_4.
        \end{equation}
        Since $v\in W^{1,2}(B_4;\R^m)$, we have verified \eqref{eq:v-W12}. 

        A simple consequence of this is that 
        \begin{equation}
            \label{eq:signgamma}
            \gamma\geq 0.
        \end{equation} 
        Indeed, by \eqref{eq:uk-av} and the strong convergence $v_k\to v$ in $L^2(B_4)$, we see that $\int_{B_4} v\, dx =0$, implying 
        $0 = \int_{B_4} v \nu_0 \, dx \geq -\gamma |B_4|$. This shows $\gamma \geq 0$.

        \medskip
        \textbf{Step 3} (construction of a competitor)\textbf{.} 
        With the goal of proving that $v$ is a minimizer of $\cE_Q$ over the class $W^{1,2}(B_4;\overline H)$, where $Q$ is as in \eqref{eq:Qk-Q}, we let $r\in(0,4)$ be arbitrary, and let $w\in W^{1,2}(B_4;\overline H)$ be any map satisfying $\supp(v- w)\subset B_r$. In this step, we will construct a competitor for a truncated version of $w$.

        Let us consider the truncation (as in \cite{L}) of $\gamma\nu_0 + w$ given by 
        \begin{equation}\label{eq:wk}
        w_k : =  \frac{(\gamma \nu_0 + w) \tilde R_k}{\max\{|\gamma\nu_0 + w|,\tilde R_k\}},
        \end{equation}
        with $\tilde R_k$ as in \eqref{eq:tek-tRk}, and define also $\tilde u_k \in W^{1,2}(B_r;\overline M_k)$ by 
        \begin{equation}\label{eq:tuk}
        \tilde u_k := \Phi_k^{-1} \circ \cR\circ(\e_k w_k).
        \end{equation}
        By \eqref{eq:phikinverse} and the definition of $w_k$, which ensures that $\e_k w_k \in B_{\e_k \tilde R_k}(0)$ a.e., we have the estimate
        \begin{equation}\label{eq:Dtuk-L2}
        \int_{B_4\setminus B_r} |D\tilde u_k|^2\,dx \leq c\e_k^2 \int_{B_4\setminus B_r} |Dw_k|^2\,dx \leq c\e_k^2, 
        \end{equation}
        where in the last inequality we used that $|Dw_k| \leq |Dw|$ a.e., $\supp(v-w)\subset B_r$ and $\int_{B_4} |Dv|^2\,dx \leq 1$ by \eqref{eq:Duk-L2} and \eqref{eq:vk-v}. 

        Moreover, since $a_k\in\partial M_k$ satisfies \eqref{eq:ak}, we have $a_k + |\rho_k(0)| (\nu_k)_{a_k} = 0$. Since \eqref{eq:rk-ek} implies $|\rho_k(0)| = (\gamma + o(1))\e_k$, and we know that $(\nu_k)_{a_k} = \nu_0 + o(1)$, 
        we obtain 
        \begin{equation*}\label{eq:ak2}
        a_k + \e_k \gamma \nu_0 = - |\rho_k(0)| (\nu_k)_{a_k} + \e_k \gamma \nu_0 = o(\e_k). 
        \end{equation*}

        Thus, as $v_k\to w$ strongly in $L^2(B_4\setminus B_r;\R^m)$, $w_k \to (\gamma\nu_0 + w)$ strongly in $L^2(B_4;\R^m)$, $u_k = \e_k v_k$ and $w_k(B_4)\subset B_{\tilde R_k}(0)$, we can estimate 
        
        \begin{equation}\label{eq:tuk-uk-L2-re}
        \begin{aligned}
        &\int_{B_4\setminus B_r} |u_k - \tilde u_k|^2\,dx \\
        & \leq o(\e_k^2) + 2\e_k^2 \int_{B_4\setminus B_r} |w - \e_k^{-1}  \tilde u_k|^2\,dx \\
        & \leq o(\e_k^2) + 4\e_k^2 \int_{B_4\setminus B_r} |w_k - \gamma\nu_0 - \e_k^{-1}  \tilde u_k|^2\,dx \\
        & \leq o(\e_k^2) + 8 \int_{B_4\setminus B_r} |\e_k w_k + a_k -  \Phi_k^{-1}\circ\cR\circ (\e_k w_k) |^2\,dx \\
        & = o(\e_k^2), 
        \end{aligned}
        \end{equation}
        where in the last line we used \eqref{eq:Rk-Phik-inverses}.

        Collecting \eqref{eq:Duk-L2}, \eqref{eq:Dtuk-L2} and \eqref{eq:tuk-uk-L2-re}, we can find a radius $s\in (r,4)$, via the Fubini theorem, such that

        \begin{equation}\label{eq:luckhaus-asmp}
        \int_{\partial B_s} ( |Du_k|^2 + |D\tilde u_k|^2 )\,d\sigma \leq c\e_k^2,\quad \delta_k^2 :=\int_{\partial B_s} \frac{|u_k - \tilde u_k|^2}{\e_k^2}\, d\sigma \to 0. 
        \end{equation}
        Let us also choose  $\alpha\in(0,1)$ small such that \begin{equation}\label{eq:alpha} 
        \alpha \left( 1- \frac{n}{2}\right) + \frac{1}{4} > 0, 
        \end{equation}
        and set 
        \begin{equation}\label{eq:lk-app}
        \lambda_k = \delta_k^\alpha \to 0.
        \end{equation}
        Then we invoke \cite[Lemma 1]{L} for $\tilde u_k$ and $u_k$, (with $p =2$, $\lambda = \lambda_k$, $\e = \delta_k$ and $\beta = 3/4$ there) which yields a map $\vp_k \in W^{1,2}(B_4;\R^m)$ such that
        \begin{equation}\label{eq:vpk-app}
        \vp_k(x) = \begin{cases}
        \tilde u_k ((1-\lambda_k)^{-1}x), & \text{if }|x| \leq (1-\lambda_k)s, \\ 
        u_k (x), & \text{if } |x| \geq s.
        \end{cases}
        \end{equation}
        Moreover, $\varphi_k$ satisfies the estimates
        \begin{equation}\label{eq:Dvpk-L2-app}
        \int_{B_s\setminus B_{(1-\lambda_k)s}} |D\vp_k|^2\,dx \leq c\e_k^2\left( 1 + \frac{\delta_k^2}{\lambda_k^2}\right)\lambda_k \leq 2c\lambda_k\e_k^2, 
        \end{equation}
        and 
        \begin{equation}\label{eq:dist-vpk-app}
        \dist(\vp_k,\overline{M_k}) \leq c \e_k \delta_k^{\frac{1}{4}}\lambda_k^{1-\frac{n}{2}} \leq c \delta_k^{\alpha(1-\frac{n}{2}) + \frac{1}{4}} \e_k  \ \to 0 \quad \hbox{as } k \to \infty, \quad  \text{in }B_s,
        \end{equation}
        where we have used \eqref{eq:alpha}, to conclude that the estimate tends to zero.

        Let us further consider
        \begin{equation}\label{eq:tvpk}
        \phi_k := \frac{1}{\e_k}\pi_k\circ \vp_k.
        \end{equation}
        This will be our competitor for $v_k$.

        \medskip
        
        \textbf{Step 4} (minimality of the limit map)\textbf{.} 
        By the definition of $\pi_k$, we have $w_k \in W^{1,2}(B_4;\overline N_k)$, where $N_k := \e_k^{-1} M_k$. Moreover, $\supp(\phi_k - v_k) \subset B_s\Subset B_4$. Since $u_k = \e_k v_k$, $q_k = \e_k Q_k$, and $u_k$ is a minimizer of $\cE_{q_k}$ over the class $W^{1,2}(B_s;\overline{M_k})$, then also $v_k$ is a minimizer of $\cE_{Q_k}$ over $W^{1,2}(B_s;\overline{N_k})$. Thus
        \begin{equation}\label{eq:min-uk-tuk}
        \begin{aligned}
         \int_{B_s} \Big(|Dv_k|^2 + Q_k^2 \chi_{v_k^{-1}(N_k)} \Big)dx & \leq  \int_{B_s} \Big(|D\phi_k|^2 + Q_k^2 \chi_{\phi_k^{-1}(N_k)} \Big) dx. 
        \end{aligned}
        \end{equation}
        In view of \eqref{eq:vpk-app} and \eqref{eq:tvpk}, we observe that $\e_k \phi_k(x) = \vp_k (x) = \tilde u_k(\frac{x}{1-\lambda_k})$ for $x\in B_{(1-\lambda_k)s}$, so
        \begin{equation}\label{eq:vpk-re-app}
        \begin{aligned}
        &\int_{B_{(1-\lambda_k)s}}  \Big(|D\phi_k|^2 + Q_k^2 \chi_{\phi_k^{-1}(N_k)} \Big) dx = \frac{(1-\lambda_k)^n}{\e_k^2} \int_{B_s} \bigg( \frac{|D\tilde u_k|^2}{(1-\lambda_k)^2} + q_k^2 \chi_{\tilde u_k^{-1}(M_k)}\bigg)dx. 
        \end{aligned}
        \end{equation}

        Note that, by  \eqref{eq:luckhaus-asmp} and \eqref{eq:lk-app}, $(1-\lambda_k)^{n-2} = 1 + o(1)$ as $k\to\infty$. We now compute each of the two terms on the right-hand side of \eqref{eq:vpk-re-app}.          
        As $|D(\cR\circ (\e_kw_k))|  = \e_k |Dw_k|$ a.e., recalling the definition of $\tilde u_k$ from \eqref{eq:tuk} and applying \eqref{eq:Rk-Phik-inverses} we see that
        \begin{equation}\label{eq:Dutk-L2}
        \frac{(1-\lambda_k)^{n-2}}{\e_k^2}\int_{B_s} |D\tilde u_k|^2 \,dx = (1 + o(1))\int_{B_s} |Dw_k|^2\,dx. 
        \end{equation}
        For the second term in \eqref{eq:vpk-re-app}, recall the definitions of $\Phi_k$ from \eqref{eq:Phi} and of $\cR$, to see that
        \begin{align*}
            \tilde u_k (x) \in M_k & \iff
            \cR(\e_kw_k(x)) \in \{y  : y^m > 0\}\\ 
            & \iff w_k(x) \in \{a : a\cdot \nu_0 > 0 \} = \gamma \nu_0 + H.
        \end{align*}
        Thus we compute
        \begin{equation}\label{eq:chiutk}
        \frac{(1-\lambda_k)^n}{\e_k^2}\int_{B_s} q_k^2 \chi_{\tilde u_k^{-1}(M_k)}\,dx = (1 + o(1))\int_{B_s} Q_k^2\chi_{w_k^{-1}(\gamma\nu_0 + H)}\,dx. 
        \end{equation}
        Combining \eqref{eq:Dutk-L2} with \eqref{eq:chiutk} and returning to \eqref{eq:vpk-re-app}, we deduce from $|Dw_k| \leq |Dw|$ a.e., $w_k^{-1}(\gamma\nu_0 + H)\subset w^{-1}(H)$ and $Q_k\to Q$ that 
        \begin{equation}\label{eq:vpk-re-app2}
        \begin{aligned}
         \int_{B_{(1-\lambda_k)s}}  \Big(|D\phi_k|^2 + Q_k^2 \chi_{\phi_k^{-1}(N_k)} \Big) dx \leq  \int_{B_s} \Big(|Dw|^2 + Q^2\chi_{w^{-1}(H)}\Big) dx + o(1).
        \end{aligned}
        \end{equation}

        To complete the computation of the right-hand side of \eqref{eq:min-uk-tuk}, we also need to estimate the integral over the annulus $B_s\setminus B_{(1-\lambda_k)s}$. However, by \eqref{eq:Qk-Q}, \eqref{eq:id-pk}, \eqref{eq:luckhaus-asmp}, \eqref{eq:lk-app}, \eqref{eq:Dvpk-L2-app}, \eqref{eq:dist-vpk-app}, and that 
        $|B_s\setminus B_{(1-\lambda_k)s} |\leq c\lambda_k$, 
        we deduce that 
        \begin{equation}\label{eq:vpk-re-app3}
        \begin{aligned}
        \int_{B_s\setminus B_{(1-\lambda_k)s}} \Big (|D\phi_k|^2 + Q_k^2 \chi_{\phi_k^{-1}(H_k)}\Big) dx  & \leq \frac{c}{\e_k^2} \int_{B_s\setminus B_{(1-\lambda_k)s}} |(\nabla\pi_k)\circ \vp_k|^2|D\vp_k|^2 \,dx + c\lambda_k \\
        & \leq \frac{c}{\e_k^2} \int_{B_s\setminus B_{(1-\lambda_k)s}} |D\vp_k|^2 \,dx + c\lambda_k \leq c\lambda_k \to 0. 
        \end{aligned}
        \end{equation} 
        Combining \eqref{eq:vpk-re-app2} with \eqref{eq:vpk-re-app3} and returning to \eqref{eq:min-uk-tuk}, we proceed to get 
        \begin{equation}\label{eq:min-uk-tuk-re}
        \int_{B_s} \Big(|Dv_k|^2 + Q_k^2\chi_{v_k^{-1}(N_k)} \Big)dx \leq  \int_{B_s} \Big(|Dw|^2 + Q^2\chi_{w^{-1}(H)}\Big) dx + o(1).
        \end{equation}
        
        To complete the proof of minimality of $v$ it suffices to use weak lower semicontinuity: firstly, due to \eqref{eq:vk-v}, we obtain that
        \begin{equation}\label{eq:lsc-1}
        \int_{B_s} |Dv|^2 \leq \liminf_{k\to\infty} \int_{B_s} |Dv_k|^2 dx.
        \end{equation}
        Secondly, arguing as in the proof for \eqref{eq:v-W12} above, we also have 
        \begin{equation}\label{eq:lsc-2}
        \int_{B_s} \chi_{v^{-1}(H)} \,dx \leq \liminf_{k\to\infty} \int_{B_s} \chi_{v^{-1}(N_k)} \,dx.
        \end{equation}
        Consequently, combining \eqref{eq:Qk-Q}, \eqref{eq:min-uk-tuk-re}, \eqref{eq:lsc-1}  and \eqref{eq:lsc-2} yields 
        \begin{equation}\label{eq:min-v-w}
        \int_{B_s} \Big(|Dv|^2 + Q^2 \chi_{v^{-1}(H)}\Big)dx \leq  \int_{B_s} \Big(|Dw|^2 + Q^2\chi_{w^{-1}(H)}\Big) dx. 
        \end{equation}
        Recall that $w\in W^{1,2}(B_4;\overline H)$ was an arbitrary map for which $\supp(v - w)\subset B_r$, and $r\in (0,4)$ was also arbitrary (although $s \in (r,4)$ may have depended on $w$). Thus $v$ is a minimizing constraint map of $\cE_Q$ over $W^{1,2}(B_4;\overline H)$.

        \medskip
        \textbf{Step 5} (strong convergence)\textbf{.} 
        Finally, we prove the strong convergence of $Dv_k\to Dv$ in $L^2(B_r;\R^{mn})$, for every $r\in(0,4)$. Due to its weak convergence in $L^2(B_r;\R^{mn})$, and the lower semicontinuity \eqref{eq:lsc-1}, it suffices to prove that 
        $$
        \limsup_{k\to\infty} \int_{B_r} |Dv_k|^2\,dx \leq \int_{B_r} |Dv|^2\,dx.
        $$
        However, this follows by simply taking $w = v$ in \eqref{eq:min-uk-tuk-re} (since $v$ is a competitor of itself), and utilizing \eqref{eq:lsc-2} and \eqref{eq:Qk-Q}. This finishes the proof. 
\end{proof}

The next lemma treats case (ii) above, i.e.,  when $\gamma = \infty$. The proof is similar in spirit to that of Lemma \ref{lem:AC-min} but, in fact, it is much simpler. This is due to the fact that we no longer need to involve the boundary flattening map $\Phi$, as there are no constraints acting on the limit map. In particular, we can follow essentially the argument in \cite[Page 358]{L} and, for this reason, we shall only mention necessary changes, omitting tedious details.

\begin{lem}\label{lem:harm}
Suppose that $\gamma  = \infty$ in \eqref{eq:rk-ek}, and let $v_k$, $v$ and $Q$ be as in \eqref{eq:vk-v} and \eqref{eq:Qk-Q} respectively. Then $\Delta v = 0$ in $B_4$, and $v_k\to v$ strongly in $W_\loc^{1,2}(B_4;\R^m)$.
\end{lem}

\begin{proof} 
As we assume $\gamma = \infty$, this means that $\lim_{k\to \infty} \tilde \rho_k(0)=+\infty,$
where $\tilde \rho_k=\e_k^{-1}\rho_k$ is the signed distance to $N_k:= \e_k^{-1} M_k$. By definition, $B_{\tilde\rho_k(0)}(0)\subset N_k$, and thus, by \eqref{eq:uk-L2} and \eqref{eq:vk-v}, 
\begin{equation}\label{eq:vk-msr}
|B_4\cap v_k^{-1}(\partial N_k)| \leq \frac{1}{(\tilde\rho_k(0))^2} \int_{B_4\cap v_k^{-1}(\partial N_k)} |v_k|^2\,dx \leq \frac{c}{(\tilde\rho_k(0))^2} \to 0,
\end{equation}
which says that the volume of the coincidence set tends to zero.

Recall from \eqref{eq:vk-v} that $v\in W^{1,2}(B_4;\R^m)$. Let $r\in(0,4)$ be arbitrary, and let $w \in W^{1,2}(B_4;\R^m)$ be any map for which $\supp(v-w)\Subset B_r$. We will prove that $w$ has at least the same Dirichlet energy as $v$. 

As in \cite[Page 358]{L}, we define the truncation of $w$, 
\begin{equation}\label{eq:wk-re}
w_k := \frac{wR_k}{\max\{|w|,R_k\}},
\end{equation}
in place of \eqref{eq:wk}, where $R_k$ is as in \eqref{eq:Rk}. Then we set 
\begin{equation}\label{eq:tuk-re}
\tilde u_k := \pi_k \circ (\e_kw_k), 
\end{equation}
where $\pi_k$ is the nearest point projection onto $\overline M_k$, defined as in \eqref{eq:pi}. 

Now it is straightforward from \eqref{eq:Duk-L2}, \eqref{eq:uk-L2}, \eqref{eq:vk-v}, and $\supp(v-w)\Subset B_r$, to verify \eqref{eq:luckhaus-asmp} for some radius $s\in(r,4)$, with $\tilde u_k$ as in \eqref{eq:tuk-re}. Thus, we obtain again from \cite[Lemma 1]{L} the extension map $\vp_k$ defined as in \eqref{eq:vpk-app}, which as before satisfies the estimates \eqref{eq:Dvpk-L2-app} and \eqref{eq:dist-vpk-app}; we choose $\delta_k$ and $\lambda_k$ exactly as in \eqref{eq:luckhaus-asmp} and \eqref{eq:lk-app}. Hence, setting $\phi_k$ as in \eqref{eq:tvpk}, we also obtain \eqref{eq:min-uk-tuk}. Now by the minimality of $v_k$ in $W^{1,2}(B_s;\overline N_k)$, we obtain the following chain of inequalities:
\begin{equation}\label{eq:min-vk-phik-harm}
\begin{aligned}
\int_{B_s} |Dv_k|^2 \,dx &= \int_{B_s} \Big( |Dv_k|^2 + Q_k^2 \Big)dx - Q^2 |B_s| + o(1)  \\
& = \int_{B_s} \Big( |Dv_k|^2 + Q_k^2\chi_{v_k^{-1}(N_k)}\Big) dx - Q^2 |B_s| + o(1) \\
& \leq \int_{B_s} \Big( |D\phi_k|^2 + Q_k^2\chi_{\phi_k^{-1}(N_k)}\Big) dx - Q^2 |B_s| + o(1) \\
& \leq \int_{B_s} |D\phi_k|^2 + o(1);
\end{aligned}
\end{equation}
note that the first equality is due to \eqref{eq:Qk-Q}, the second equality is due to \eqref{eq:vk-msr}, and in the last line we used $|\phi_k^{-1}(N_k)\cap B_s| \leq |B_s|$, in addition to \eqref{eq:Qk-Q} again. By the definitions \eqref{eq:vpk-app}, \eqref{eq:tvpk} and \eqref{eq:wk-re},  we obtain, just as in \eqref{eq:vpk-re-app}--\eqref{eq:vpk-re-app2}, that
\begin{equation}\label{eq:vpk-re-harm}
\begin{aligned}
\int_{B_{(1-\lambda_k)s}} |D\phi_k|^2 \,dx & = \frac{(1-\lambda_k)^{n-2}}{\e_k^2}\int_{B_s} |D\tilde u_k|^2\,dx \\
 & = (1 + o(1))\int_{B_s} |Dw_k|^2 \,dx  \leq \int_{B_s} |Dw|^2\,dx + o(1).
\end{aligned}
\end{equation}
Moreover, by \eqref{eq:id-pk}, \eqref{eq:Dvpk-L2-app} and \eqref{eq:dist-vpk-app}, 
\begin{equation}\label{eq:vpk-re-harm2}
\begin{aligned}
\int_{B_s\setminus B_{(1-\lambda_k)s}} |D\phi_k|^2\,dx &\leq \frac{c}{\e_k^2}\int_{B_s\setminus B_{(1-\lambda_k)s}} |D\vp_k|^2\,dx \leq c\lambda_k\to 0. 
\end{aligned}
\end{equation}
Combining \eqref{eq:vpk-re-harm} and \eqref{eq:vpk-re-harm2}, we may return to \eqref{eq:min-vk-phik-harm} to obtain 
\begin{equation}\label{eq:min-vk-w}
\int_{B_s} |Dv_k|^2 \,dx \leq \int_{B_s} |Dw|^2 \,dx + o(1).
\end{equation}
Finally, by the lower semicontinuity of the Dirichlet energy, we arrive at 
\begin{equation}\label{eq:min-v-w-harm}
\int_{B_s} |Dv|^2 \,dx \leq \int_{B_s} |Dw|^2\,dx. 
\end{equation}
As $w\in W^{1,2}(B_4;\R^m)$ was an arbitrary map with $\supp(v-w)\Subset B_r$, and $r\in(0,4)$ was also arbitrary, we deduce that $v$ minimizes the Dirichlet energy locally in $B_4$, without any constraint; thus $\Delta v = 0$ in $B_4$. 

The strong convergence of $Dv_k\to Dv$ in $L_\loc^2(B_4;\R^{mn})$ can also be proved by putting $w = v$ in \eqref{eq:min-vk-w}, exactly as in Lemma \ref{lem:AC-min}. 
\end{proof}

We now prove an interior gradient  estimate on the limit map. 
In case (ii), i.e.,  when  $\gamma = \infty$, the interior gradient estimate for the limit map $v$ is straightforward, by Lemma \ref{lem:harm}. The next lemma gives the corresponding estimate in case (i), i.e.,  when $\gamma \in \R$. We exploit the fact that the image constraint is given by a half-space, as shown in Lemma \ref{lem:AC-min}, as this decouples completely the normal and tangential components of the maps and thus makes our problem scalar. 

\begin{lem}\label{lem:apprx}
In either case of Lemmas \ref{lem:AC-min} and \ref{lem:harm}, the map $v$ verifies  
$$
\| Dv \|_{L^\infty(B_1)} \leq c (1+\Lambda),
$$
for some constant $c$ depending only on $n$ and $m$. 
\end{lem}

\begin{proof}
It suffices to treat the case where $v$ is as in the statement of Lemma \ref{lem:AC-min}. Let $H$ and $Q$ be as in the lemma. Let $\cR$ be a rigid motion (that is, rotation and translation) that maps $H$ onto $G := \{ a : a \cdot e_m > 0\}$. Set $w := \cR\circ v$ and write $w \equiv (\hat w,w_m)$. Let us observe that $\Delta \hat w = 0$ in $B_4$, and that $w_m$ is a scalar minimizer of 
$$
\widehat{\cE_Q}[\vp] := \int_{B_4} \Big( |D\vp|^2 + Q^2\chi_{\{\vp > 0\}} \Big)dx,
$$
among all nonnegative function $\vp\in W^{1,2}(B_4)$ with $\supp(w_m - \vp)\Subset B_4$, since $\cE_Q [w] = \cE_0[\hat w] +\widehat{\cE_Q}[w_m]$, and $\cE_Q [v] = \cE_Q[w]$ (here we let the energy $\cE_Q$ systematically vary upon the target constraint, which is  a slight abuse of notation). 

By \eqref{eq:Duk-L2} and \eqref{eq:vk-v}, we have $\int_{B_4} |Dv|^2 \,dx \leq 1$. Since $|Dv| = |Dw|$ a.e.,
$$
\int_{B_4} |Dw|^2 \,dx \leq 1.
$$
However, as $|Dw|^2 = |D\hat w|^2 + |Dw_m|^2$, the interior gradient estimate for harmonic functions yields 
$$
\| D\hat w \|_{L^\infty(B_1)} \leq c_1\int_{B_4} |D\hat w|^2 \,dx \leq c_1, 
$$
for some $c_1$ depending only on $n$ and $m$. Next it follows from the interior Lipschitz estimate for the scalar minimizers of the Alt-Caffarelli energy, see \cite[Corollary 3.3]{AC}, that
$$
\| Dw_m \|_{L^\infty(B_1)} \leq c_2( 1 + Q) \leq c_2(1 + \Lambda), 
$$
where the second inequality is due to \eqref{eq:Qk-Q}, and $c_2$ is a constant depending only on $n$. Combining the above two estimates together and using again that $|Dv| = |Dw|$ a.e., we arrive at the desired conclusion.
\end{proof}

We are now ready to conclude the proof of Proposition \ref{prop:apprx}. 

\begin{proof}[Proof of Proposition \ref{prop:apprx}]

By Lemma \ref{lem:apprx}, $h_k := \e_k v \in W^{1,2}(B_4;\R^m)$ verifies  \eqref{eq:Dh-Linf}. However, by Lemmas \ref{lem:AC-min} and \ref{lem:harm}, $Dv_k\to Dv$ strongly in $L^2(B_2;\R^{mn})$, regardless of the value of $\gamma \in (-\infty,\infty]$ in \eqref{eq:rk-ek}; here we also used Lemma \ref{lem:gammainfty} to rule out $\gamma=-\infty$. Thus
$$
\int_{B_2} |D(u_k - h_k)|^2 \,dx =  \e_k^2 \int_{B_2} |D(v_k - v)|^2 \,dx <\e_k^2\eta_0,
$$
for all large $k$, which gives a contradiction against \eqref{eq:Duk-Dh-L2} with $h = h_k$, for $k$ large enough.
\end{proof}


\section{The {$\varepsilon$}-regularity theorem}\label{sec:e-reg}

The purpose of this section is to  prove the $\e$-regularity theorem, namely Theorem \ref{thm:e-reg}, for minimizing constraint maps for $\cE_q$.

For the rest of this section, when we write a constant, the subscript(s), if any, will be the parameter(s) on which the constant depends on, in addition to the obvious parameters $n$, $m$ and the $C^3$-characters of $\partial M$; e.g., $\e_\delta$ depends on $\delta$, in addition to $n$, $m$, and $\partial M$.

Due to the Lipschitz approximation (Proposition \ref{prop:apprx}) for maps with small energy, we obtain a measure estimate for the superlevel set of the maximal function of $|Du|^2$. Such a passage was already shown in the proof of \cite[Lemma 3.7]{BW}. Still, to keep our exposition self-contained (as much as possible), we include  the argument here.

\begin{lem}\label{lem:cpt}
For all $\delta>0$ there is a small positive constant $\e_\delta$  and a large constant $N$ (independent of $\delta$) such that, if $u\in W^{1,2}(B_4;\overline M)$ is a minimizing constraint map of $\cE_q$ and $E(u,0,4) \leq \e^2$, where $q\leq \e\leq \e_\delta$, then 
    $$
    |\{ \cM(|Du|^2\chi_{B_4}) > N^2\e^2\} \cap B_1| \leq \delta |B_1|.
    $$ 
\end{lem}

\begin{proof}
Fix $\eta \in (0,1)$ and  apply Proposition \ref{prop:apprx} (with $\Lambda = 1$ there, as we assume $q \in [0,\e]$): for $\e\in (0,\e_\eta)$, there is a map $h_\e\in W^{1,2}(B_3;\R^m)$ such that 
\begin{equation}\label{eq:harm-app}
\int_{B_3} | D(u - h_\e)|^2\, dx \leq \e^2\eta^2,
\end{equation}
and 
\begin{equation}\label{eq:Dhe-Linf}
\| Dh_\e \|_{L^\infty(B_2)} \leq N_0\e ,
\end{equation}
where $N_0$ depends only on $n$ and $m$. 

We claim that with a suitable choice $N \equiv N(n,m) > N_0$, 
\begin{equation}\label{eq:max-claim}
\{ \cM(|Du|^2 \chi_{B_4}) > N^2\e^2 \}\cap B_1 \subset \{ \cM(|D(u - h_\e)|^2\chi_{B_2}) > N_0^2\e^2 \}.
\end{equation}
To prove this claim, let $x_0\in B_1$ be such that \begin{equation}
    \label{eq:claimepsreg}
    \cM(|D(u - h_\e)|^2\chi_{B_2})(x_0) \leq N_0^2\e^2.
\end{equation} Then for any $r\in (0,2)$, we have from \eqref{eq:Dhe-Linf} and the triangle inequality that 
$$
\int_{B_r(x_0)} |Du|^2 \,dx \leq 2N_0^2 \e^2 r^n + 2\int_{B_r(x_0)} |Dh_\e|^2\,dx \leq 4N_0^2\e^2 r^n.
$$ 
For any $r \geq 2$, we derive from the assumption $E(u,0,4)\leq \e^2$,
$$
\int_{B_r(x_0)\cap B_4} |Du|^2 \,dx \leq \int_{B_4} |Du|^2 \,dx \leq 4^{n-2}\e^2 \leq 4^{n-2}\e^2 r^n,
$$
for some $C_0 \equiv C_0(n)$. Thus, we have proved that $\cM(|Du|^2\chi_{B_4})(x_0) \leq N^2\e^2$, for $N := \max\{2N_0,2^{n-2}\}$; note that $N$  depends on $n$ and $m$ only, as so does $N_0$. Since $x_0\in B_1$ was an arbitrary point satisfying \eqref{eq:claimepsreg}, the claim in \eqref{eq:max-claim} is proved. 

Finally, the weak-$L^1$ estimate from Theorem \ref{thm:max}(\ref{it:weak1}) along with \eqref{eq:harm-app} implies that
\begin{equation}\label{eq:weak11}
|\{ \cM(|D(u - h_\e)|^2\chi_{B_2} ) > N_0^2 \}| \leq \frac{C}{N_0^2}\int_{B_2} |D(u - h_\e)|^2\,dx \leq \frac{C\e^2\eta^2}{N_0^2},
\end{equation}
with $C\equiv C(n)$. Now given $\delta$ we may choose $\eta$ sufficiently small such that $C\eta^2 \leq N_0^2 \delta$, and then identify $\e_\delta$ with the small constant $\e_\eta$ chosen above. Combining \eqref{eq:max-claim} and \eqref{eq:weak11} then yields 
$$
|\{ \cM(|Du|^2\chi_{B_4} ) > N^2\e^2\}\cap B_1| \leq \delta |B_1|,
$$
and we arrive at the desired conclusion. 
\end{proof}

Now we are ready to prove a geometric decay of the measure of the set where $|Du|$ is large.

\begin{lem}\label{lem:CZ}
    Let $u\in W^{1,2}(B_8;\overline M)$ be a minimizing constraint map for $\cE_q$ with $q \in [0,\e]$, and suppose that $E(u,0,8) \leq \e^2$. Then for every $\delta\in(0,1)$, there exists a small constant $\bar\e_\delta$ and a large constant $\bar N$ (which is independent of $\delta$) such that the following holds: setting 
    $$
    A_k := \{ \cM( |Du|^2\chi_{B_4}) > \bar N^{2k}\e^2\}\cap B_1,
    $$
    if $\e\leq \bar\e_\delta$ and $A_0 \neq \emptyset$, then $|A_1|\leq \delta |B_1|$, and 
    $$
    |A_{k+1}| \leq 10^n\delta |A_k|, 
    \quad\,\forall k\in \N.
    $$
\end{lem}

\begin{proof}
The conclusion will follow from the Vitali covering lemma, Lemma \ref{lem:CZ-pre}, hence we verify its two conditions. Let $\e_\delta$ and $N$ be as in Lemma \ref{lem:cpt}. We shall choose $\bar\e_\delta$ and $\bar N$ in such a way that $\bar\e_\delta \leq \e_\delta$ and $\bar N \geq N$. Then the assertion that $|A_1| \leq \delta |B_1|$ follows immediately from Lemma \ref{lem:cpt}. Let us also remark that $A_{k+1}\subset A_k$ for every $k \in \N\cup\{0\}$.

We are now going to verify the second condition in Lemma \ref{lem:CZ-pre}: we want to show that $|A_{k+1}\cap B_r(x_0)|\leq \delta |B_r(x_0)|$, under the assumption that for some $x_0\in B_1$ and $r\in (0,1]$, we have $B_r(x_0)\cap B_1\setminus A_k \neq\emptyset$, i.e., we assume that there is a point $x_1\in B_r(x_0)\cap B_1$ such that 
\begin{equation}\label{eq:max1}
\cM(|Du|^2\chi_{B_4})(x_1) \leq \bar N^{2k} \e^2. 
\end{equation}
Note that $|x_1 - x_0| < r$. Hence, \eqref{eq:max1} implies that if $r \in (0,\frac{4}{5})$, then  
\begin{equation}\label{eq:max2}
\begin{aligned}
E(u,x_0,4r) & \leq \frac{4^{2-n}}{5^{2-n}} E(u,x_1,5r)\leq \frac{5^n}{4^{n-2}} r^2  \cM(|Du|^2\chi_{B_4})(x_1) \leq c_1 r^2 \bar N^{2k}\e^2,
\end{aligned}
\end{equation}
for some $c_1 \equiv c_1(n)$. Also since $B_4(x_0)\subset B_8$ and $r\in(0,1]$, it follows from the almost monotinicty lemma  (Lemma \ref{lem:monot} with $q \leq \e$) and the assumption $E(u,0,8) \leq \e^2$ that
\begin{equation}\label{eq:max3}
\begin{aligned}
E(u,x_0,4r) &\leq E(u,x_0,4) + c_2 \e^2 \leq \frac{4^{2-n}}{8^{2-n}} E(u,0,8) + c_2 \e^2 \leq c_3 \e^2,
\end{aligned}
\end{equation}
for some $c_i \equiv c_i(n) \geq 1$, for $i\in\{2,3\}$. Recall that $\e_\delta$ is chosen as in Lemma \ref{lem:cpt}. Now write 
\begin{equation}\label{eq:max4}
\e(r) := \min\{ \sqrt{c_1} r \bar N^k \e, \sqrt{c_3}\e\},
\end{equation}
and finally choose $\bar\e_\delta$ small such that 
\begin{equation}\label{eq:max5}
\sqrt{c_3}\bar\e_\delta \leq \e_\delta.
\end{equation}
For $r \in [\frac{4}{5},1]$ we can choose $\bar N > 1$ large (independent of $r$) such that $\e(r) = \sqrt{c_3}\e$.
Then combining \eqref{eq:max2}, \eqref{eq:max3}, \eqref{eq:max4} and \eqref{eq:max5} altogether yield, for $r\in (0,1]$,
\begin{equation}\label{eq:max6}
E(u_{x_0,r}, 0 ,4) = E(u,x_0,4r) \leq \e(r)^2 \leq \e_\delta^2,
\end{equation}
where $u_{x_0,r}(x)\equiv u(rx + x_0)$. By the scaling relation, we see that $u_{x_0,r}$ is a minimizing constraint map for $\cE_{q(r)}$, with
$$
q(r) : = rq \leq \e(r),
$$
where the last inequality follows directly from \eqref{eq:max4} because $q\leq \e$, $N> 1$ and $c_i>1$, $i\in\{1,3\}$. Therefore, Lemma \ref{lem:cpt} applies to $u_{x_0,r}$, which yields 
$$
|\{ \cM(|Du_{x_0,r}|^2\chi_{B_4}) > N^2\e(r)^2 \} \cap B_1| \leq \delta |B_1|.
$$
Rewriting the last estimate in terms of $u$, we arrive at 
\begin{equation}
    \label{eq:auxCZlem}
    |\{ r^2\cM(|D u|^2\chi_{B_{4r}(x_0)}) > N^2 \e(r)^2 \}\cap B_r(x_0)| \leq \delta |B_r(x_0)|.
\end{equation}

Finally, we are left with replacing the set in \eqref{eq:auxCZlem} with $A_{k+1}$. To do so, we observe that in fact we can replace $x_0$ and $r$ in \eqref{eq:max2} and \eqref{eq:max3} with any $x\in B_r(x_0)$ and any $s \in [r,5]$ (i.e., the argument holds as long as $B_s(x)\subset B_8$). Defining $\e(s)$ as in \eqref{eq:max4} with $s$ in place of $r$, we observe that 
$$
E(u,x,s) \leq \e(s)^2 \leq s^2 \left[\frac{\e(r)^2}{r^2}\right],\quad\forall s \in [r,5],\forall x\in B_r(x_0).
$$
Since $B_4\subset B_5(x)$ for any $x\in B_r(x_0)$, we conclude that
\begin{equation}\label{eq:max9}
\begin{aligned}
&\cM(|Du|^2\chi_{B_{4r}(x_0)})(x)\leq \frac{N^2\e(r)^2}{r^2} \\ & \Longrightarrow \cM(|Du|^2\chi_{B_4})(x) \leq \frac{N^2\e(r)^2}{r^2}, \quad\forall x\in B_r(x_0),
\end{aligned}
\end{equation}
as can be seen by estimating $s^{-n}\int_{B_s(x)} |Du|^2 \chi_{B_4}\,  dx = s^2 E(u,x,s)$ separately for $s\leq r$, $r
\leq s\leq 4$ and $s\geq 4$.
By \eqref{eq:auxCZlem} and \eqref{eq:max9}, we obtain 
\begin{equation}\label{eq:max10}
|\{r^2 \cM(|Du|^2\chi_{B_4}) > N^2\e(r)^2\}\cap B_r(x_0)| \leq \delta|B_r(x_0)|.
\end{equation}

At this point, we choose $\bar N$ large such that $\bar N > c_1 N$. Then from the choice \eqref{eq:max4} of $\e(r)$, we see that 
$$
\frac{N^2\e(r)^2}{r^2} \leq c_1\bar N^{2k} N^2 \e^2 \leq \bar N^{2(k+1)}\e^2,
$$
implying $A_{k+1} \subset \{ (4r\sqrt n)^2 \cM(|Du|^2\chi_{B_4}) > N^2 \e^2 \}$. By \eqref{eq:max10}, we arrive at 
$$
|A_{k+1}\cap B_r(x_0)| \leq \delta |B_r(x_0)|.
$$
Since we started from the assumption that $B_r(x_0)\cap B_1\setminus A_k\neq\emptyset$, we have just proved the contraposition of the second requirement for the Vitali covering lemma (Lemma \ref{lem:CZ-pre}), and the conclusion follows. 
\end{proof}

Having the above two lemmas at our disposal, it is now easy to prove a suboptimal $\e$-regularity result, asserting that minimizing constraint maps are in $W^{1,p}$, for all $p<\infty$, near points of small energy. Again we remark that the assertion itself works in a much general setting beyond the minimizing maps for the Alt-Caffarelli functional.

\begin{lem}\label{lem:e-reg}
Let $p \in (1,\infty)$ be given. There is a small constant $\e_p > 0$ and a large constant $c_p > 1$ such that the following holds: if $u\in W^{1,2}(B_8;\overline M)$ is a minimizing constraint map for $\cE_q$ with $q\in[0,\e]$ and $E(u,0,8) \leq\e^2$ for some $\e \in [0,\e_p)$, then $|Du| \in L^p(B_1)$ with 
$$
\| Du\|_{L^p(B_1)} \leq c_p.
$$
\end{lem}

\begin{proof}
    Let $\delta$ be small enough so that $10^n\bar N^{2p} \delta \leq \frac{1}{2}$, and let $\e_p$ be small such that $\e_p \leq \bar\e_\delta$, with $\bar N$ and $\bar\e_\delta$ as in Lemma \ref{lem:CZ}. Then the lemma yields $|A_{k+1}| \leq 10^n\delta |A_k|$, for every $k \in \N$ and $|A_1| \leq \delta|B_1|$. Therefore, $|A_k| \leq (10^n\delta)^k|B_1|$ for every $k\in\N$, whence 
    $$
    \sum_{k=1}^\infty \bar N^{2pk} |A_k| \leq \sum_{k=1}^\infty 10^n \bar N^{2pk}\delta^k \leq |B_1|,
    $$
    due to our choice of $\delta$.  Hence, $\cM(|Du|^2\chi_{B_4}) \in L^p(B_1)$ and 
    $$
    \| \cM(|Du|^2\chi_{B_4}) \|_{L^p(B_1)} \leq c_p.
    $$ 
    The conclusion now follows from  Theorem \ref{thm:max}\eqref{it:stronlp}.
\end{proof}

It remains to upgrade the $\e$-regularity result from Lemma \ref{lem:e-reg}, yielding $W^{1,p}$-regularity, to a result yielding the desired $W^{1,\infty}$-regularity.
This will follow by a bootstrap argument, based on the observations of Section \ref{sec:basicanalysis}, where we showed that the tangential components (with respect to $\partial M$) of a minimizing constraint map have an additional derivative when compared to the distance part.

We will first prove a rescaled statement of Theorem \ref{thm:e-reg}. Here we exploit that the $0$-homogeneous rescaling, say of factor $r$, along which the normalized energy $E$ is invariant, decreases the weight $q$ to $r q$. 

\begin{prop}\label{prop:e-reg}
There exists a small constant $\bar\e > 0$ and a large constant $c> 1$, both depending only on $n$, $m$, and the $C^3$-character of $\partial M$, such that if $u\in W^{1,2}(B_4;\overline M)$ is a minimizing constraint map for $\cE_q$ with $q\in[0,\e]$, and $E(u,0,4) \leq \e^2$, for some $\e\in [0,\bar\e]$, then $|Du| \in L^\infty(B_1)$, and 
$$
\| Du \|_{L^\infty(B_1)} \leq c. 
$$
\end{prop}

\begin{proof}
Throughout this proof, we shall use $c$ to denote a generic, positive constant depending at most on $n$, $m$, and the $C^3$-character of $\partial M$; $c$ may vary at each occurrence. 

Let $\bar\e$ be the small energy constant as in Lemma \ref{lem:e-reg} corresponding to $p= 4n$. Assume that $\e\leq \bar \e$. Then by a standard rescaling and covering argument, one can deduce from Lemma \ref{lem:e-reg} that $|Du| \in L^{4n}(B_3)$ with 
\begin{equation}\label{eq:Du-Lp}
\| Du \|_{L^{4n}(B_3)} \leq c.
\end{equation}
The Sobolev embedding theorem yields $u\in C^{0,\frac{3}{4}}(B_3)$ and  
\begin{equation}\label{eq:osc-u}
[u]_{C^{0,\frac{3}{4}}(B_3)} \leq c.
\end{equation}
It follows from \eqref{eq:osc-u} that the non-coincidence set $u^{-1}(M)\cap B_3$ is open. Taking outer variations of the functional $\cE_q$ in this open set, one can easily observe that $\Delta u = 0$ in $u^{-1}(M)\cap B_3$. Hence, if $B_1\subset u^{-1}(M)$, then the desired $C^{0,1}$-estimate follows immediately from \eqref{eq:osc-u}. 

For the rest of the proof, let us assume that $B_1\cap u^{-1}(\partial M) \neq\emptyset$, and choose any $x_0\in B_1\cap u^{-1}(\partial M)$. Due to \eqref{eq:osc-u} and $\partial M$ being of class $C^3$, we can find a small constant $\delta > 0$, depending only on $n$, $m$ and the $C^3$-character of $\partial M$, such that $u ( B_{4\delta}(x_0)) \subset \cN(\partial M)$; we remark that $\delta$ is independent of the choice of $x_0$. Since $\Pi \in C^2(\cN(\partial M);\partial M)$ and $\nu\cdot\nabla \Pi = 0$ on $\partial M$, we deduce from Lemma \ref{lem:proj-W2p}, along with \eqref{eq:Du-Lp} and \eqref{eq:osc-u}, that $|D^2(\Pi\circ u)| \in L^{2n}(B_{2\delta}(x_0))$ and 
\begin{equation}\label{eq:proj-W2p}
\| D^2 (\Pi\circ u) \|_{L^{2n}(B_{2\delta}(x_0))} \leq c;
\end{equation}
here $c$ of course depends on $\delta$, but as $\delta$ depends on the fixed parameters, we suppress the dependence of $c$ on $\delta$. Again by the Sobolev embedding, $\Pi\circ u\in C^{0,1}(B_{2\delta}(x_0))$, and especially by \eqref{eq:proj-W2p},
\begin{equation}\label{eq:proj-C01}
[\Pi\circ u ]_{C^{1,\alpha }(B_{2\delta}(x_0))} \leq c \qquad (\hbox{some small } \alpha >0).
\end{equation}
Thus, we can use Lemma \ref{lem:min} to deduce that $\rho\circ u \in C^{0,1}(B_\delta(x_0))$ and 
\begin{equation}\label{eq:dist-C01}
\|\rho\circ u \|_{C^{0,1}(B_\delta(x_0))} \leq c.
\end{equation}
Combining \eqref{eq:proj-C01} with \eqref{eq:dist-C01}, and utilizing the identity \eqref{eq:decompu}, we obtain that $u \in C^{0,1}(B_\delta(x_0))$ and 
\begin{equation}\label{eq:Du-Linf2}
[u]_{C^{0,1}(B_\delta(x_0))} \leq c.  
\end{equation}

Recalling that $x_0$ was an arbitrary point in $B_1\cap u^{-1}(\partial M)$, and $\delta$ was chosen independently of $x_0$, we obtain 
\begin{equation}\label{eq:u-C01-fb}
[u]_{C^{0,1}(B_1\cap B_\delta(u^{-1}(\partial M)))} \leq c. 
\end{equation}
Since $\Delta u = 0$ in $B_3\cap u^{-1}(M)$, and $u = \Pi\circ u$ in $B_3\cap u^{-1}(\partial M)$, one can finally extend \eqref{eq:u-C01-fb} to all of $B_1$ by virtue of \eqref{eq:proj-C01} and \eqref{eq:osc-u} (along with the interior estimate for harmonic functions). We leave out the details of the last part as it is standard.  
\end{proof}

We are now ready to prove Theorem \ref{thm:e-reg}. 

\begin{proof}[Proof of Theorem \ref{thm:e-reg}]
Let $\bar\e$ be as in Proposition \ref{prop:e-reg}. Let us choose $\bar 
 r$ small enough such that $\bar r q \leq \bar\e$. By the assumption that $E(u,x_0,2r) \leq \e^2$, with $\e \leq \bar\e$ and $r \leq \bar r$, the rescaled map $u_{x_0,r}$ verifies $E(u_{x_0,r}, 0, 2) \leq \e^2$. Moreover, as $u$ is a minimizing constraint map of $\cE_q$, $u_{x_0,r}$ is a minimizing constraint map of $\cE_{rq}$, with $r q \leq \bar \e$. Now we may apply Proposition \ref{prop:e-reg} to $u_{x_0,r}$ (via standard covering argument) and then rescaling back to obtain the desired estimate. We omit the details. 
\end{proof}

As a corollary to Theorem \ref{thm:e-reg} we can characterize  the set of singularities $\Sigma(u)$, defined in \eqref{eq:Sing}, as the set where the normalized energy density of $u$ is positive:

\begin{cor}\label{cor:e-reg}
Let $u\in W^{1,2}(\Omega;\overline M)$ be a minimizing constraint map for the functional $\cE_q$, $q\geq 0$. Then $\Sigma(u) = \Omega\setminus \Theta_u^{-1}(0)$. 
\end{cor}

Let us also state another standard corollary from the above characterization of the singular points and the upper semicontinuity of the energy density, obtained in Corollary \ref{cor:semi}:

\begin{cor}\label{cor:e-reg2}
Let $u\in W^{1,2}(\Omega;\overline M)$ be a minimizing constraint map for the functional $\cE_q$, $q\geq 0$.  If $x_0\in\Sigma(u)$, then $0\in \Sigma(\phi)$ for any $\phi\in T_{x_0}u$,  i.e. $\phi$ has a singularity at the origin.
\end{cor}

Finally, we are able to completely justify \eqref{eq:main-sys} (cf.\ \cite{AC,GS}). For this we need to define the {\it reduced boundary} $\partial_{{\rm red}} E$ of a set $E$ with locally finite perimeter. Indeed, the reduced boundary of $E$
is  the set of points $x \in \R^n$ for which the density 
$\lim_{r \to 0} \frac1{|B(x,r)|} \int_{B(x,r)} \nu_E\, dx $
exists and has modulus one.

\begin{cor}\label{cor:main-sys}
Let $u\in W^{1,2}(\Omega;\overline M)$ be a minimizing constraint map for the functional $\cE_q$, $q\geq 0$. Then \eqref{eq:main-sys} holds in the sense of distributions: in particular, $q\, 
\mathcal{H}^{n-1}|_{\partial_{\rm red}u^{-1}(M)}$ is a distribution (hence a Radon measure) and $\partial u^{-1}(M)$ has locally finite perimeter in $\Omega$, satisfying $\cH^{n-1}(\Omega\cap (\partial u^{-1}(M)\setminus\partial_{\rm red} u^{-1}(M))) = 0$.
\end{cor}


\section{Regularity of the distance part}
\label{sec:contrho}

In this section, we study the continuity of $\rho\circ u$, which is  the distance part of the minimizing constraint map $u$. Let us present the main result of this section:

\begin{thm}\label{thm:dist-reg}
Let the complement of $M$ be a closed convex set in $\R^m$ with non-empty interior, and let $u \in W^{1,2}\cap L^\infty(\Omega;\overline M)$ be a minimizing constraint map for $\cE_q$, where $q\geq 0$. Then $\rho\circ u \in C(\Omega)$. 
\end{thm}

In our proof of the above assertion the convexity assumption enters crucially through the subharmonicity of $\rho \circ u$, obtained in Corollary \ref{cor:subharm}. We remark that the subharmonicity needs not be true for general target $M$.

To deduce everywhere continuity of $\rho\circ u$, we resort to the monotonicity of the normalized energy and the strong convergence of $0$-homogeneous blowups. Therefore, the distance part of any blowup is also subharmonic. However, as the distance part is also $0$-homogeneous, it must be a constant. Then our assertion follows from a case-study of the value of the constant. Let us now present the full argument. 

\begin{proof}[Proof of Theorem \ref{thm:dist-reg}]
Fix $x_0\in\Omega$. We aim at proving the following $(\e, \delta)$-argument: there is some $a \geq 0$ such that for every $\e > 0$, there exists some $\delta > 0$ such that $|\rho\circ u - a|\leq \e $ a.e.\ in $B_\delta(x_0)$.

For this purpose, let $\phi \in T_{x_0}u$ be given (as per Definition \ref{defn:tan}). According to Corollary \ref{cor:blowup}, $T_{x_0}u$ is non-empty, and $\phi$ is $0$-homogeneous, and so is $\rho\circ \phi$. However, based on Lemma \ref{lem:strong} and Corollary \ref{cor:subharm}, $\rho\circ \phi$ is subharmonic. Consequently, as shown in Lemma \ref{lem:0-hom}, $\rho\circ \phi$ is constant almost everywhere in $\R^n$. Now, we proceed to analyze the situation based on the value of this constant.

\begin{case}
$\rho\circ \phi = a$ a.e.\ for some $a >0$.
\end{case}

In this case, we have $\phi^{-1}(\partial M) =  \emptyset$ and hence 
$|\phi^{-1}(\partial M)| = 0$. Recalling from Corollary \ref{cor:blowup} that $\phi$ itself is a minimizing constraint map for $\cE_0$, it follows from \cite{D} that
$$
\int_B D\phi :Dv\,dy = - \int_{B\cap \phi^{-1}(\partial M)} v \cdot A_\phi(D\phi,D\phi)\,dy = 0, 
$$
for every ball $B\subset\R^n$ and $v\in C_c^\infty(B;\R^m)$. 
That is, $\phi$ is a vector-valued harmonic function in the entire space. Consequently, $\Sigma(\phi) = \emptyset $,  and by Corollary \ref{cor:e-reg2}, it follows that $x_0$ does not belong to $\Sigma(u)$. Now, we can set $a = \rho\circ u(x_0)$, and with this choice, the desired $(\epsilon, \delta)$-argument becomes effective.

\begin{case}
$\rho\circ \phi= 0$ a.e.
\end{case}

Since $\phi\in T_{x_0}u$, we can find a sequence $r_k\to 0$ such that $u_{x_0,r_k}\to \phi$ in $W^{1,2}(B_2;\R^m)$, and thus also $\rho\circ u_{x_0,r_k} \to \rho\circ \phi = 0$ strongly in $L^2(B_2)$. Utilizing the subharmonicity of $\rho \circ u_{x_0,r_k}$ from Lemma \ref{lem:subharm}, we obtain that in fact $\rho\circ u_{x_0,r_k}\to 0$  in $L^\infty(B_1)$. Rephrasing this observation in terms of $u$, we obtain
$$
\lim_{k\to\infty} \| \rho\circ u \|_{L^\infty(B_{r_k}(x_0))} = \lim_{k\to\infty} \| \rho\circ u_{x_0,r_k} \|_{L^\infty(B_1)} = 0,
$$
which proves the $(\e, \delta)$-argument with $a = 0$. Hence, our proof is complete. 
\end{proof}

Thanks to the above proposition, we obtain a well-defined free boundary, $\partial u^{-1}(M)\cap\Omega$, as well as a coincidence set $u^{-1}(\partial M)$, and non-coincidence set $u^{-1}(M)$, which is open. The following assertion is 
then an immediate consequence of the above. 

\begin{cor}\label{cor:dist-reg}
Under the setting of Theorem \ref{thm:dist-reg}, $u$ is harmonic in the open set $u^{-1}(M)$. In particular, $\Sigma(u)\subset u^{-1}(\partial M)$. 
\end{cor}

As another by-product, we observe that tangent maps at singular points are locally minimizing harmonic maps into $\partial M$.

\begin{cor}\label{cor:dist-reg2}
In the setting of Theorem \ref{thm:dist-reg}, if $x_0\in\Sigma(u)$ then every $\phi\in T_{x_0}u$ satisfies $\phi \in W_{\loc}^{1,2}(\R^n;\partial M)$ and it is a minimizing harmonic map into $\partial M$, i.e., for every bounded domain $\Omega\subset \R^n$, $\int_\Omega |D\phi|^2\,dy \leq \int_\Omega |Dv|^2\,dy$ for every $v\in \phi+W_0^{1,2}(\Omega;\partial M)$. 
\end{cor}

\begin{proof}

Let $x_0 \in \Sigma(u)$ be given. By Theorem \ref{thm:dist-reg} and Corollary \ref{cor:dist-reg},  $\rho\circ u(x_0) = 0$, and for every $\e > 0$, there is some $r > 0$ such that $0\leq \rho\circ u \leq \e$ in $B_r(x_0)$. Thus, with $u_{x_0,r}(y)= u(x_0+ry)$,  we have $\rho\circ u_{x_0,r} \to 0$, in the $L^\infty$, as was remarked in the proof  Theorem \ref{thm:dist-reg}.

Now choose an arbitrary map $\phi\in T_{x_0} u$. By the above observation, $\rho\circ \phi = 0$ a.e.\ in $B_1$. Then by the $0$-homogeneity of $\phi$, $\rho\circ \phi = 0$ a.e.\ in $\R^n$. Thus, $|\phi^{-1}(M)| = 0$, so $\phi\in W_{\loc}^{1,2}(\R^n;\partial M)$. On the other hand, by Corollary \ref{cor:blowup}, $\phi$ is a minimizer for the Dirichlet energy in every bounded domain 
$\Omega\subset\R^n$, among all competitors $v \in \phi+W_0^{1,2}(\Omega;\overline M)$. Clearly $W^{1,2}(\Omega;\partial M)\subset W^{1,2}(\Omega;\overline M)$, so the minimality continues to hold among competitors in $W^{1,2}(\Omega;\partial M)$. Thus $\phi$ is a minimizing harmonic map from every bounded domain $\Omega$ into $\partial M$. 
\end{proof}

 
\section{Density of the non-coincidence set near singularities}\label{sec:dense}

In this section, we prove that the non-coincidence set, denoted as $u^{-1}(M)$, exhibits null density at singular points of the map $u$. Instead of assuming only that the complement of the set $M$ is convex, we will impose the more robust condition \eqref{eq:boundshessrho}. Notably, our density estimate relies solely on the decay rate of $\rho\circ u$ at the specific chosen singular point, in addition to our universal parameters, mentioned earlier. As a consequence of this result, we will also establish that $\rho\circ u$ vanishes to an infinite degree at these singular points.

\begin{prop}\label{prop:dense}
Let $M$ be a $C^3$ domain satisfying  \eqref{eq:boundshessrho}, and let  $u \in W^{1,2}(B_2;\overline M)$ be a minimizing constraint map for $\cE_q$ with $q\in[0,q_0]$ such that $\int_{B_2} |Du|^2 \,dx \leq \Lambda$. Then for every $\e \in (0,1)$, there exists a constant $0<\eta =\eta(\e,\Lambda,q_0,n,m,\kappa)$ such that if $0\in\Sigma(u)$ and $\rho\circ u \leq \eta$ in $B_2$, then 
$$
|B_1 \cap u^{-1}(M)| \leq \e |B_1|.
$$
\end{prop}

\begin{proof}
Assume towards a contradiction that for some $\e_0 > 0$, there exists a sequence $u_k \in W^{1,2}(B_2;\overline M)$ of a minimizing constraint map associated with $q_k \in [0,q_0]$ such that $\int_{B_2} |Du_k|^2\,dx \leq \Lambda$, and $\rho\circ u_k \leq \eta_k\to 0$ in $B_2$, and $0\in\Sigma(u_k)$, but 
\begin{equation}\label{eq:ck}
\means{B_1} \bar\chi_k\,dx   \geq \e_0,\quad \text{where}\quad  \bar\chi_k := \chi_{u_k^{-1}(M)}.
\end{equation}

 Since $\{\bar\chi_k\}_{k=1}^\infty$ is bounded in $L^\infty(B_1)$, we have $\bar\chi_k\overset{\ast}{\rightharpoonup} \bar \chi$  in $L^\infty(B_1)$ along a subsequence, for some measurable function $\bar\chi$  with $0\leq \bar \chi\leq 1$ a.e.\ on $B_1$; to keep our exposition simple, we shall continue to use the same subscript for the subsequence. 

Since the complement of $M$ is uniformly convex (see \eqref{eq:boundshessrho}), $\partial M$ is necessarily bounded. Hence, it follows from our assumption $\rho\circ u_k \to 0$ in $B_2$ that $\{u_k\}_{k=1}^\infty$ is a bounded sequence in $L^\infty(B_2;\overline M)$. This combined with 
$$
\int_{B_2} |Du_k|^2\,dx \leq \Lambda
$$
shows that $\{u_k\}_{k=1}^\infty$ is a bounded sequence in $W^{1,2}(B_2;\overline M)$. Therefore, Lemma \ref{lem:strong} implies that $u_k\to u$ strongly in $W^{1,2}(B_1;\R^m)$ up to a subsequence, which for simplicity we do not relabel, as this does not affect our analysis below. Since $q_k\in [0,q_0]$ for all $k$, we may also assume without loss of generality that $q_k\to q$ for some $q\in[0,q_0]$. Then $u\in W^{1,2}(B_1;\overline M)$ and it is a minimizing constraint map for $\cE_q$. However, as we also assume that
$$
\rho\circ u_k \to 0 \quad\text{in }L^\infty(B_1),
$$
we have $u \in W^{1,2}(B_1;\partial M)$. Thus, regardless of the value of $q$, $u$ is a minimizing harmonic map into $\partial M$. Moreover,  since $0\in\Sigma(u_k)$ (by our assumption), Corollaries \ref{cor:semi} and \ref{cor:e-reg} imply that $0\in\Sigma(u)$. 

As a minimizing harmonic map, $u$ satisfies 
$$
\Delta u = A_u(Du,Du) \quad\text{in }B_1\setminus \Sigma(u),
$$
in the classical sense. Note that $B_1\setminus\Sigma(u)$ is a single open, connected component, as $\Sigma(u)$ being a set of Hausdorff dimension at most $n-3$ (see \cite{SU}) cannot enclose any open set. Thus, the unique continuation principle, cf.\ \cite{A}, along with $0\in\Sigma(u)$, implies that  
\begin{equation}\label{eq:Du}
|Du| > 0\quad\text{a.e.\ in }B_1.
\end{equation}

Let us write $H_k:= \Hess\rho_{u_k}(Du_k,Du_k)$ and $H := \Hess\rho_u(Du,Du)$. By Lemma \ref{lem:subharm}, we have 
\begin{equation}\label{eq:uk''-pde}
\int_{B_1} (\rho\circ u_k)\Delta \vp\,dx \geq \int_{B_1} H_k\vp \,dx
\end{equation}
for every $\vp\in C_c^\infty(B_1)$ with $\vp\geq 0$.  Since $u_k \to u$ strongly in $W^{1,2}(B_1;\R^m)$, we have that $u_k\to u$ a.e.\ in $B_1$ and $Du_k\to Du$ strongly in $L^2(B_1;\R^{mn})$. As $\partial M$ is of class $C^2$, we have  $\Hess\rho_{u_k}(\xi,\xi) \to \Hess \rho_u(\xi,\xi)$ a.e.\ in $B_1$ 
for all  $\xi\in\R^m$. Therefore an application of the Dominated Convergence Theorem, made possible by the upper bound in \eqref{eq:boundshessrho}, shows that 
$$
H_k \to H\quad\text{strongly in }L^1(B_1).
$$
 Recalling that $\bar\chi_k\overset{\ast}{\rightharpoonup} \bar\chi$ in $L^\infty(B_1)$, we deduce that 
$$
H_k\bar\chi_k \overset{\ast}{\rightharpoonup} H\bar\chi\quad\text{ in }L^\infty(B_1).
$$
This  along with $\rho\circ u_k \to 0$ in $L^\infty(B_1)$ now yields in \eqref{eq:uk''-pde}  that
$$
\begin{aligned}
\int_{B_1} H\bar\chi \vp \,dx & = \lim_{k\to\infty} \int_{B_1} H_k\bar\chi_k\vp\,dx  \leq \limsup_{k\to\infty} \int_{B_1} (\rho\circ u_k)\Delta\vp\,dx = 0.
\end{aligned}
$$
Since $\vp\geq 0$ is arbitrary in $C_c^\infty(B_1)$, we deduce that 
\begin{equation}\label{eq:Hbc}
H\bar\chi = 0\quad\text{a.e.\ in }B_1. 
\end{equation}
At this point, we use the strict convexity of the complement of $M$, which along with \eqref{eq:Du} implies that $H = \Hess\rho_u(Du,Du) > 0$ a.e.\ in $B_1$. Thus, by \eqref{eq:Hbc}, we must in fact have 
$$
\bar\chi = 0 \quad\text{a.e.\ in }B_1.
$$

Summarizing, we have proved that 
$\bar\chi_k\overset{\ast}{\rightharpoonup} 0$  in $L^\infty(B_1)$, which now by $\bar\chi_k\geq 0$ a.e.\ in $B_1$ implies $\bar\chi_k\to 0$ in $L^1(B_1)$. 
This is clearly a violation of \eqref{eq:ck}. This finishes the proof. 
\end{proof}

\begin{cor}\label{cor:dense2}
Under the assumptions of Proposition \ref{prop:dense}, for every modulus of continuity $\omega$, there exists another modulus of continuity $\sigma$, which is determined {\it a priori} by $n$, $m$, $q_0$, $\Lambda$, $\kappa$, and $\omega$, such that 
 if $0\in\Sigma(u)$ and $\rho\circ u \leq \omega(r)$ in $B_r$, for every $r\in(0,1)$, then 
$$
|B_r \cap u^{-1}(M)| \leq \sigma(r) |B_r|,
$$
for every $r\in(0,1)$. 
\end{cor}

\begin{proof}
Consider the rescaled map $u_r:= u_{0,r} \in W^{1,2}(B_2;\overline M)$, which is a minimizing constraint map for $\cE_{q_r}$, with $q_r := r q\in [0,q_0]$. By Lemma \ref{lem:monot}, 
$$
\int_{B_2} |Du_r|^2\,dy = (2r)^{2-n}\int_{B_{2r}} |Du|^2\,dx \leq \int_{B_2} |Du|^2\,dx + c_0q^2 \leq \Lambda + c_0 q^2,
$$
for every $r\in(0,1)$, with $c_0 \equiv c_0(n)$. Hence, we can apply Proposition \ref{prop:dense} (with $\Lambda + c_0q^2$ in place of $\Lambda$) to $u_r$ for every $r\in(0,\frac{1}{2})$. For each $r\in(0,\frac{1}{2})$, define $\sigma(r)$ as the supremum of $|B_1\cap v^{-1}(M)|/|B_1|$ over all minimizing constraint map $v\in W^{1,2}(B_2;\overline M)$ with $\int_{B_2}|Dv|^2\,dx \leq \Lambda + c_0q^2$, $0\in\Sigma(v)$ and $\rho\circ v\leq\omega(r)$ in $B_2$. By definition, $\sigma(r) \leq 1$ for every $r\in(0,\frac{1}{2})$ (so the supremum is well-defined). Since $\omega(r) \to 0$ as $r\to 0$, we may repeat the proof of Proposition \ref{prop:dense} and observe that $\sigma(r) \to 0$ as $r\to 0$. We can also rearrange $\sigma(r)$ to be nondecreasing in $r$. This finishes the proof. 
\end{proof}

As mentioned at the beginning of this section, due to the subsolution property, we can now deduce that $\rho\circ u$ vanishes to an infinite degree at every singular point of $u$. Recall that, according to Corollary \ref{cor:dist-reg}, such points are located within 
$u^{-1} (\partial M)$.
Furthermore, the speed of this vanishing depends solely on the degree and the modulus of continuity of $\rho \circ u$,   in addition to the universal parameters. Indeed, it is known that any integrable function, satisfying a weak reverse H\"older inequality,
vanishes of infinite order almost everywhere on the zero set of the  function; see, for instance, \cite[Theorem 14.5.1]{IM}. In our case we use the stronger property
of $\rho \circ u$ being subharmonic, and 
 present its short proof here.

\begin{cor}\label{cor:dense}
Under the assumptions of Proposition \ref{prop:dense}, for every $k= 1,2,\ldots$, there exist constants $\eta_k > 0$ and $C_k > 1$, all depending on $n$, $m$, $\Lambda$, $q_0$, $\kappa$, and $k$ only, such that if $\rho\circ u \leq \eta_k$ in $B_1$, then 
$$
\| \rho\circ u\|_{L^\infty(B_r)} \leq C_k r^k,
$$
for every $r\in(0,\frac{1}{2})$.
\end{cor}

\begin{proof}
    Since $\rho\circ u$ is weakly subharmonic in $B_1$ (see Lemma \ref{lem:subharm}), the local $L^\infty$ estimate along with $\{\rho\circ u> 0 \} = u^{-1}(M)$ yields, for every $r\in (0,\frac{1}{2})$,
   \begin{equation}\label{eq:sub}
    \| \rho\circ u \|_{L^\infty(B_r)}^2 \leq c_0\means{B_{2r}} (\rho\circ u)^2\,dx \leq c_0 \bigg[\frac{|B_{2r}\cap u^{-1}(M)|}{|B_{2r}|}\bigg]\| \rho\circ u \|_{L^\infty(B_{2r})}^2,
    \end{equation}
    where $c_0$ depends only on $n$; in the second inequality, we also use $\rho\circ u\geq 0$ everywhere. 
    Now, fix $k \in \mathbb N$ and consider a modulus of continuity $\omega(r) \leq \eta_k$ for every $r\in(0,1]$. Thanks to Corollary~\ref{cor:dense2}, there exists a modulus of continuity $\sigma$ such that
    $$
|B_{2r}\cap u^{-1}(M)| \leq \sigma(2r)|B_{2r}|
$$
for all $r \in (0,1/2)$. 
Choosing $r_k$ small enough so that $\sigma(2r)\leq \frac{2^{-2k}}{c_0}$ for all $r \in (0,r_k)$,
\eqref{eq:sub} implies
$$
\| \rho\circ u\|_{L^\infty(B_r)} \leq 2^{-k} \| \rho\circ u\|_{L^\infty(B_{2r})},
$$
for every $r\in (0,r_k)$. By a standard iteration on every dyadic scale, choosing $r = 2^{-\ell}r_k$, $\ell =1,2,\ldots$,
we deduce that
$$
\| \rho\circ u\|_{L^\infty(B_r)} \leq 2^k \left(\frac{r}{r_k}\right)^k \| \rho\circ u\|_{L^\infty(B_{r_k})} \leq C_k r^k
$$
for every $r\in (0,r_k)$.
Since the desired inequality for $r \in [r_k,1/2)$ is trivially satisfied by choosing a sufficiently large constant $C_k$, the result follows.
\end{proof}

 
\section{Nondegeneracy at the free boundary}\label{sec:ndeg}

In this section we prove the nondegeneracy of the distance part of minimizing constraint maps. The nondegeneracy is very different in the cases $q>0$ and $q=0$, and we only treat the former (but see \cite{FGKS} for the latter case). When $q>0$, in order to minimize the volume of the non-coincidence set, the distance part needs to grow at least linearly from the free boundary. Following the ideas of \cite{AC}, we construct a competitor that attaches completely to the boundary of the constraint in the interior, and show that it has smaller $\cE_q$ energy.

\begin{prop}\label{prop:ndeg}
Let the complement of $M$ be a convex, closed set with nonempty interior whose principal curvatures  are bounded uniformly by $\kappa$, and let $u\in W^{1,2}(\Omega;\overline M)$ be a minimizing constraint map for $\cE_q$ with $q > 0$. Then there exists $\e > 0$, depending only on $n$, $m$, $q$, and $\kappa$, such that 
$$
\sup_{B_{2r}(x_0)} (\rho\circ u) \leq \e r \quad\Longrightarrow \quad B_r(x_0)\subset u^{-1}(\partial M)
$$
whenever $B_{2r}(x_0)\subset \Omega$ and $r\in(0,1)$.
\end{prop}

\begin{proof}
By translation invariance, we shall assume that $x_0 = 0$. Fix $\e>0$, to be determined later, and let us assume, as in the requirement of the proposition, that 
\begin{equation}\label{eq:w-rk}
\sup_{B_{2r}} (\rho\circ u) \leq \e r. 
\end{equation}

Let $\psi_0$ be the truncated fundamental solution for the Laplacian such that $\psi_0 = 0$ on $\partial B_1$, $\psi_0 = 1$ on $\partial B_2$ and $\psi_0 > 0$ in $B_2\setminus \overline{B_1}$. We define the function
\begin{equation}\label{eq:psi}
\psi(x) := \e r \psi_0\left(\frac{x}{r}\right). 
\end{equation}
By construction, $\Delta\psi = 0$ in $B_{2r}\setminus \overline{B_r}$, $\psi = \e r$ on $\partial B_{2r}$ and $\psi = 0$ on $\partial B_{r}$. Let us continuously extend $\psi$ as $0$ inside $B_{r}$, and define
$$
v := \Pi\circ u + (\psi\wedge (\rho\circ u))\nu_u.
$$
Since both $\rho\circ u$ and $\psi$ are nonnegative, and $\psi \geq \rho\circ u$ on $\partial B_{2r}$, $v$ is an admissible competitor for $u$ in $B_{2r}$. Hence, the minimality of $u$ through the fact that $\psi =0$ in $B_{r}$ (i.e., $v(B_{r} )\subset \partial M$) as well as that $v^{-1}(M)\subset u^{-1}(M)$ (since $\rho\circ v \leq \rho \circ u$ in $B_{2r}$) implies that
\begin{equation}\label{eq:ndeg-var}
\int_{B_{r}} ( |Du|^2 - |D(\Pi\circ u)|^2 + q^2\chi_{u^{-1}(M)}) \,dx \leq I = -2I' + I'',
\end{equation}
where
\begin{equation}\label{eq:I-I'-I''}
\begin{aligned}
I &:=  \int_{B_{2r} \setminus B_{r} } ( |Dv|^2 - |Du|^2)\,dx, \\
I' &:= \int_{B_{2r} \setminus B_{r} } D(\Pi\circ u) : D((\rho\circ u - \psi)_+ (\nu\circ u) )\,dx, \\
I'' &:= \int_{B_{2r} \setminus B_{r} } ( |D((\psi\wedge (
\rho\circ u)) \nu\circ u)|^2 - |D (
(\rho\circ u)\nu\circ u)|^2 ) \,dx.
\end{aligned}
\end{equation}

Let us compute $I'$. As the complement of $M$ is assumed to be convex, whose principal curvatures  are uniformly bounded by $\kappa$, we can choose $\e$ sufficiently small, depending only on $m$ and $\kappa$, such that Lemma \ref{lem:con} and \eqref{eq:w-rk} yields 
\begin{equation}\label{eq:r-w-n}
D(\Pi\circ u): D(\nu\circ u) \geq 0 \quad\text{a.e.\ in }B_{2r} . 
\end{equation}
Now utilizing $\nu^k \partial_i \Pi^k = 0$, 
and  the chain rule  $D(\Pi^k\circ u) = \partial_i\Pi^k|_u Du^i$,
implies
$$(\nu^k\circ u)D(\Pi^k\circ u) = \nu^k|_u \partial_i\Pi^k|_u Du^i = 0.$$
In particular  $|(\nu^k \circ u)D (\Pi^k\circ u)| = 0$,
which  along with \eqref{eq:r-w-n} implies
\begin{equation}\label{eq:I'}
\begin{aligned}
I' & = \int_{B_{2r} \setminus B_{r} } (\rho\circ u - \psi)_+ D(\Pi\circ u) : D(\nu\circ u) \,dx \geq 0.
\end{aligned}
\end{equation}

Expanding the terms of $I''$ yields 
$$
\begin{aligned}
I'' & = -2 \int_{B_{2r} \setminus B_{r} } ( D \psi \cdot D(\rho\circ u - \psi)_+ + \psi (\rho\circ u - \psi)_+ |D(\nu\circ u)|^2 ) \,dx \\
&\quad - \int_{B_{2r} \setminus B_{r} } ( |D(\rho\circ u - \psi)_+|^2 + (\rho\circ u - \psi)_+^2 |D(\nu\circ u)|^2 )\,dx.
\end{aligned}
$$
Noting that $\psi > 0$ in $B_{2r} \setminus B_{r} $, we obtain 
\begin{equation}\label{eq:I''}
I '' \leq - 2  \int_{B_{2r} \setminus B_{r} } D \psi \cdot D(\rho\circ u - \psi)_+ \,dx \leq 2c \e \int_{\partial B_{r} } \rho\circ u\,d\sigma,
\end{equation}
where the second inequality follows from integration by parts, and $\partial\psi/\partial\vec n \leq c \e$ on $\partial B_{r} $, with $c \equiv c(n)$. Combining \eqref{eq:ndeg-var}, \eqref{eq:I'} and \eqref{eq:I''} yields 
\begin{equation}\label{eq:ndeg-var-re}
\begin{aligned}
\int_{B_r } (|Du|^2 - |D(\Pi\circ u)|^2 + q^2\chi_{u^{-1}(M)})\,dx & \leq 2c \e \int_{\partial B_r } \rho\circ u\,d\sigma. 
\end{aligned} 
\end{equation}

However, by Lemma \ref{lem:con} again (more precisely \eqref{eq:con}), we obtain from \eqref{eq:w-rk} (with $\e$ chosen small, depending only on $m$ and $\kappa$ again) that 
\begin{equation}\label{eq:Dw-DV-Du}
|D(\rho\circ u)|^2 + |D(\Pi\circ u)|^2 \leq |Du|^2 \quad\text{a.e.\ in }B_{2r} . 
\end{equation}
Combining \eqref{eq:ndeg-var-re} and \eqref{eq:Dw-DV-Du}, we obtain 
\begin{equation}\label{eq:sing}
\int_{B_{r} } ( |D(\rho\circ u)|^2 + q^2\chi_{u^{-1}(M)})\,dx \leq 2c_2\e\int_{\partial B_{r} } (\rho\circ u)\,d\sigma. 
\end{equation}
However, by \eqref{eq:w-rk} and \eqref{eq:W11} (as well as that $\rho\circ u\in W^{1,2}(B_r )$ and $\{\rho\circ u > 0\} = u^{-1}(M)$), 
$$
\int_{\partial B_{r} } \rho\circ u  \,d\sigma \leq n\e \int_{B_{r} } \chi_{u^{-1}(M)}\,dx + \int_{B_{r} } |D(\rho\circ u)|\,dx.
$$
Since $|D(\rho\circ u)| = 0$ a.e.\ in $u^{-1}(\partial M))$, Young's inequality yields 
$$
\int_{B_r } |D(\rho\circ u)|\,dx \leq \frac{1}{q}\int_{B_r } |D(\rho\circ u)|^2 \,dx + \frac{q}{4}\int_{B_r } \chi_{u^{-1}(M)}\,dx;
$$
here is indeed the very place we use the strict inequality $q > 0$. 
Combining the last two displayed formulas  yield 
\begin{equation}\label{eq:sing-re}
\int_{B_r } \rho\circ u\,d\sigma \leq  \frac{1}{q}\int_{B_r } \left(|D(\rho\circ u)|^2 + \left(\frac{q^2}{4} + n\e q\right) \chi_{u^{-1}(M)}\right)dx.
\end{equation}
Finally, combining \eqref{eq:sing} with \eqref{eq:sing-re}, we conclude that if  $4n\e<3q$ and $2c\e<q$, then
$$
\int_{\partial B_{r} } \rho\circ u \,d\sigma = 0.
$$
Since $\Delta (\rho\circ u) \geq 0$ in $B_{r} $ (see Lemma \ref{lem:subharm}), the maximum principle now implies that $\rho\circ u = 0$ in $B_{r} $, i.e., $B_{r} \subset u^{-1}(\partial M)$ as desired. 
\end{proof}


\section{Regularity near the non-coincidence set}\label{sec:reg}

In this section we prove Theorem \ref{thm:main}, which gives the universal Lipschitz regularity of minimizing constraint maps for the Alt--Caffarelli energy, in a universal neighborhood of their non-coincidence set. By a simple rescaling, we may take $\Omega=B_2$ and $\eta=1$. Precisely, we prove the following proposition.

\begin{prop}\label{prop:reg}
Let $M$ be as in Theorem \ref{thm:main}, and let $u \in W^{1,2}(B_2;\overline M)$ be a minimizing constraint map for the functional $\cE_q$ with $ q > 0$ such that  $\int_{B_2} |Du|^2 \,dx \leq \Lambda$. Then there are constants $c,\delta > 0$, depending only on $n$, $m$, $q$, $\Lambda$ and $\partial M$, such that 
$$
\| Du \|_{L^\infty(B_1\cap B_\delta(u^{-1}(M)))} \leq c.
$$
\end{prop}

 This result follows from the $\e$-regularity theorem and Proposition \ref{prop:ndeg} above, which yields a universal non-degeneracy at free boundary points. Indeed, the latter proposition yields a universal decay of the normalized energy in a universal neighborhood of the non-coincidence 
 set, stated in the next lemma.

\begin{lem}\label{lem:reg}
There exists a modulus of continuity $\omega$, which is determined by $n$, $m$, $q$, $\Lambda$ and $\partial M$ only, such that if $x_0\in B_1\cap \overline{u^{-1}(M)}$, then $E(u,x_0,r) \leq \omega(r)$ for every $r \in (0,1)$.  
\end{lem}

\begin{proof}
Assume that the assertion is not true, for some $q > 0$ fixed. Then given any small $\omega_0 > 0$, one would be able to find, for each $k=1,2,$, a minimizing constraint map $u_k\in W^{1,2}(B_4;\overline M)$ for the functional $\cE_q$, a point $x_k \in B_1\cap \overline{u_k^{-1}(M)}$ and a radius $r_k \to 0$ such that 
\begin{equation}\label{eq:Euk}
E(u_k,x_k,r_k) \geq \omega_0^2,\quad\text{but}\quad \int_{B_4} |Du_k|^2 \,dx \leq\Lambda.
\end{equation}
Since $x_k \in B_1\cap \overline{u_k^{-1}(M)}$, $B_r(x_k)\setminus u_k^{-1}(\partial M)\neq\emptyset$ for every $r\in(0,\frac{1}{2})$. Then as the contrapositive of Proposition \ref{prop:ndeg}, there exists a small constant $\e_0 > 0$, which can be determined {\it a priori} by $n$, $m$, $q$ and $\partial M$, such that 
\begin{equation}\label{eq:ruk-sup}
\sup_{B_{2r}(x_k)} (\rho\circ u_k) \geq \e_0 r,  
\end{equation}
for every $r\in(0,\frac{1}{2})$. 

In light of the almost monotonicity of the normalized energy (Lemma \ref{lem:monot}), one can find from \eqref{eq:Euk} a small constant $r_{\omega_0} \in(0,\frac{1}{2})$, depending only on $n$, $q$ and $\omega_0$, such that for all large $k$,
\begin{equation}\label{eq:Euk-re}
E(u_k,x_k,r) \geq \frac{1}{2}\omega_0^2,\quad\forall r \in (r_k,r_{\omega_0}).
\end{equation}
On the other hand, by the compactness of minimizing constraint maps (Lemma \ref{lem:strong}), there exists a minimizing constraint map $u\in W^{1,2}(B_2;\overline M)$ for the functional $\cE_q$ such that $u_k \to u$ strongly in $W^{1,2}(B_2;\R^m)$, up to a subsequence which we do not relabel. Extracting a further subsequence if necessary, we may also assume that $x_k\to x_0 \in \overline{B_1}$. Then since $r_k\to 0$, passing to the limit in \eqref{eq:Euk-re} yields
\begin{equation}\label{eq:Euk-re2}
E(u,x_0,r) \geq \frac{1}{2}\omega_0^2 >0,\quad\forall r \in(0,r_{\omega_0}). 
\end{equation}
Thus, by Corollaries \ref{cor:e-reg} and \ref{cor:dist-reg}, $x_0\in\Sigma(u)\subset u^{-1}(\partial M)$. Then by Theorem \ref{thm:dist-reg}, given $\eta > 0$, we can find a small $r_\eta > 0$ such that $\rho\circ u \leq \eta$ in $B_{r_\eta}(x_0)$. Now choosing $\eta$ sufficiently small as in Corollary \ref{cor:dense}, we obtain 
\begin{equation}\label{eq:ru-sup}
\| \rho\circ u \|_{L^\infty(B_r(x_0))} \leq C\left(\frac{r}{r_\eta}\right)^2,\quad\forall r \in (0,r_\eta). 
\end{equation}
However, since $u_k\to u$ a.e.\ in $B_2$, passing to the limit in \eqref{eq:ruk-sup} for every $r\in(0,r_\eta)$ yields along with \eqref{eq:ru-sup} that 
$$
\e_0 r \leq C\left( \frac{r}{r_\eta} \right)^2,
$$
which yields a contradiction for small $r$. Hence the proof is finished. 
\end{proof}

We are now ready to establish the universal Lipschitz regularity of the mapping. 

\begin{proof}[Proof of Proposition \ref{prop:reg}]
Let $\omega$ be as in Lemma \ref{lem:reg}, and fix a point $x_0\in B_1\cap \overline{u^{-1}(M)}$. Choose $\e$ and $r_q$ as in the $\e$-regularity theorem, Theorem \ref{thm:e-reg}; both parameters depend at most on $n$, $m$, $q$ and $\partial M$. Find $\delta\leq r_q$ such that $\omega(2\delta) \leq \e^2$; clearly $\delta$ depends only on $n$, $m$, $q$, $\Lambda$ and $\partial M$. Then by the aforementioned lemma $E(u,x_0,2\delta) \leq \e^2$, and thus the $\e$-regularity theorem yields $|Du|\in L^\infty(B_\delta(x_0))$ and 
$$
\| Du \|_{L^\infty(B_\delta(x_0))} \leq \frac{c_q}{\delta},
$$
where $c_q$ depends on $n$, $m$, $q$ and $\partial M$ only. We shall write $c := \delta^{-1}c_q$, which now depends further on $\Lambda$. Since this estimate holds uniformly for all $x_0\in B_1\cap \overline{u^{-1}(M)}$, our assertion is proved. 
\end{proof}


\section{Regularity of flat free boundaries}\label{sec:fb-reg}

 Knowing that minimizing constraint maps are Lipschitz regular in a universal neighborhood of their non-coincidence sets, we are ready to study the regularity of their free boundaries. In particular, in this section we will prove Theorem \ref{thm:fb-reg}.

 Let us begin with the basic regularity theory, which in view of Lemma \ref{lem:min} is more or less a corollary to \cite{AC, GS}. 

 \begin{lem}\label{lem:fb-reg}
 Let $\Lambda > 1$ and $q > 0$ be given, and  $M$ be as in Theorem \ref{thm:main}. There exist small constants $\e , \bar\gamma > 0$, both depending only on $n$, $m$, $q$, $\Lambda$ and $\partial M$, such that if $u\in W^{1,2}(B_4;\overline M)$ is a minimizing constraint map for the functional $\cE_q$ such that $\int_{B_4} |Du|^2\,dx \leq \Lambda$, and if 
\eqref{eq:flat} holds, then $B_1\cap \partial u^{-1}(M)$ is a $C^{1,\bar\gamma}$-hypersurface, whose $C^{1,\bar\gamma}$-character
 depends only on $n$, $m$, $q$, $\Lambda$ and the $C^3$-character of $\partial M$. 
 \end{lem}

 \begin{proof}
 By Proposition \ref{prop:reg} we know that, 
 for some $\delta := \delta(n,m,q,\Lambda,\partial M)$,  
 \begin{equation}\label{eq:Du-Linf}
 \| Du \|_{L^\infty(B_\delta(u^{-1}(M))\cap B_3)} \leq c_0,
 \end{equation}
 for some constant $c_0\equiv c_0(n,m,q,\Lambda,\partial M)$. Thus Lemma \ref{lem:min} is applicable, and shows that $\rho\circ u\in C^{0,1}(B_\delta(u^{-1}(M))\cap B_3)$ is a weak solution to 
 \begin{equation}\label{eq:dist-pde}
 \Delta (\rho\circ u) - \Hess\rho_u(Du,Du)\chi_{u^{-1}(M)} = q \cH^{n-1}|_{\partial_{\rm red}u^{-1}(M)}.
 \end{equation} 
 By Proposition \ref{prop:ndeg} and \eqref{eq:Du-Linf}, 
 \begin{equation}\label{eq:ndeg-re}
 \e_0 r \leq \means{\partial B_r(x_0)} \rho\circ u\,d\sigma \leq c_0 r,
 \end{equation}
 for some $\e_0\equiv \e_0(n,m,q,\partial M)$, for every $x_0\in \partial u^{-1}(M)\cap B_2$, and every $r\in(0,\delta)$. Note also that $\Hess\rho_u(Du,Du) \in C^{0,\gamma}(B_2)$ for every $\gamma\in (0,1)$, such that 
 \begin{equation}\label{eq:Hessrho-Ca}
 \| \Hess\rho_u(Du,Du) \|_{C^{0,\gamma}(B_2)} \leq c(\gamma),
 \end{equation}
 where $c(\gamma) \equiv c(n,m,q,\gamma,\Lambda,\partial M)$. To see this, $\Pi\circ u \in C^{1,\gamma}(B_2)$ by Lemma \ref{lem:proj-W2p} and the Sobolev embedding; then \eqref{eq:Hessrho-Ca} is an immediate consequence of identities \eqref{eq:identityhessrho} and \eqref{eq:tangpart}, as $\Hess\rho \in C^1(\overline M)$. 
 
 By \eqref{eq:boundshessrho}, \eqref{eq:dist-pde}, \eqref{eq:ndeg-re} and \eqref{eq:Hessrho-Ca}, the implication that the flatness gives $C^{1,\gamma}$ for the free boundary follows directly by \cite[Theorem 2.17]{GS} (see also \cite[Theorem 8.1]{AC}). 
 \end{proof}

However, achieving higher regularity of the free boundary, as well as the mapping around flat points, is no longer an incidental outcome of the scalar theory. This is primarily due to the presence of the term $\Hess\rho_u(Du, Du)$ in \eqref{eq:dist-pde}, where the regularity is affected by the Jacobian of the projected image, denoted as $D(\Pi\circ u)$; this term is vectorial. This issue was initially encountered by three of the authors in their recent work \cite{FKS}, which focuses on constraint maps of obstacle type. In this paper, we employ a bootstrap argument, identical to the one used in the earlier work, and we provide a brief overview of the argument.

 \begin{prop}\label{prop:fb-reg}
 Let $k\geq 2$ be an integer and $M$ be a domain of class $C^{k+1}$, and $u\in C^{0,1}(B;\overline M)$ a minimizing constraint map for the functional $\cE_q$, $q>0$. If $B\cap \partial u^{-1}(M)$ is a $C^1$-hypersurface, then for every $\sigma,\gamma\in(0,1)$, $\sigma B\cap \partial u^{-1}(M)$ is locally a $C^{k,\gamma}$-hypersurface, whose $C^{k,\gamma}$-character is determined only by $n$, $m$, $q$, $\gamma$, $\sigma$, $k$, $\diam B$, $\| Du \|_{L^\infty(B)}$, the $C^{k+1,\gamma}$-character of $\partial M$ and the $C^1$-character of $B\cap \partial u^{-1}(M)$. 
 \end{prop}

 \begin{proof}

Given the assumption that $u\in C^{0,1}(B)$ and that $B\cap \partial u^{-1}(M)$ is of class $C^1$, the conditions specified in Lemma \ref{lem:fb-reg} are met. As a result, $B\cap \partial u^{-1}(M)$ can be considered a locally $C^{1,\bar\gamma}$-hypersurface, where $\bar\gamma\in(0,1)$. The value of $\bar\gamma $ depends solely on parameters such as $n$, $m$, 
$ |Du |_{L^\infty(B)}$, and the $C^3$-character of $\partial M$. We will maintain this value of $\bar \gamma $  throughout the rest of the proof.

 By Lemma \ref{lem:min}, $\rho\circ u$ verifies 
 \begin{equation}\label{eq:dist-pde-re}
 \begin{cases}
 \Delta (\rho\circ u) = \Hess \rho_u(Du,Du) & \text{in } B\cap u^{-1}(M), \\
 \rho\circ u =0,\, |D(\rho\circ u)| = q & \text{on }B\cap \partial u^{-1}(M).
 \end{cases}
 \end{equation}
 The crucial observation is that $\Hess\rho_u(Du,Du) \in C_{\textup{loc}}^{0,\gamma}(B)$ for every $\gamma\in (0,1)$. To see this, we first note that $\rho\circ u\in C^{0,1}(B)$ by  assumption and $\Pi\circ u \in C_{\textup{loc}}^{1,\gamma}(B)$ by Lemma \ref{lem:proj-W2p} and the Sobolev embedding; then the regularity of the right-hand side in \eqref{eq:dist-pde-re} is an immediate consequence of identities \eqref{eq:identityhessrho} and \eqref{eq:tangpart}, as $\Hess\rho \in C^1(\cN(\partial M)\cap\overline M)$.
 
As $B\cap \partial u^{-1}(M)$ is of class $C^{1,\bar\gamma}$, the boundary regularity estimates for elliptic equations with Neumann boundary conditions now yield $\rho\circ u\in C_{\textup{loc}}^{2,\bar\gamma}(B\cap\overline{u^{-1}(M)})$. Thus, we can further employ the partial hodograph and Legendre transformation \cite[Theorem 2]{KN} to deduce that $B\cap \partial u^{-1}(M)$ is locally a $C^{2,\gamma}$-hypersurface, for every $\gamma\in(0,1)$. Applying the boundary regularity estimates again, we also deduce that $\rho\circ u\in C_{\textup{loc}}^{2,\gamma}(B\cap \overline {u^{-1}(M)})$ for every $\gamma\in(0,1)$. 

The conclusion now follows by a bootstrap argument as in \cite{FKS}. Indeed, if for some integer $1\leq \ell \leq k-2$ we have that
\begin{itemize}
    \item $\Pi\circ u\in C_{\textup{loc}}^{\ell,\gamma}(B\cap \overline{u^{-1}(M)})\cap C_{\textup{loc}}^{\ell,\gamma}(B\setminus u^{-1}(M)),$
    \item $\rho\circ u \in C_{\textup{loc}}^{\ell+1,\gamma}(B\cap\overline{u^{-1}(M)})$,
    \item $B\cap \partial u^{-1}(M)$ is locally a $C^{\ell+1,\gamma}$-hypersurface,
\end{itemize}
  then by applying the regularity theory for elliptic transmission problems \cite[Theorem 1.1]{Z} to \eqref{eq:proj-pde}, we deduce that in fact 
  $$\Pi\circ u \in C_{\textup{loc}}^{\ell+1,\gamma}(B\cap \overline{u^{-1}(M)})\cap C_{\textup{loc}}^{\ell+1,\gamma}(B\setminus u^{-1}(M)),$$
  thus gaining a derivative. We can then repeat the argument in the previous two paragraphs to conclude. See \cite{FKS} for further details.

 \end{proof}

 Finally, Theorem \ref{thm:fb-reg} now follows immediately from Lemma \ref{lem:fb-reg} and Proposition \ref{prop:fb-reg}. We shall omit the obvious details.


\section{Regularity of the image map}\label{sec:proj}

This final section concerns the regularity of the projected image $\Pi\circ u$ of minimizing constraint maps: we will prove Theorem \ref{thm:proj-reg}. In Remark \ref{rmk:proj-W2BMO} we already saw that the projected image of a Lipschitz constraint map has second derivatives in $BMO$.

Here we observe that this $BMO$-regularity can be sharpened to $L^\infty$ around regular free boundary points. This result follows by closer inspection of the proof of Proposition \ref{prop:fb-reg} 

\begin{lem}\label{lem:proj-C11}
Let $\partial M$ be of class $C^{3,\gamma}$, for some $\gamma\in(0,1)$, and $u\in C^{0,1}(B_4;\cN(\partial M)\cap\overline M)$ be a minimizing constraint map for the functional $\cE_q$. Suppose that $0\in \partial u^{-1}(M)$ and $B_4\cap \partial u^{-1}(M)$ is a $C^1$-hypersurface. Then $|D^2(\Pi\circ u)| \in L^\infty(B_1)$ and 
$$
\| D^2(\Pi\circ u) \|_{L^\infty(B_1)} \leq c,
$$
where $c$ depends only on $n$, $m$, $\gamma$, $\|Du\|_{L^\infty(B_4)}$, the $C^1$-character of $B_4\cap \partial u^{-1}(M)$ and the $C^{3,\gamma}$-character of $\partial M$. Moreover, $D^2(\Pi\circ u) \in C^{0,\gamma'}(B_1\cap\overline{u^{-1}(M)})\cap C^{0,\gamma'}(B_1\setminus u^{-1}(M))$, for every $\gamma'\in(0,\gamma)$. 
\end{lem}

\begin{proof}
First of all, with $\partial M$ being of class $C^3$ only, the proof of Proposition \ref{prop:fb-reg} yields that $\rho\circ u \in C^{2,\alpha}(B_2\cap \overline{u^{-1}(M)})$ and $B_2\cap \partial u^{-1}(M)$ is of class $C^{2,\alpha}$, for every $\alpha\in(0,1)$. Next, we notice that the first step of the bootstrap argument works under the (weaker) assumption that $\partial M$ is of class $C^{3,\gamma}$, instead of $C^4$, as this still ensures that $\Hess\Pi \in C^{0,\gamma}(\cN(\partial M))$. In particular, the regularity theory for transmission problems \cite{Z}, as in the above proof, yields that $\Pi\circ u \in C^{2,\alpha\gamma}(B_1\cap \overline{u^{-1}(M)})\cap C^{2,\alpha\gamma}(B_1\setminus u^{-1}(M))$, so in particular $|D^2(\Pi\circ u)| \in L^\infty(B_1)$. 

As in the proof of Proposition \ref{prop:fb-reg}, we shall not present the details on the dependence of the $L^\infty$-norm. 
\end{proof}

Observe that our proof above shows 
$C^{2,\alpha \gamma}$-regularity for projected map from each side,
but only $C^{1,1}$ across the the free boundary.

Finally, note that Theorem \ref{thm:proj-reg} now follows immediately from Theorem \ref{thm:main}, Corollary \ref{rmk:proj-W2BMO} and Lemmas \ref{lem:fb-reg} and \ref{lem:proj-C11}, upon a suitable rescaling. 


 \appendix


 \section{Proof for the almost monotonicity formula}\label{sec:tech}


In this appendix we present, for   readers' convenience  the proofs of  technical lemmas in Section \ref{sec:cpt}.

\begin{proof}[Proof of Lemma \ref{lem:monot}]

Fix any pair $(r,s)$ of constants with $0<r<s< \dist(x_0,\partial \Omega)$. For any test function $\zeta\in C_c^\infty(B_s(x_0);\R^n)$, we set $Q_t(x) := x + t \zeta(x)$. Then for all $|t|$ sufficiently small, $Q_t$ is a diffeomorphism of $B_s(x_0)$ onto itself. Set 
$$
U_t := u\circ Q_t^{-1} \subset W^{1,2}(B_R(x_0);\overline M),
$$
that agrees with $u$ in a neighborhood of $\partial B_R(x_0)$, for all small $|t|$. Thus the family $\{U_t\}$ consists of admissible competitors for $u$, with 
 $U_0 = u$ on $B_s(x_0)$. Utilizing $\det DQ_t = 1 + t \ddiv \zeta + o(t)$ and $DQ_t^{-1} = I - t D\zeta + o(t)$ as $t\to 0$, we get by the change of variables formula
$$
\begin{aligned}
&\int_{B_s(x_0)} ( |DU_t|^2 + q^2 \chi_{U_t^{-1}(M)} )\,dx \\
&= \int_{B_s(x_0)} (|Du \,DQ_t^{-1}\circ Q_t|^2 + q^2\chi_{u^{-1}(M)})\det DQ_t\,dy \\
&= \int_{B_s(x_0)} ((|Du|^2 + q^2\chi_{u^{-1}(M)})( 1 + t\ddiv\zeta)  - 2t D_\alpha u^i D_\beta u^i D_\alpha \zeta^\beta)\,dy + o(t).
\end{aligned}
$$
Hence, by the minimality of $u$, 
$$
\begin{aligned}
0 & \leq \cE_q[U_t] - \cE_q[u] \leq t \int_{B_s(x_0)} ((|Du|^2 + q^2\chi_{u^{-1}(M)})\ddiv\zeta - 2 D_\alpha u^i D_\beta u^i D_\alpha \zeta^\beta)\,dy + o(t).
\end{aligned}
$$
Since we allow $t$ to change sign, the linear term must vanish, which yields 
$$
\int_{B_s(x_0)}( (|Du|^2 + q^2 \chi_{u^{-1}(M)}) \ddiv \zeta - 2 D_\alpha u^i D_\beta u^i D_\alpha\zeta^\beta)\,dy = 0.
$$
Now choosing $\zeta(x) := (x-x_0)\eta(|x-x_0|)$, with a smooth cutoff function $\eta \in C_c^\infty([0,s])$ satisfying $\eta(0) = 1$ and $0\leq \eta\leq 1$, we find 
$$
\begin{aligned}
&\int_{B_s(x_0)} [ (n-2)|Du|^2 + nq^2 \chi_{u^{-1}(M)} ]\eta \,dx + \int_{B_s(x_0)} \left[|Du|^2 + q^2\chi_{u^{-1}(M)} - 2 \left|\frac{\partial u}{\partial 
N}\right|^2\right] R \eta'  dx  =0,
\end{aligned}
$$
where $R := |x-x_0|$, and $\partial/\partial N$ denotes the directional derivative in the radial direction $R^{-1}(x-x_0)$, and $\eta,\eta'$ are evaluated at $R$.  Selecting a sequence $\{\eta_j\}_{j=1}^\infty$ that approximates the characteristic function of $[0,s]$ and passing to the limit, we obtain
$$
\begin{aligned}
&\int_{B_s(x_0)} [(n-2) |Du|^2 + n q^2 \chi_{u^{-1}(M)}]\,dx  - \frac{1}{s}\frac{d}{ds} \int_{B_s(x_0)} \left[ |Du|^2 + q^2 \chi_{u^{-1}(M)}- 2\left| \frac{\partial u}{\partial N} \right|^2\right]dx = 0.
\end{aligned}
$$
Multiplying both sides by $s^{1-n}$ and rearranging terms, we arrive at 
$$
\begin{aligned}
& \frac{d}{d s} \left[s^{2-n}\int_{B_s(x_0)} [ |Du|^2 +  q^2 \chi_{u^{-1}(M)} ]\,dx\right] \\
& = 2s^{1-n} \int_{B_s(x_0)} q^2 \chi_{u^{-1}(M)}\,dx + 2\frac{d}{ds}\left[\int_{B_s(x_0)} R^{2-n} \left| \frac{\partial u}{\partial N}\right|^2 \,dx\right] ,
\end{aligned}
$$
for a.e.\ $s \in (0,\dist(x_0,\partial\Omega))$. Hence, for $0<r<s<\dist(x_0,\partial\Omega)$, upon integration we find that
\begin{align*}
\begin{split}
\label{eq:integratedenergy}
    & 2 \int_{B_s(x_0)\setminus B_r(x_0)} R^{2-n}\left| \frac{\partial u}{\partial N}\right|^2 dx \\
    & \leq s^{2-n} \int_{B_s(x_0)} [ |Du|^2 + q^2 \chi_{u^{-1}(M)} ]\,dx - r^{2-n} \int_{B_r(x_0)} [|Du|^2 + q^2\chi_{u^{-1}(M)}] \,dx  \\
    & \leq s^{2-n} \int_{B_s(x_0)} |Du|^2\,dx - r^{2-n} \int_{B_r(x_0)} |Du|^2 \,dx + \omega_n q^2 s^2,
\end{split}
\end{align*}
where $\omega_n$ is the volume of the $n$-dimensional ball. This finishes the proof of the almost monotonicity of the normalized energy. Moreover, if we have $s^{2-n} \int_{B_s(x_0)} |Du|^2\,dx = r^{2-n} \int_{B_r(x_0)} \,dx $ for some $r<s$ and $q = 0$, then $|\partial u/\partial N| =0$ a.e.\ in $B_s(x_0)\setminus B_r(x_0)$, which yields the $0$-homogeneity of $u$ about $x_0$ in the annulus.

\end{proof}

\section{Proof for the compactness lemma}\label{app:compact}

\begin{proof}[Proof of Lemma \ref{lem:strong}]

Let's consider $B_r(x_0)\Subset \Omega$. Since our argument is translation-invariant, we can assume $x_0 = 0$ to simplify our explanation. By selecting subsequences, we can find a map $u\in W^{1,2}(\Omega;\overline M)$ such that $u_k\to u$ weakly in $W^{1,2}(\Omega;\R^m)$ and strongly in $L^2(\Omega;\R^m)$.

To prove that $u$ is a minimizing constraint map for $\cE_q$, with $q = \lim_{k\to\infty}q_k$, 
let us fix a competitor $v\in W^{1,2}(\Omega;\overline M)$ such that $\supp (u - v)\subset B_r$. Then by Fubini's Theorem as well as the fact that $Du_k \to Dv$ weakly in $L^2(\Omega\setminus B_r;\R^{mn})$ and $u_k\to v$ strongly in $L^2(\Omega\setminus B_r;\R^m)$, we can find $s > r$ with $B_s\subset \Omega$ and $K > 1$ such that 
\begin{equation}\label{eq:K}
 \int_{\partial B_{{s}}} \left(|D u_k|^2+|Dv|^2 + \frac{|u_k - v|^2}{\e_k^2}\right) d\sigma \leq K^2,
\end{equation}
for all $k$ large, where $\e_k>0$ is chosen so that 
\begin{equation}\label{eq:ek}
\e_k^2 = \int_{\partial B_{{s}}} |u_k - v|^2\,d\sigma  \to 0.
\end{equation}
Let us choose $\alpha \in(0,1)$ small as in \eqref{eq:alpha}, and set 
\begin{equation}\label{eq:lk}
\lambda_k := \e_k^\alpha.
\end{equation}
Applying \cite[Lemma 1]{L} to the maps $u_k$ and $v$ (with $p =2$, $\lambda = \lambda_k$, $\e = \e_k$ and $\beta = 3/4$ there), we obtain from \eqref{eq:K} and \eqref{eq:ek} an extension $\vp_k \in W^{1,2}(B_{{s}};\R^m)$ such that  
\begin{equation}\label{eq:vpk}
\begin{aligned}
&\vp_k (x) = \begin{cases}
v((1-\lambda_k)^{-1}x),& \text{if } |x|\leq (1-\lambda_k){s},\\
 u_k(x),& \text{if } |x|={s},
 \end{cases}
 \end{aligned}
\end{equation}
which also satisfies 
\begin{equation}\label{eq:Dvpk-L2}
\int_{B_{{s}} \setminus B_{(1-\lambda_k){s}}} |D\vp_k|^2\,dx \leq cK^2\left( 1 + \frac{\e_k^2}{\lambda_k^2}\right)\lambda_k \leq 2cK^2\lambda_k,
\end{equation}
and 
\begin{equation}\label{eq:dist-vpk}
\dist(\vp_k,\overline M) \leq cK\e_k^{\frac{1}{4}}\lambda_k^{1-\frac{n}{2}} \leq cK\e_k^{\alpha (1-\frac{n}{2}) + \frac{1}{4}}\quad\text{in }B_{{s}},
\end{equation}
where the last inequalities in \eqref{eq:Dvpk-L2} and \eqref{eq:dist-vpk} are due to \eqref{eq:alpha} and \eqref{eq:lk}.

Let $\cN(\partial M)$ be the tubular neighborhood where the nearest point projection $\Pi\colon\cN(\partial M)\to \partial M$ is well-defined and of class $C^1$ (recall that we assume $\partial M$ to be of class $C^2$). Denoting $\cN(M)$ by the union of $\cN(\partial M)$ and $M$ (which again becomes a tubular neighborhood of $M$) we may define a Lipschitz projection $\pi :\cN(M)\to \overline M$ such that $\pi(w) = w$ if $w\in\overline M$ and $\pi(w) = \Pi(w)$ if $w\in \cN(\partial M)\setminus \overline M$. Clearly, we have 
\begin{equation}\label{eq:pi-Lip}
\| \nabla \pi \|_{L^\infty(\cN(M))}\leq c. 
\end{equation}

Now we are ready to prove $\cE_q[u] \leq \cE_q[v]$ over $B_s$. First, we observe from the lower semicontinuity of $\cE_q$ and $q_k\to q$ that 
\begin{equation}\label{eq:lsc}
\begin{aligned}
\int_{B_r} [|Du|^2 + q^2\chi_{u^{-1}(M)}]\,dx 
&\leq  \liminf_{k\to\infty} \int_{B_r} [|Du_k|^2 + q_k^2\chi_{u_k^{-1}(M)}]\,dx.
\end{aligned}
\end{equation}
Note that due to \eqref{eq:dist-vpk} and \eqref{eq:ek}, we obtain a well-defined map 
$$
v_k := \pi\circ\vp_k \in W^{1,2}(B_s;\overline M). 
$$
Since $\vp_k - u_k \in W_0^{1,2}(B_{{s}};\R^m)$, and $\pi\circ\vp_k = \pi\circ u_k = u_k$ on $\partial B_{s}$ (in the trace sense) due to \eqref{eq:vpk}, we have $v_k - u_k \in W_0^{1,2}(B_{{s}};\R^m)$. Therefore, $v_k$ is an admissible map, and the minimality of $u_k$ for $\cE_{q_k}$, yields
\begin{equation}\label{eq:min}
\int_{B_{{s}}} [|Du_k|^2 + q_k^2\chi_{u_k^{-1}(M)}]\,dx \leq \int_{B_{{s}}} [|Dv_k|^2 + q_k^2\chi_{v_k^{-1}(M)}]\,dx.
\end{equation}
In view of \eqref{eq:vpk}, \eqref{eq:lk} and \eqref{eq:ek}, we can compute that 
\begin{equation}\label{eq:vpk-re}
\begin{split}
\lim_{k\to\infty}\int_{B_{(1-\lambda_k){s}}} [|Dv_k|^2 + q_k^2\chi_{v_k^{-1}(M)}] \,dx & = \lim_{k\to\infty}(1 - \lambda_k)^n\int_{B_{{s}}}\left[\frac{|Dv|^2}{(1-\lambda_k)^2} + q_k^2 \chi_{v^{-1}(M)} \right]dx \\
&= \int_{B_{{s}}} [|Dv|^2 + q^2\chi_{v^{-1}(M)}]\,dx,
\end{split}
\end{equation}
where in the derivation of the last equality we also used $q_k\to q$. On the other hand, by \eqref{eq:Dvpk-L2} and \eqref{eq:pi-Lip},
\begin{equation}\label{eq:vpk-re2}
\begin{split}
  &\int_{B_{{s}}\setminus B_{(1-\lambda_k){s}}} [|Dv_k|^2 + q_k^2\chi_{v_k^{-1}(M)}] \,dx \\
  & \leq \int_{B_{{s}}\setminus B_{(1-\lambda_k){s}}} (|\nabla \pi_{v_k}|^2|D\vp_k|^2 + q_k^2)\,dx \leq c(K^2+q^2)\lambda_k.
  \end{split}
\end{equation}
Combining \eqref{eq:min}, \eqref{eq:lsc}, \eqref{eq:vpk-re} and \eqref{eq:vpk-re2} altogether, we arrive at 
\begin{equation}\label{eq:min-re}
\int_{B_{{s}}} [|Du|^2 + q^2\chi_{u^{-1}(M)}]\,dx \leq \int_{B_{{s}}} [|Dv|^2 + q^2\chi_{v^{-1}(M)}]\,dx. 
\end{equation}

Recall that $v\in W^{1,2}(\Omega;\overline M)$ was an arbitrary competitor, and we had $\supp(u- v) \subset B_{r}\Subset\Omega$. Additionally, $B_r\Subset B_s\subset\Omega$ (where $s$ may depend on $v$). Therefore, based on equation \eqref{eq:min-re}, we can conclude that $u$ is a minimizing constraint map of $\cE_q$ in $B_r$. Since we can replace $B_r$ with any ball $B_r(x_0)$ compactly contained in $\Omega$, we can deduce that $u$ is a minimizing constraint map of $\cE_q$ in any subdomain $\Omega'\Subset\Omega$.

To prove the strong convergence, we choose $v = u$ in \eqref{eq:vpk} and run the same argument above (more specifically, \eqref{eq:min} -- \eqref{eq:vpk-re2}), which yields 
$$
\lim_{k\to\infty} \int_{B_r} [|Du_k|^2 + q_k^2\chi_{u_k^{-1}(M)}]\,dx = \int_{B_r} [|Du|^2 + q^2\chi_{u^{-1}(M)}]\,dx.
$$
However, due to the weak lower semicontinuity of each term in the integrand, specifically $\int_{B_r} |Du_k|^2 ,dx$ and $|u_k^{-1}(M)\cap B_r|$, we can conclude that
$$
\lim_{k\to\infty} \int_{B_r} |Du_k|^2\,dx = \int_{B_r} |Du|^2\,dx.
$$
Then the strong convergence of $Du_k\to Du$ in $L^2(B_r;\R^{mn})$ follows easily. Again replacing $B_r$ with any $B_r(x_0)\subset \Omega$, we arrive at $Du_k\to Du$ in $L_\loc^2(\Omega;\R^{mn})$ as desired.
\end{proof}


\section{Removable singularities in special cases}\label{sec:remove}

We observe that singular points do not arise in the case of minimizing constraint maps for $\cE_q$, where $q > 0$, when the boundary of the constraint is one-sided, such as when 
$\partial M$ is a graph in $\R^m $; precisely, we assume that condition \eqref{eq:one-sided} below holds. It is worth noting that the case with $q = 0$ (although without an {\it a priori} estimate) was addressed by Fuchs \cite{F2} a long time ago.

\begin{prop}

Let $M$ be a $C^3$-domain with uniformly bounded principal curvatures, satisfying 
\begin{equation}
    \label{eq:one-sided}
    (a-\xi)\cdot\nu_a\leq 0\quad \text{ for all } a\in\partial M \text{ for some } \xi\in\R^m.
\end{equation} Suppose $u\in W^{1,2}(\Omega;\overline M)$ is a minimizing constraint map for $\cE_q$, where $q \in [0,q_0]$, and $| u |_{W^{1,2}(\Omega)}\leq \Lambda$. Then $u\in C_{loc}^{0,1}(\Omega)$ and for every subdomain $\Omega'\Subset\Omega$,
$$
[u]_{C^{0,1}(\Omega')} \leq c,
$$
where $c$ depends only on $n$, $m$, $\Lambda$, $\dist(\Omega',\partial\Omega)$, and $\partial M$.
\end{prop}

\begin{proof}

As established in Lemma \ref{lem:reg}, and in accordance with Theorem \ref{thm:e-reg}, let's demonstrate that the rescaled energy uniformly decays in the interior. For simplicity, we will consider $\Omega = B_4$ and $\Omega' = B_1$.

Let $\Lambda$ be a given constant, and let $\bar\e$ represent the small energy constant from Theorem \ref{thm:e-reg}. To arrive at a contradiction, we can construct a sequence analogous to that in Lemma \ref{lem:reg}. In other words, for $u_k\in W^{1,2}(B_4;\overline M)$, we assume that it is a minimizing constraint map for $\cE_{q_k}$, where $q_k\in [0,q_0]$, $x_k\in B_1$, and $r_k\to 0$ such that $| u_k |_{W^{1,2}(B_4)} \leq \Lambda$, but $E(u_k,x_k,2r_k) \geq \bar\e^2$. As we let $k\to\infty$, we can employ Lemma \ref{lem:strong} to extract a subsequence (denoted with the same subscript) such that $q_k\to q\in [0,q_0]$, $x_k\to x_0\in \overline{B_1}$, and $u_k \to u$ strongly in $W^{1,2}(B_2;\R^m)$. This results in a minimizing constraint map $u\in W^{1,2}(B_2;\overline M)$ for $\cE_q$, and it satisfies the condition:

$$
\Theta_u(x_0) \geq \bar\e^2 > 0. 
$$

Select any tangent map $\phi\in T_{x_0}u$. According to Corollary \ref{cor:blowup}, $\phi \in W_\loc^{1,2}(\R^n;\overline M)$ is a $0$-homogeneous minimizing constraint map of $\cE_0$. Since $\phi$ is energy-minimizing, the one-sided condition on $\partial M$ allows us to invoke \cite[Theorem 1]{F2}, which demonstrates that $\Sigma(\phi) = \emptyset$, meaning that $\phi$ is continuous everywhere. Because $\phi$ is $0$-homogeneous, its everywhere continuity necessitates that it must be constant almost everywhere. In particular, this implies that $\Theta_u(x_0) = \int_{B_1} |D\phi|^2,dy = 0$, which leads to a contradiction.

Hence, we conclude that if $u\in W^{1,2}(B_4;\overline M)$ is a minimizing constraint map for $\cE_q$, with $q\in[0,q_0]$, and $| u |_{W^{1,2}(B_4)} \leq \Lambda$, then there exists a small constant $r \in (0,1)$, depending only on $n$, $m$, $q_0$, $\Lambda$, and $\partial M$, such that $E(u,x_0,2r) < \bar\e^2$ for all $x_0\in B_1$. By applying Theorem \ref{thm:e-reg}, we obtain $u\in C^{0,1}(B_r(x_0))$ with

$$
[u]_{C^{0,1}(B_r(x_0))} \leq \frac{c_0}{r},
$$
where $c_0$ depends only on $n$, $m$, and $\partial M$. Since $x_0$ is an arbitrary point in $B_1$ and $r$ is independent of $x_0$, we deduce that $u\in C^{0,1}(B_1)$, and 
$$
[u]_{C^{0,1}(B_1)} \leq c,
$$
where $c$ depends only on $n$, $m$, $q_0$, $\Lambda$, and $\partial M$. This proves the proposition for the case when $\Omega = B_4$ and $\Omega'=B_1$. Using the standard scaling and covering argument, we can extend this result to the general case.

\end{proof}


 \section*{Acknowledgements}
A.F. has been partially supported by the European Research Council under the Grant Agreement No.~721675 ``Regularity and Stability in Partial Differential Equations (RSPDE)''.
A.G. was supported by Dr.~Max Rössler, the Walter Haefner Foundation, and the ETH Zürich Foundation, and he thanks Uppsala University and KTH for their hospitality during a visit where part of this research was conducted.
H.S. was supported by the Swedish Research Council (grant no.~2021-03700).

\section*{Declarations}

\noindent {\bf  Data availability statement:} All necessary data are included in the manuscript.

\medskip
\noindent {\bf  Funding and/or Conflicts of interests/Competing interests:} 
 The authors declare that there are no financial, competing, or conflicts of interests.


\end{document}